\documentclass[a4paper,12pt]{article}
\usepackage{amsmath,amssymb,amsthm}
\usepackage{amscd}
\usepackage{color}
\usepackage{graphicx}

\makeatletter

\@addtoreset{equation}{subsection}
\makeatother

\newtheorem{dfn}{Definition}[subsection]
\newtheorem{prop}[dfn]{Proposition}
\newtheorem{thm}[dfn]{Theorem}
\newtheorem{lem}[dfn]{Lemma}
\newtheorem{cor}[dfn]{Corollary}

\newtheorem{rem}[dfn]{Remark}
\newtheorem{conj}[dfn]{Conjecture}
\newtheorem{quest}[dfn]{Question}

\newcommand\gfg{{\mathfrak{g}}}

\newcommand\Ker{{\rm Ker}}

\newcommand\ISE{{\mathcal{I}(S,E)}}

\begin{document}

\title{Intersections of curves on surfaces
and their applications to mapping class groups}

\author{Nariya Kawazumi and Yusuke Kuno}
\date{}

\maketitle

\begin{abstract}
We introduce an operation that measures the self intersections of paths on a surface.
As applications, we give a criterion of the realizability of a generalized Dehn twist,
and derive a geometric constraint on the image of the Johnson homomorphism.
\end{abstract}

\section{Introduction}

In study of the mapping class group of a surface, it is sometimes
convenient and crucial to work with curves, i.e., loops or paths, on the surface.
Since the mapping class group is the group of isotopy classes of diffeomorphisms
of the surface, it acts on the homotopy classes of curves. As is
illustrated in the classical theorem of Dehn-Nielsen, this action
distinguishes elements of the mapping class group well. Moreover, this
action preserves intersections of curves.

In this paper we introduce an operation that measures the self intersections
of paths on a surface, and discuss its applications to the mapping class groups.
Let $S$ be an oriented surface, and $*_0,*_1\in \partial S$ points
on the boundary. We denote by $\Pi S(*_0,*_1)$ the set of homotopy classes
of paths from $*_0$ to $*_1$, and by $\hat{\pi}^{\prime}(S)$ the set
of homotopy classes of non-trivial free loops on $S$.
In \S \ref{subsec:BM} we introduce a $\mathbb{Q}$-linear map
$$\mu\colon \mathbb{Q}\Pi S(*_0,*_1) \to \mathbb{Q}\Pi S(*_0,*_1)
\otimes \mathbb{Q}\hat{\pi}^{\prime}(S),$$
by looking at the self intersections of a given path. This map is inspired
by Turaev's self intersection \cite{T1}, and is actually a refinement of it.

One motivation to introduce $\mu$ comes from the Goldman-Turaev Lie bialgebra.
The free $\mathbb{Q}$-vector space $\mathbb{Q}\hat{\pi}^{\prime}(S)$
spanned by the set $\hat{\pi}^{\prime}(S)$
is an involutive Lie bialgebra with respect to the Goldman bracket \cite{Go} and the Turaev cobracket \cite{T2}.
In \cite{KK1} \cite{KKp}, we showed that the $\mathbb{Q}$-vector space $\mathbb{Q}\Pi S(*_0,*_1)$ is a (left)
$\mathbb{Q}\hat{\pi}^{\prime}(S)$-module with respect to a structure map
$\sigma\colon \mathbb{Q}\hat{\pi}^{\prime}(S) \otimes \mathbb{Q}\Pi S(*_0,*_1)
\to \mathbb{Q}\Pi S(*_0,*_1)$. In fact, by investigating the properties of $\mu$
we arrive at the notion of a {\it comodule of a Lie coalgebra}, and that of a
{\it bimodule of a Lie bialgebra} (see Appendix). We show that $\mathbb{Q}\Pi S(*_0,*_1)$
is a $\mathbb{Q}\hat{\pi}^{\prime}(S)$-bimodule with respect to $\sigma$ and $\mu$.

In \cite{KKp}, we introduced a filtration on $\mathbb{Q}\hat{\pi}^{\prime}(S)$
and showed that the Goldman bracket induces a Lie bracket on the completion
$\widehat{\mathbb{Q}\hat{\pi}}(S)$, which we called the {\it completed Goldman Lie algebra}.
We also defined a completion of $\mathbb{Q}\Pi S(*_0,*_1)$ and showed
that it is a complete $\widehat{\mathbb{Q}\hat{\pi}}(S)$-module.
As we will see, the operation $\mu$ recovers the Turaev cobracket. Analyzing the
behavior of $\mu$ under the conjunction of paths, we show that $\mu$ naturally
extends to completions and the Turaev cobracket extends to a complete cobracket
on $\widehat{\mathbb{Q}\hat{\pi}}(S)$. Thus we could call
$\widehat{\mathbb{Q}\hat{\pi}}(S)$ the {\it completed Goldman-Turaev Lie bialgebra}.
Also, we show that the completion $\widehat{\mathbb{Q}\Pi S}(*_0,*_1)$ is
a complete $\widehat{\mathbb{Q}\hat{\pi}}(S)$-bimodule.

Since $\mu$ is defined in terms of the intersections of curves, it is automatically
compatible with the action of the mapping class group. In \S \ref{sec:App} and \S \ref{sec:App2} we give
two applications of this fact to study of the mapping class group. The first one is
an application to generalized Dehn twists \cite{KKp} \cite{Ku2} \cite{MT},
which are elements of a certain enlargement
of the mapping class group. We can ask whether a generalized Dehn twist is realized by
a diffeomorphism of the surface, and we give a criterion of the realizability of a generalized
Dehn twist using $\mu$ (Proposition \ref{criterion}). This criterion is powerful enough so that we can extend results
about a figure eight \cite{KKp} \cite{Ku2} to loops in wider classes (Theorem \ref{main.thm}).
The second one is an application to the Johnson homomorphism, which is an embedding
of the `smallest' Torelli group (in the sense of Putman \cite{P}) into a pro-nilpotent group.
Using the fact that a diffeomorphism preserves $\mu$, we derive a geometric
constraint on the image of the Johnson homomorphism (Theorem \ref{zero}). That this constraint is
non-trivial can be seen from examples of null-homologous, non-simple loops
whose generalized Dehn twists are not realized by diffeomorphisms.
Further, based on a tensorial description of the homotopy intersection
form given by Massuyeau and Turaev \cite{MT}, we study some lower terms of a tensorial description of
the constraint (Theorem \ref{laurent}). In particular we show that the Morita trace \cite{MoAQ}
is recovered from the lowest term of our geometric constraint (Theorem \ref{54trace}).

\medskip
\noindent \textbf{Acknowledgments.}
The first-named author is partially supported by the Grant-in-Aid for
Scientific Research (A) (No.20244003), (A) (No.22244005) and (B) (No.24340010) from the
Japan Society for Promotion of Sciences. The second-named
author is supported by JSPS Research Fellowships
for Young Scientists (22$\cdot$4810) and the Grant-in-Aid for
Research Activity Start-up (No.24840038). They would like to thank
Naoya Enomoto for informing them his recent result about the Enomoto-Satoh traces,
Robert Penner for fruitful discussions at QGM, Aarhus University,
and Vladimir Turaev for useful comments on a draft of this paper.
A part of this work has been done
during the second author's stay at QGM from September 2011 to November 2011.
He would like to thank QGM for kind hospitality.

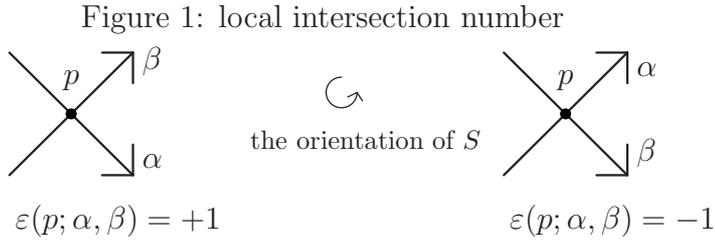
\tableofcontents

\section{The Goldman-Turaev Lie bialgebra and its bimodule}
\label{sec:GT}

Let $S$ be a connected oriented surface.
We denote by $\hat{\pi}(S)=[S^1,S]$ the homotopy
set of oriented free loops on $S$. In other words, $\hat{\pi}(S)$
is the set of conjugacy classes of $\pi_1(S)$. We denote
by $|\ |\colon \pi_1(S)\to \hat{\pi}(S)$ the natural projection,
and we also denote by
$|\ |\colon \mathbb{Q}\pi_1(S)\to \mathbb{Q}\hat{\pi}(S)$
its $\mathbb{Q}$-linear extension.

\subsection{The Goldman-Turaev Lie bialgebra}

We recall the Goldman-Turaev Lie bialgebra \cite{Go} \cite{T2}.
\begin{figure}
\begin{center}
\caption{local intersection number}

\vspace{0.2cm}

\unitlength 0.1in
\begin{picture}( 32.4800,  8.4000)(  2.0000,-10.3000)
%
{\color[named]{Black}{%
\special{pn 13}%
\special{pa 200 840}%
\special{pa 840 200}%
\special{fp}%
\special{pa 200 200}%
\special{pa 840 840}%
\special{fp}%
}}%
%
{\color[named]{Black}{%
\special{pn 4}%
\special{sh 1}%
\special{ar 520 520 26 26 0  6.28318530717959E+0000}%
}}%
%
{\color[named]{Black}{%
\special{pn 13}%
\special{pa 840 840}%
\special{pa 840 680}%
\special{fp}%
\special{pa 840 840}%
\special{pa 680 840}%
\special{fp}%
}}%
%
{\color[named]{Black}{%
\special{pn 13}%
\special{pa 840 200}%
\special{pa 680 200}%
\special{fp}%
\special{pa 840 200}%
\special{pa 840 360}%
\special{fp}%
}}%
\put(4.8000,-4.0000){\makebox(0,0)[lb]{$p$}}%
\put(8.8800,-8.0000){\makebox(0,0)[lb]{$\alpha$}}%
\put(8.8800,-3.2000){\makebox(0,0)[lb]{$\beta$}}%
\put(2.3200,-11.6000){\makebox(0,0)[lb]{$\varepsilon(p;\alpha,\beta)=+1$}}%
%
{\color[named]{Black}{%
\special{pn 13}%
\special{pa 2760 840}%
\special{pa 3400 200}%
\special{fp}%
\special{pa 2760 200}%
\special{pa 3400 840}%
\special{fp}%
}}%
%
{\color[named]{Black}{%
\special{pn 4}%
\special{sh 1}%
\special{ar 3080 520 26 26 0  6.28318530717959E+0000}%
}}%
%
{\color[named]{Black}{%
\special{pn 13}%
\special{pa 3400 840}%
\special{pa 3400 680}%
\special{fp}%
\special{pa 3400 840}%
\special{pa 3240 840}%
\special{fp}%
}}%
%
{\color[named]{Black}{%
\special{pn 13}%
\special{pa 3400 200}%
\special{pa 3240 200}%
\special{fp}%
\special{pa 3400 200}%
\special{pa 3400 360}%
\special{fp}%
}}%
\put(30.4000,-4.0000){\makebox(0,0)[lb]{$p$}}%
\put(34.4800,-8.0000){\makebox(0,0)[lb]{$\beta$}}%
\put(34.4800,-3.2000){\makebox(0,0)[lb]{$\alpha$}}%
\put(27.9200,-11.6000){\makebox(0,0)[lb]{$\varepsilon(p;\alpha,\beta)=-1$}}%
%
{\color[named]{Black}{%
\special{pn 8}%
\special{ar 1920 408 80 80  6.2831853 6.2831853}%
\special{ar 1920 408 80 80  0.0000000 4.7123890}%
}}%
%
{\color[named]{Black}{%
\special{pn 8}%
\special{pa 2000 408}%
\special{pa 1952 432}%
\special{fp}%
\special{pa 2000 408}%
\special{pa 2028 452}%
\special{fp}%
}}%
\put(14.4800,-7.0400){\makebox(0,0)[lb]{{\footnotesize the orientation of $S$}}}%
\end{picture}%
\end{center}
\end{figure}

Let $\alpha$ and $\beta$ be oriented immersed loops on $S$ such that
their intersections consist of transverse double points. For each
$p \in \alpha \cap \beta$, let $\alpha_p\beta_p\in \pi_1(S,p)$ be the
loop going first along the loop $\alpha$ based at $p$, then going
along $\beta$ based at $p$. Also, let $\varepsilon(p;\alpha,\beta)\in \{ \pm 1 \}$ be the
local intersection number of $\alpha$ and $\beta$ at $p$.
See Figure 1. The Goldman bracket of $\alpha$ and $\beta$ is defined as
\begin{equation}
\label{Gbra}
[\alpha,\beta]:=\sum_{p\in \alpha \cap \beta}
\varepsilon(p; \alpha,\beta)|\alpha_p\beta_p| \in \mathbb{Q}\hat{\pi}(S).
\end{equation}
The free $\mathbb{Q}$-vector space $\mathbb{Q}\hat{\pi}(S)$ spanned by the set $\hat{\pi}(S)$
equipped with this bracket is a Lie algebra. See \cite{Go}.
Let $1\in \hat{\pi}(S)$ be the class of a constant loop, then its
linear span $\mathbb{Q}1$ is an ideal of $\mathbb{Q}\hat{\pi}(S)$.
We denote by $\mathbb{Q}\hat{\pi}^{\prime}(S)$
the quotient Lie algebra $\mathbb{Q}\hat{\pi}(S)/\mathbb{Q}1$, and let
$\varpi \colon \mathbb{Q}\hat{\pi}(S) \to \mathbb{Q}\hat{\pi}^{\prime}(S)$ be the projection.
We write $|\ |^{\prime}:=\varpi \circ |\ |\colon \mathbb{Q}\pi_1(S)
\to \mathbb{Q}\hat{\pi}^{\prime}(S)$.

Let $\alpha\colon S^1 \to S$ be an oriented immersed loop such that its self intersections
consist of transverse double points. Set
$D=D_{\alpha}:=\{ (t_1,t_2) \in S^1\times S^1; t_1\neq t_2, \alpha(t_1)=\alpha(t_2) \}$.
For $(t_1,t_2)\in D$, let $\alpha_{t_1t_2}$ (resp. $\alpha_{t_2t_1}$) be
the restriction of $\alpha$ to the interval $[t_1,t_2]$ (resp. $[t_2,t_1]$) $\subset S^1$
(they are indeed loops since $\alpha(t_1)=\alpha(t_2)$).
Also, let $\dot{\alpha}(t_i) \in T_{\alpha(t_i)}S$ be the velocity vectors of
$\alpha$ at $t_i$, and set $\varepsilon(\dot{\alpha}(t_1),\dot{\alpha}(t_2))=+1$ if
$(\dot{\alpha}(t_1),\dot{\alpha}(t_2))$ gives the orientation of $S$, and
$\varepsilon(\dot{\alpha}(t_1),\dot{\alpha}(t_2))=-1$ otherwise.
The Turaev cobracket of $\alpha$ is defined as
\begin{equation}
\label{Tcob}
\delta(\alpha):=\sum_{(t_1,t_2)\in D}
\varepsilon(\dot{\alpha}(t_1),\dot{\alpha}(t_2)) |\alpha_{t_1t_2}|^{\prime} \otimes
|\alpha_{t_2t_1}|^{\prime}\in \mathbb{Q}\hat{\pi}^{\prime}(S)
\otimes \mathbb{Q}\hat{\pi}^{\prime}(S).
\end{equation}
This gives rise to a well-defined Lie cobracket
$\delta\colon \mathbb{Q}\hat{\pi}^{\prime}(S) \to
\mathbb{Q}\hat{\pi}^{\prime}(S) \otimes \mathbb{Q}\hat{\pi}^{\prime}(S)$
(note that $\delta(1)=0$). Moreover, $\mathbb{Q}\hat{\pi}^{\prime}(S)$
is an involutive Lie bialgebra with respect to the Goldman bracket and the
Turaev cobracket. See \cite{T2}. The involutivity is due to Chas \cite{Chas}.
We call $\mathbb{Q}\hat{\pi}^{\prime}(S)$ the {\it Goldman-Turaev Lie bialgebra}.

\subsection{A $\mathbb{Q}\hat{\pi}^{\prime}(S)$-bimodule}
\label{subsec:BM}

Hereafter we assume that the boundary of $S$ is not empty.
Take distinct points $*_0, *_1 \in \partial S$,
and let $\Pi S(*_0,*_1)$ be the homotopy set $[([0,1],0,1),(S,*_0,*_1)]$.
We shall show that the free $\mathbb{Q}$-vector space $\mathbb{Q}\Pi S(*_0,*_1)$ spanned
by the set $\Pi S(*_0,*_1)$ has
the structure of an involutive right $\mathbb{Q}\hat{\pi}^{\prime}(S)$-bimodule.
For the definition of a bimodule, see Appendix.
In \S \ref{*_1=*_2} we discuss the case $*_0=*_1$.

\textit{A left $\mathbb{Q}\hat{\pi}^{\prime}(S)$-module structure}.
Let $\alpha$ be an oriented immersed loop on $S$,
and $\beta\colon [0,1]\to S$ an immersed path from $*_0$ to $*_1$ such that
their intersections consist of transverse double points. Then the formula
\begin{equation}
\label{sigma}
\sigma(\alpha \otimes \beta):=\sum_{p\in \alpha \cap \beta}
\varepsilon(p;\alpha,\beta)\beta_{*_0 p} \alpha_p \beta_{p*_1} \in
\mathbb{Q}\Pi S(*_0,*_1)
\end{equation}
gives rise to a well-defined $\mathbb{Q}$-linear map $\sigma\colon
\mathbb{Q}\hat{\pi}^{\prime}(S) \otimes \mathbb{Q}\Pi S(*_0,*_1)
\to \mathbb{Q}\Pi S(*_0,*_1)$. Here, $\varepsilon(p;\alpha,\beta)\in \{ \pm 1 \}$
has the same meaning as before, and $\beta_{*_0p} \alpha_p \beta_{p*_1}$
means the path going first from $*_0$ to $p$ along $\beta$, then
going along $\alpha$ based at $p$, and finally going from $p$ to $*_1$
along $\beta$. By the same proof as that of \cite{KK1} Proposition 3.2.2,
we see that $\mathbb{Q}\Pi S(*_0,*_1)$ is a left
$\mathbb{Q}\hat{\pi}^{\prime}(S)$-module with respect to $\sigma$.
See also \cite{KKp} \S 4.

\textit{A right $\mathbb{Q}\hat{\pi}^{\prime}(S)$-comodule structure}.
Let $\gamma\colon [0,1]\to S$ be an immersed path from $*_0$ to $*_1$
such that its self intersections consist of transverse double points.
Let $\Gamma=\Gamma_{\gamma} \subset {\rm Int}(S)$ be the set of double points of $\gamma$.
For $p\in \Gamma$, we denote $\gamma^{-1}(p)=\{t_1^p,t_2^p \}$, so that
$t_1^p< t_2^p$. Set
\begin{equation}
\label{dfn-mu}
\mu(\gamma):=-\sum_{p\in \Gamma}
\varepsilon(\dot{\gamma}(t_1^p),\dot{\gamma}(t_2^p))
( \gamma_{0 t_1^p} \gamma_{t_2^p 1} )
\otimes |\gamma_{t_1^p t_2^p}|^{\prime} \in \mathbb{Q}\Pi S(*_0,*_1) \otimes
\mathbb{Q}\hat{\pi}^{\prime}(S).
\end{equation}
Here $\varepsilon(\dot{\gamma}(t_1^p),\dot{\gamma}(t_2^p))\in \{ \pm 1 \}$
has the same meaning as before, $\gamma_{0 t_1^p} \gamma_{t_2^p 1}$
is the conjunction of the restrictions of $\gamma$ to $[0,t_1^p]$ and
$[t_2^p,1]$, and $\gamma_{t_1^p t_2^p}$ is the restriction of $\gamma$
to $[t_1^p,t_2^p]$. The map $\mu$ is closely related to Turaev's self intersection
\cite{T1}. See Remark \ref{rem:Tu}.

\begin{prop}
\label{comod}
The formula {\rm (\ref{dfn-mu})} gives rise to a well-defined $\mathbb{Q}$-linear map
$$\mu\colon \mathbb{Q}\Pi S(*_0,*_1) \to
\mathbb{Q}\Pi S(*_0,*_1) \otimes \mathbb{Q}\hat{\pi}^{\prime}(S).$$
Moreover, $\mathbb{Q}\Pi S(*_0,*_1)$ is a right
$\mathbb{Q}\hat{\pi}^{\prime}(S)$-comodule
with respect to $\mu$.
\end{prop}

\begin{proof}
Any immersions $\gamma$ and $\gamma^{\prime}$
with $\gamma(0)=\gamma^{\prime}(0)=*_0$ and $\gamma(1)=\gamma^{\prime}(1)=*_1$,
homotpic to each other relative to $\{ 0,1 \}$, such that their self intersections consist
of transverse double points, are related by a sequence of
three local moves ($\omega 1$), ($\omega 2$), ($\omega 3$),
and an ambient isotopy of $S$. See Goldman \cite{Go} \S 5 and Figure 2.
To prove that $\mu$ is well-defined, it is sufficient to verify
that $\mu(\gamma)=\mu(\gamma^{\prime})$ if $\gamma$ and $\gamma^{\prime}$
are related by one of the three moves.

\begin{figure}
\begin{center}
\caption{local moves ($\omega 1$), ($\omega 2$), and ($\omega 3$)}

\vspace{0.2cm}

\unitlength 0.1in
\begin{picture}( 36.8000, 12.5000)(  2.0000,-12.9400)
%
{\color[named]{Black}{%
\special{pn 8}%
\special{ar 1320 814 480 480  0.0000000 0.1250000}%
\special{ar 1320 814 480 480  0.2000000 0.3250000}%
\special{ar 1320 814 480 480  0.4000000 0.5250000}%
\special{ar 1320 814 480 480  0.6000000 0.7250000}%
\special{ar 1320 814 480 480  0.8000000 0.9250000}%
\special{ar 1320 814 480 480  1.0000000 1.1250000}%
\special{ar 1320 814 480 480  1.2000000 1.3250000}%
\special{ar 1320 814 480 480  1.4000000 1.5250000}%
\special{ar 1320 814 480 480  1.6000000 1.7250000}%
\special{ar 1320 814 480 480  1.8000000 1.9250000}%
\special{ar 1320 814 480 480  2.0000000 2.1250000}%
\special{ar 1320 814 480 480  2.2000000 2.3250000}%
\special{ar 1320 814 480 480  2.4000000 2.5250000}%
\special{ar 1320 814 480 480  2.6000000 2.7250000}%
\special{ar 1320 814 480 480  2.8000000 2.9250000}%
\special{ar 1320 814 480 480  3.0000000 3.1250000}%
\special{ar 1320 814 480 480  3.2000000 3.3250000}%
\special{ar 1320 814 480 480  3.4000000 3.5250000}%
\special{ar 1320 814 480 480  3.6000000 3.7250000}%
\special{ar 1320 814 480 480  3.8000000 3.9250000}%
\special{ar 1320 814 480 480  4.0000000 4.1250000}%
\special{ar 1320 814 480 480  4.2000000 4.3250000}%
\special{ar 1320 814 480 480  4.4000000 4.5250000}%
\special{ar 1320 814 480 480  4.6000000 4.7250000}%
\special{ar 1320 814 480 480  4.8000000 4.9250000}%
\special{ar 1320 814 480 480  5.0000000 5.1250000}%
\special{ar 1320 814 480 480  5.2000000 5.3250000}%
\special{ar 1320 814 480 480  5.4000000 5.5250000}%
\special{ar 1320 814 480 480  5.6000000 5.7250000}%
\special{ar 1320 814 480 480  5.8000000 5.9250000}%
\special{ar 1320 814 480 480  6.0000000 6.1250000}%
\special{ar 1320 814 480 480  6.2000000 6.2832853}%
}}%
%
{\color[named]{Black}{%
\special{pn 13}%
\special{ar 1320 814 480 160  6.2831853 6.2831853}%
\special{ar 1320 814 480 160  0.0000000 3.1415927}%
}}%
%
{\color[named]{Black}{%
\special{pn 8}%
\special{ar 3400 814 480 480  0.0000000 0.1250000}%
\special{ar 3400 814 480 480  0.2000000 0.3250000}%
\special{ar 3400 814 480 480  0.4000000 0.5250000}%
\special{ar 3400 814 480 480  0.6000000 0.7250000}%
\special{ar 3400 814 480 480  0.8000000 0.9250000}%
\special{ar 3400 814 480 480  1.0000000 1.1250000}%
\special{ar 3400 814 480 480  1.2000000 1.3250000}%
\special{ar 3400 814 480 480  1.4000000 1.5250000}%
\special{ar 3400 814 480 480  1.6000000 1.7250000}%
\special{ar 3400 814 480 480  1.8000000 1.9250000}%
\special{ar 3400 814 480 480  2.0000000 2.1250000}%
\special{ar 3400 814 480 480  2.2000000 2.3250000}%
\special{ar 3400 814 480 480  2.4000000 2.5250000}%
\special{ar 3400 814 480 480  2.6000000 2.7250000}%
\special{ar 3400 814 480 480  2.8000000 2.9250000}%
\special{ar 3400 814 480 480  3.0000000 3.1250000}%
\special{ar 3400 814 480 480  3.2000000 3.3250000}%
\special{ar 3400 814 480 480  3.4000000 3.5250000}%
\special{ar 3400 814 480 480  3.6000000 3.7250000}%
\special{ar 3400 814 480 480  3.8000000 3.9250000}%
\special{ar 3400 814 480 480  4.0000000 4.1250000}%
\special{ar 3400 814 480 480  4.2000000 4.3250000}%
\special{ar 3400 814 480 480  4.4000000 4.5250000}%
\special{ar 3400 814 480 480  4.6000000 4.7250000}%
\special{ar 3400 814 480 480  4.8000000 4.9250000}%
\special{ar 3400 814 480 480  5.0000000 5.1250000}%
\special{ar 3400 814 480 480  5.2000000 5.3250000}%
\special{ar 3400 814 480 480  5.4000000 5.5250000}%
\special{ar 3400 814 480 480  5.6000000 5.7250000}%
\special{ar 3400 814 480 480  5.8000000 5.9250000}%
\special{ar 3400 814 480 480  6.0000000 6.1250000}%
\special{ar 3400 814 480 480  6.2000000 6.2832853}%
}}%
%
{\color[named]{Black}{%
\special{pn 13}%
\special{ar 3240 814 320 160  6.2831853 6.2831853}%
\special{ar 3240 814 320 160  0.0000000 3.1415927}%
}}%
%
{\color[named]{Black}{%
\special{pn 13}%
\special{ar 3560 814 320 160  6.2831853 6.2831853}%
\special{ar 3560 814 320 160  0.0000000 3.1415927}%
}}%
%
{\color[named]{Black}{%
\special{pn 13}%
\special{ar 3400 814 160 160  3.1415927 6.2831853}%
}}%
\put(22.5600,-8.1400){\makebox(0,0)[lb]{$\longleftrightarrow$}}%
\put(2.0000,-1.7400){\makebox(0,0)[lb]{($\omega 1$) \textit{birth-death of monogons}}}%
%
{\color[named]{Black}{%
\special{pn 4}%
\special{sh 1}%
\special{ar 3400 958 26 26 0  6.28318530717959E+0000}%
}}%
\end{picture}%

\vspace{0.2cm}

\unitlength 0.1in
\begin{picture}( 36.8000, 12.5000)(  2.0000,-12.9400)
%
{\color[named]{Black}{%
\special{pn 8}%
\special{ar 1320 814 480 480  0.0000000 0.1250000}%
\special{ar 1320 814 480 480  0.2000000 0.3250000}%
\special{ar 1320 814 480 480  0.4000000 0.5250000}%
\special{ar 1320 814 480 480  0.6000000 0.7250000}%
\special{ar 1320 814 480 480  0.8000000 0.9250000}%
\special{ar 1320 814 480 480  1.0000000 1.1250000}%
\special{ar 1320 814 480 480  1.2000000 1.3250000}%
\special{ar 1320 814 480 480  1.4000000 1.5250000}%
\special{ar 1320 814 480 480  1.6000000 1.7250000}%
\special{ar 1320 814 480 480  1.8000000 1.9250000}%
\special{ar 1320 814 480 480  2.0000000 2.1250000}%
\special{ar 1320 814 480 480  2.2000000 2.3250000}%
\special{ar 1320 814 480 480  2.4000000 2.5250000}%
\special{ar 1320 814 480 480  2.6000000 2.7250000}%
\special{ar 1320 814 480 480  2.8000000 2.9250000}%
\special{ar 1320 814 480 480  3.0000000 3.1250000}%
\special{ar 1320 814 480 480  3.2000000 3.3250000}%
\special{ar 1320 814 480 480  3.4000000 3.5250000}%
\special{ar 1320 814 480 480  3.6000000 3.7250000}%
\special{ar 1320 814 480 480  3.8000000 3.9250000}%
\special{ar 1320 814 480 480  4.0000000 4.1250000}%
\special{ar 1320 814 480 480  4.2000000 4.3250000}%
\special{ar 1320 814 480 480  4.4000000 4.5250000}%
\special{ar 1320 814 480 480  4.6000000 4.7250000}%
\special{ar 1320 814 480 480  4.8000000 4.9250000}%
\special{ar 1320 814 480 480  5.0000000 5.1250000}%
\special{ar 1320 814 480 480  5.2000000 5.3250000}%
\special{ar 1320 814 480 480  5.4000000 5.5250000}%
\special{ar 1320 814 480 480  5.6000000 5.7250000}%
\special{ar 1320 814 480 480  5.8000000 5.9250000}%
\special{ar 1320 814 480 480  6.0000000 6.1250000}%
\special{ar 1320 814 480 480  6.2000000 6.2832853}%
}}%
%
{\color[named]{Black}{%
\special{pn 8}%
\special{ar 3400 814 480 480  0.0000000 0.1250000}%
\special{ar 3400 814 480 480  0.2000000 0.3250000}%
\special{ar 3400 814 480 480  0.4000000 0.5250000}%
\special{ar 3400 814 480 480  0.6000000 0.7250000}%
\special{ar 3400 814 480 480  0.8000000 0.9250000}%
\special{ar 3400 814 480 480  1.0000000 1.1250000}%
\special{ar 3400 814 480 480  1.2000000 1.3250000}%
\special{ar 3400 814 480 480  1.4000000 1.5250000}%
\special{ar 3400 814 480 480  1.6000000 1.7250000}%
\special{ar 3400 814 480 480  1.8000000 1.9250000}%
\special{ar 3400 814 480 480  2.0000000 2.1250000}%
\special{ar 3400 814 480 480  2.2000000 2.3250000}%
\special{ar 3400 814 480 480  2.4000000 2.5250000}%
\special{ar 3400 814 480 480  2.6000000 2.7250000}%
\special{ar 3400 814 480 480  2.8000000 2.9250000}%
\special{ar 3400 814 480 480  3.0000000 3.1250000}%
\special{ar 3400 814 480 480  3.2000000 3.3250000}%
\special{ar 3400 814 480 480  3.4000000 3.5250000}%
\special{ar 3400 814 480 480  3.6000000 3.7250000}%
\special{ar 3400 814 480 480  3.8000000 3.9250000}%
\special{ar 3400 814 480 480  4.0000000 4.1250000}%
\special{ar 3400 814 480 480  4.2000000 4.3250000}%
\special{ar 3400 814 480 480  4.4000000 4.5250000}%
\special{ar 3400 814 480 480  4.6000000 4.7250000}%
\special{ar 3400 814 480 480  4.8000000 4.9250000}%
\special{ar 3400 814 480 480  5.0000000 5.1250000}%
\special{ar 3400 814 480 480  5.2000000 5.3250000}%
\special{ar 3400 814 480 480  5.4000000 5.5250000}%
\special{ar 3400 814 480 480  5.6000000 5.7250000}%
\special{ar 3400 814 480 480  5.8000000 5.9250000}%
\special{ar 3400 814 480 480  6.0000000 6.1250000}%
\special{ar 3400 814 480 480  6.2000000 6.2832853}%
}}%
\put(22.5600,-8.1400){\makebox(0,0)[lb]{$\longleftrightarrow$}}%
\put(2.0000,-1.7400){\makebox(0,0)[lb]{($\omega 2$) \textit{birth-death of bigons}}}%
%
{\color[named]{Black}{%
\special{pn 13}%
\special{pa 872 974}%
\special{pa 1768 974}%
\special{fp}%
}}%
%
{\color[named]{Black}{%
\special{pn 13}%
\special{pa 872 654}%
\special{pa 1768 654}%
\special{fp}%
}}%
%
{\color[named]{Black}{%
\special{pn 13}%
\special{ar 3400 974 448 320  3.1415927 6.2831853}%
}}%
%
{\color[named]{Black}{%
\special{pn 13}%
\special{ar 3400 654 448 320  6.2831853 6.2831853}%
\special{ar 3400 654 448 320  0.0000000 3.1415927}%
}}%
%
{\color[named]{Black}{%
\special{pn 4}%
\special{sh 1}%
\special{ar 3008 814 26 26 0  6.28318530717959E+0000}%
}}%
%
{\color[named]{Black}{%
\special{pn 4}%
\special{sh 1}%
\special{ar 3784 814 26 26 0  6.28318530717959E+0000}%
}}%
\end{picture}%

\vspace{0.2cm}

\unitlength 0.1in
\begin{picture}( 36.8000, 12.5000)(  2.0000,-12.9400)
%
{\color[named]{Black}{%
\special{pn 8}%
\special{ar 1320 814 480 480  0.0000000 0.1250000}%
\special{ar 1320 814 480 480  0.2000000 0.3250000}%
\special{ar 1320 814 480 480  0.4000000 0.5250000}%
\special{ar 1320 814 480 480  0.6000000 0.7250000}%
\special{ar 1320 814 480 480  0.8000000 0.9250000}%
\special{ar 1320 814 480 480  1.0000000 1.1250000}%
\special{ar 1320 814 480 480  1.2000000 1.3250000}%
\special{ar 1320 814 480 480  1.4000000 1.5250000}%
\special{ar 1320 814 480 480  1.6000000 1.7250000}%
\special{ar 1320 814 480 480  1.8000000 1.9250000}%
\special{ar 1320 814 480 480  2.0000000 2.1250000}%
\special{ar 1320 814 480 480  2.2000000 2.3250000}%
\special{ar 1320 814 480 480  2.4000000 2.5250000}%
\special{ar 1320 814 480 480  2.6000000 2.7250000}%
\special{ar 1320 814 480 480  2.8000000 2.9250000}%
\special{ar 1320 814 480 480  3.0000000 3.1250000}%
\special{ar 1320 814 480 480  3.2000000 3.3250000}%
\special{ar 1320 814 480 480  3.4000000 3.5250000}%
\special{ar 1320 814 480 480  3.6000000 3.7250000}%
\special{ar 1320 814 480 480  3.8000000 3.9250000}%
\special{ar 1320 814 480 480  4.0000000 4.1250000}%
\special{ar 1320 814 480 480  4.2000000 4.3250000}%
\special{ar 1320 814 480 480  4.4000000 4.5250000}%
\special{ar 1320 814 480 480  4.6000000 4.7250000}%
\special{ar 1320 814 480 480  4.8000000 4.9250000}%
\special{ar 1320 814 480 480  5.0000000 5.1250000}%
\special{ar 1320 814 480 480  5.2000000 5.3250000}%
\special{ar 1320 814 480 480  5.4000000 5.5250000}%
\special{ar 1320 814 480 480  5.6000000 5.7250000}%
\special{ar 1320 814 480 480  5.8000000 5.9250000}%
\special{ar 1320 814 480 480  6.0000000 6.1250000}%
\special{ar 1320 814 480 480  6.2000000 6.2832853}%
}}%
%
{\color[named]{Black}{%
\special{pn 8}%
\special{ar 3400 814 480 480  0.0000000 0.1250000}%
\special{ar 3400 814 480 480  0.2000000 0.3250000}%
\special{ar 3400 814 480 480  0.4000000 0.5250000}%
\special{ar 3400 814 480 480  0.6000000 0.7250000}%
\special{ar 3400 814 480 480  0.8000000 0.9250000}%
\special{ar 3400 814 480 480  1.0000000 1.1250000}%
\special{ar 3400 814 480 480  1.2000000 1.3250000}%
\special{ar 3400 814 480 480  1.4000000 1.5250000}%
\special{ar 3400 814 480 480  1.6000000 1.7250000}%
\special{ar 3400 814 480 480  1.8000000 1.9250000}%
\special{ar 3400 814 480 480  2.0000000 2.1250000}%
\special{ar 3400 814 480 480  2.2000000 2.3250000}%
\special{ar 3400 814 480 480  2.4000000 2.5250000}%
\special{ar 3400 814 480 480  2.6000000 2.7250000}%
\special{ar 3400 814 480 480  2.8000000 2.9250000}%
\special{ar 3400 814 480 480  3.0000000 3.1250000}%
\special{ar 3400 814 480 480  3.2000000 3.3250000}%
\special{ar 3400 814 480 480  3.4000000 3.5250000}%
\special{ar 3400 814 480 480  3.6000000 3.7250000}%
\special{ar 3400 814 480 480  3.8000000 3.9250000}%
\special{ar 3400 814 480 480  4.0000000 4.1250000}%
\special{ar 3400 814 480 480  4.2000000 4.3250000}%
\special{ar 3400 814 480 480  4.4000000 4.5250000}%
\special{ar 3400 814 480 480  4.6000000 4.7250000}%
\special{ar 3400 814 480 480  4.8000000 4.9250000}%
\special{ar 3400 814 480 480  5.0000000 5.1250000}%
\special{ar 3400 814 480 480  5.2000000 5.3250000}%
\special{ar 3400 814 480 480  5.4000000 5.5250000}%
\special{ar 3400 814 480 480  5.6000000 5.7250000}%
\special{ar 3400 814 480 480  5.8000000 5.9250000}%
\special{ar 3400 814 480 480  6.0000000 6.1250000}%
\special{ar 3400 814 480 480  6.2000000 6.2832853}%
}}%
\put(22.5600,-8.1400){\makebox(0,0)[lb]{$\longleftrightarrow$}}%
\put(2.0000,-1.7400){\makebox(0,0)[lb]{($\omega 3$) \textit{jumping over a double point}}}%
%
{\color[named]{Black}{%
\special{pn 13}%
\special{pa 840 814}%
\special{pa 1080 814}%
\special{fp}%
\special{pa 1800 814}%
\special{pa 1560 814}%
\special{fp}%
}}%
%
{\color[named]{Black}{%
\special{pn 13}%
\special{ar 1320 814 240 240  3.1415927 6.2831853}%
}}%
%
{\color[named]{Black}{%
\special{pn 13}%
\special{ar 3400 814 240 240  6.2831853 6.2831853}%
\special{ar 3400 814 240 240  0.0000000 3.1415927}%
}}%
%
{\color[named]{Black}{%
\special{pn 13}%
\special{pa 2920 814}%
\special{pa 3160 814}%
\special{fp}%
\special{pa 3880 814}%
\special{pa 3640 814}%
\special{fp}%
}}%
%
{\color[named]{Black}{%
\special{pn 13}%
\special{pa 1640 494}%
\special{pa 1000 1134}%
\special{fp}%
\special{pa 1000 494}%
\special{pa 1640 1134}%
\special{fp}%
}}%
%
{\color[named]{Black}{%
\special{pn 13}%
\special{pa 3720 494}%
\special{pa 3080 1134}%
\special{fp}%
\special{pa 3080 494}%
\special{pa 3720 1134}%
\special{fp}%
}}%
%
{\color[named]{Black}{%
\special{pn 4}%
\special{sh 1}%
\special{ar 1152 646 26 26 0  6.28318530717959E+0000}%
}}%
%
{\color[named]{Black}{%
\special{pn 4}%
\special{sh 1}%
\special{ar 1488 646 26 26 0  6.28318530717959E+0000}%
}}%
%
{\color[named]{Black}{%
\special{pn 4}%
\special{sh 1}%
\special{ar 1320 814 26 26 0  6.28318530717959E+0000}%
}}%
%
{\color[named]{Black}{%
\special{pn 4}%
\special{sh 1}%
\special{ar 3400 814 26 26 0  6.28318530717959E+0000}%
}}%
%
{\color[named]{Black}{%
\special{pn 4}%
\special{sh 1}%
\special{ar 3224 990 26 26 0  6.28318530717959E+0000}%
}}%
%
{\color[named]{Black}{%
\special{pn 4}%
\special{sh 1}%
\special{ar 3568 990 26 26 0  6.28318530717959E+0000}%
}}%
\end{picture}%
\end{center}
\end{figure}

Suppose $\gamma$ and $\gamma^{\prime}$ are related by the move ($\omega 1$).
The contribution of the double point in the right picture of the move ($\omega 1$)
is zero, since the class of a null-homotopic loop is zero in $\mathbb{Q}\hat{\pi}^{\prime}(S)$.
Hence $\mu(\gamma)=\mu(\gamma^{\prime})$.

Suppose $\gamma$ and $\gamma^{\prime}$ are related by the move ($\omega 2$).
We may assume the left picture corresponds to $\gamma^{\prime}$.
Then $\gamma$ has two more double points than $\gamma^{\prime}$.
We write them by $p$ and $q$ so that $t_1^p < t_1^q$.
As in Figure 3, there are two possibilities: $t_2^q< t_2^p$ or $t_2^p< t_2^q$,
but in any case, $\gamma_{0 t_1^p} \gamma_{t_2^p 1}$ is homotopic
to $\gamma_{0 t_1^q} \gamma_{t_2^q 1}$ relative to $\{ 0,1 \}$,
$|\gamma_{t_1^p t_2^p}|=|\gamma_{t_1^q t_2^q}|$,
and $\varepsilon(\dot{\gamma}(t_1^p),\dot{\gamma}(t_2^p))
=-\varepsilon(\dot{\gamma}(t_1^q),\dot{\gamma}(t_2^q))$.
Hence the contributions from $p$ and $q$ cancel and
$\mu(\gamma)=\mu(\gamma^{\prime})$.

\begin{figure}
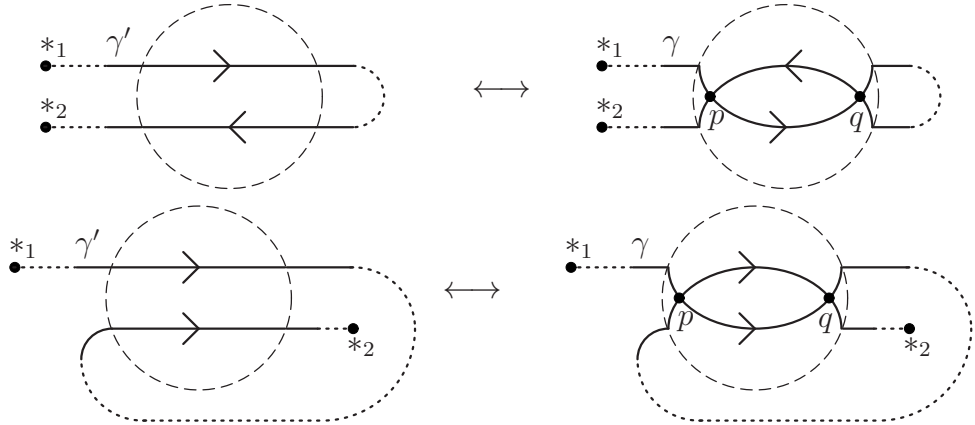

\begin{center}
\caption{invariance under the move ($\omega 2$)}

\vspace{0.2cm}

\input{invariance-move2-1.tex}

\vspace{0.2cm}

\input{invariance-move2-2.tex}
\end{center}
\end{figure}

Suppose $\gamma$ and $\gamma^{\prime}$ are related by the move ($\omega 3$).
Similarly to the case of ($\omega 2$), we see that a cancel happens and
$\mu(\gamma)=\mu(\gamma^{\prime})$. Typical cases are illustrated in Figure 4,
where the contributions from $p$ (resp. $q$) and
$p^{\prime}$ (resp. $q^{\prime}$) cancel.

\begin{figure}
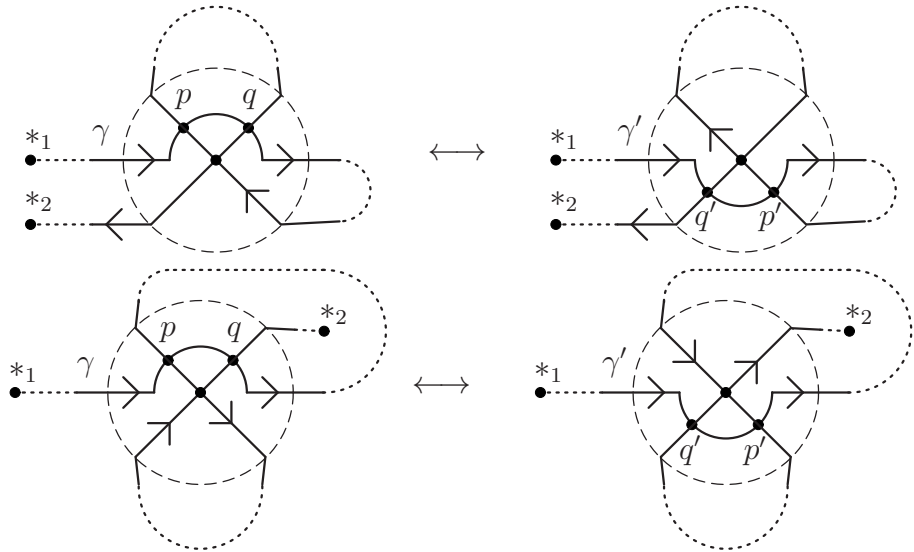

\begin{center}
\caption{invariance under the move ($\omega 3$)}

\vspace{0.2cm}

\input{invariance-move3-1.tex}

\vspace{0.2cm}

\input{invariance-move3-2.tex}

\end{center}
\end{figure}

We next show that $\mathbb{Q}\Pi S(*_0,*_1)$ is a right
$\mathbb{Q}\hat{\pi}^{\prime}(S)$-comodule, i.e.,
$(1_{\mathbb{Q}\Pi S(*_0,*_1)} \otimes (1-T))
(\mu \otimes 1_{\mathbb{Q}\hat{\pi}^{\prime}(S)}) \circ \mu
=(1_{\mathbb{Q}\Pi S(*_0,*_1)} \otimes \delta)\circ \mu$ (see Appendix).
Let $\gamma\colon [0,1]\to S$ be an immersed path from $*_0$ to $*_1$
such that its self intersections consist of transverse double points.
To compute $(1_{\mathbb{Q}\Pi S(*_0,*_1)} \otimes (1-T))
(\mu \otimes 1_{\mathbb{Q}\hat{\pi}^{\prime}(S)}) \mu(\gamma)$,
we need to compute $\mu(\gamma_{0 t_1^p} \gamma_{t_2^p 1})$
where $p\in \Gamma$. The double points of $\gamma_{0 t_1^p} \gamma_{t_2^p 1}$
come from those of $\gamma$. Let $q$ be a double point of
$\gamma_{0 t_1^p} \gamma_{t_2^p 1}$ and denote $\gamma^{-1}(q)
=\{ t_1^q, t_2^q \}$, so that $t_1^q<t_2^q$. There are three possibilities:
(i) $t_1^q<t_2^q<t_1^p<t_2^p$, (ii) $t_1^p<t_2^p<t_1^q<t_2^q$,
(iii) $t_1^q<t_1^p<t_2^p<t_2^q$. In cases (i) and (ii),
the contribution to $(\mu \otimes 1_{\mathbb{Q}\hat{\pi}^{\prime}(S)})\mu(\gamma)$
from $(p,q)$ is
\begin{eqnarray*}
& &
\varepsilon(\dot{\gamma}(t_1^p), \dot{\gamma}(t_2^p))
\varepsilon(\dot{\gamma}(t_1^q), \dot{\gamma}(t_2^q))
\left( \gamma_{0 t_1^q} \gamma_{t_2^q t_1^p} \gamma_{t_2^p 1} \right)
\otimes |\gamma_{t_1^q t_2^q}|^{\prime} \otimes |\gamma_{t_1^p t_2^p}|^{\prime}
\quad {\rm and} \\
& &
\varepsilon(\dot{\gamma}(t_1^p), \dot{\gamma}(t_2^p))
\varepsilon(\dot{\gamma}(t_1^q), \dot{\gamma}(t_2^q))
\left( \gamma_{0 t_1^p} \gamma_{t_2^p t_1^q} \gamma_{t_2^q 1} \right)
\otimes |\gamma_{t_1^q t_2^q}|^{\prime} \otimes |\gamma_{t_1^p t_2^p}|^{\prime},
\end{eqnarray*}
respectively. Here $\gamma_{t_2^q t_1^p}$ means the restriction of
$\gamma$ to $[t_2^q, t_1^p]$ and $\gamma_{0 t_1^q}
\gamma_{t_2^q t_1^p} \gamma_{t_2^p 1}$ means the conjunction.
Therefore the contributions to $(\mu \otimes 1_{\mathbb{Q}\hat{\pi}^{\prime}(S)})\mu(\gamma)$
from $(p,q)$ of type (i) or (ii) are written as a linear combination of tensors of the form
$u\otimes (v\otimes w+ w\otimes v)$. Since $(1-T)(v\otimes w +w\otimes v)=0$,
these contributions vanish on $\mathbb{Q}\Pi S(*_0,*_1) \otimes
\mathbb{Q}\hat{\pi}^{\prime}(S)\otimes \mathbb{Q}\hat{\pi}^{\prime}(S)$.
Hence we only need to consider the contributions from (iii), and
\begin{eqnarray}
&  & (1_{\mathbb{Q}\Pi S(*_0,*_1)} \otimes (1-T))
(\mu \otimes 1_{\mathbb{Q}\hat{\pi}^{\prime}(S)}) \mu(\gamma) \nonumber \\
& = & \sum_{\substack{p,q \\ t_1^q<t_1^p<t_2^p<t_2^q}}
\varepsilon(\dot{\gamma}(t_1^p), \dot{\gamma}(t_2^p))
\varepsilon(\dot{\gamma}(t_1^q), \dot{\gamma}(t_2^q))
z(\gamma, p, q) \label{z(gamma)},
\end{eqnarray}
where
$z(\gamma,p,q)=
(\gamma_{0 t_1^q} \gamma_{t_2^q 1}) \otimes
\left( |\gamma_{t_1^q t_1^p} \gamma_{t_2^p t_2^q}|^{\prime} \otimes
|\gamma_{t_1^p t_2^p}|^{\prime}-|\gamma_{t_1^p t_2^p}|^{\prime} \otimes
|\gamma_{t_1^q t_1^p} \gamma_{t_2^p t_2^q}|^{\prime}
\right)$.
On the other hand, to compute
$(1_{\mathbb{Q}\Pi S(*_0,*_1)} \otimes \delta)\mu(\gamma)$ we need to
compute $\delta(|\gamma_{t_1^p t_2^p}|^{\prime})$ where $p\in \Gamma$.
Each double point of the loop $\gamma_{t_1^p t_2^p}$ comes from
$q \in \Gamma$ such that $t_1^p < t_1^q < t_2^q < t_2^p$.
Thus $\delta(|\gamma_{t_1^p t_2^p}|^{\prime})$ is equal to
$$\sum_{q; t_1^p< t_1^q < t_2^q < t_2^p}
\varepsilon(\dot{\gamma}(t_1^q),\dot{\gamma}(t_2^q))
\left(
|\gamma_{t_1^q t_2^q}|^{\prime} \otimes |\gamma_{t_2^q t_2^p} \gamma_{t_1^p t_1^q}|^{\prime}
-|\gamma_{t_2^q t_2^p} \gamma_{t_1^p t_1^q}|^{\prime} \otimes |\gamma_{t_1^q t_2^q}|^{\prime}
\right),\quad {\rm and}$$
\begin{eqnarray}
& &
(1_{\mathbb{Q}\Pi S(*_0,*_1)} \otimes \delta)\mu(\gamma) \nonumber \\
& = & -\sum_{\substack{p,q \\ t_1^p< t_1^q < t_2^q < t_2^p}}
\varepsilon(\dot{\gamma}(t_1^p), \dot{\gamma}(t_2^p))
\varepsilon(\dot{\gamma}(t_1^q), \dot{\gamma}(t_2^q))
w(\gamma, p,q), \label{w(gamma)}
\end{eqnarray}
where
$w(\gamma,p,q)=(\gamma_{0 t_1^p} \gamma_{t_2^p 1}) \otimes
\left( |\gamma_{t_1^q t_2^q}|^{\prime} \otimes |\gamma_{t_2^q t_2^p} \gamma_{t_1^p t_1^q}|^{\prime}
-|\gamma_{t_2^q t_2^p} \gamma_{t_1^p t_1^q}|^{\prime} \otimes \gamma_{t_1^q t_2^q}
\right)$. Since $z(\gamma,p,q)=-w(\gamma,q,p)$, the right hand sides of
(\ref{z(gamma)}) and (\ref{w(gamma)}) are equal. This completes
the proof.
\end{proof}

\begin{rem}
\label{rem-*}
{\rm We have taken $*_0$ and $*_1$ from $\partial S$.
If at least one of $*_0$ and $*_1$ lies in ${\rm Int}(S)$,
we need to consider another kind of local move illustrated in Figure 5.
In this case the formula (\ref{dfn-mu}) does not work. For example,
in Figure 5, the contribution from $p$ in the left picture is non-trivial.
}
\end{rem}

\begin{figure}
\begin{center}
\caption{the formula does not work if a base point lies in ${\rm Int}(S)$}

\vspace{0.2cm}

\unitlength 0.1in
\begin{picture}( 31.7300, 10.7900)(  5.4000,-12.5300)
%
{\color[named]{Black}{%
\special{pn 8}%
\special{ar 1068 840 226 226  3.9269908 5.4977871}%
}}%
%
{\color[named]{Black}{%
\special{pn 8}%
\special{ar 1068 504 264 264  0.6435011 2.4980915}%
}}%
%
{\color[named]{Black}{%
\special{pn 4}%
\special{sh 1}%
\special{ar 588 1000 26 26 0  6.28318530717959E+0000}%
}}%
%
{\color[named]{Black}{%
\special{pn 13}%
\special{pa 588 1000}%
\special{pa 598 970}%
\special{pa 610 938}%
\special{pa 620 908}%
\special{pa 630 878}%
\special{pa 654 818}%
\special{pa 668 788}%
\special{pa 682 760}%
\special{pa 696 732}%
\special{pa 712 704}%
\special{pa 748 650}%
\special{pa 768 624}%
\special{pa 790 600}%
\special{pa 812 576}%
\special{pa 838 554}%
\special{pa 864 532}%
\special{pa 892 514}%
\special{pa 922 496}%
\special{pa 952 480}%
\special{pa 984 466}%
\special{pa 1016 456}%
\special{pa 1048 446}%
\special{pa 1080 442}%
\special{pa 1114 438}%
\special{pa 1148 440}%
\special{pa 1180 444}%
\special{pa 1214 452}%
\special{pa 1246 464}%
\special{pa 1278 480}%
\special{pa 1308 498}%
\special{pa 1336 520}%
\special{pa 1362 544}%
\special{pa 1384 570}%
\special{pa 1402 598}%
\special{pa 1416 626}%
\special{pa 1426 656}%
\special{pa 1428 684}%
\special{pa 1426 714}%
\special{pa 1418 742}%
\special{pa 1404 770}%
\special{pa 1386 796}%
\special{pa 1362 822}%
\special{pa 1336 846}%
\special{pa 1308 868}%
\special{pa 1276 890}%
\special{pa 1242 908}%
\special{pa 1206 926}%
\special{pa 1168 940}%
\special{pa 1130 952}%
\special{pa 1092 962}%
\special{pa 1056 968}%
\special{pa 1018 972}%
\special{pa 982 972}%
\special{pa 950 970}%
\special{pa 918 964}%
\special{pa 888 952}%
\special{pa 862 938}%
\special{pa 840 922}%
\special{pa 818 900}%
\special{pa 800 878}%
\special{pa 782 852}%
\special{pa 766 824}%
\special{pa 752 794}%
\special{pa 740 762}%
\special{pa 728 730}%
\special{pa 694 628}%
\special{pa 684 594}%
\special{pa 674 560}%
\special{pa 664 526}%
\special{pa 640 462}%
\special{pa 626 432}%
\special{pa 612 404}%
\special{pa 596 376}%
\special{pa 588 360}%
\special{sp}%
}}%
%
{\color[named]{Black}{%
\special{pn 4}%
\special{sh 1}%
\special{ar 588 360 26 26 0  6.28318530717959E+0000}%
}}%
%
{\color[named]{Black}{%
\special{pn 8}%
\special{ar 3308 840 226 226  3.9269908 5.4977871}%
}}%
%
{\color[named]{Black}{%
\special{pn 8}%
\special{ar 3308 504 264 264  0.6435011 2.4980915}%
}}%
%
{\color[named]{Black}{%
\special{pn 4}%
\special{sh 1}%
\special{ar 2828 1000 26 26 0  6.28318530717959E+0000}%
}}%
%
{\color[named]{Black}{%
\special{pn 4}%
\special{sh 1}%
\special{ar 2828 360 26 26 0  6.28318530717959E+0000}%
}}%
%
{\color[named]{Black}{%
\special{pn 13}%
\special{pa 2828 1000}%
\special{pa 2868 868}%
\special{pa 2890 804}%
\special{pa 2902 774}%
\special{pa 2916 744}%
\special{pa 2930 714}%
\special{pa 2944 684}%
\special{pa 2960 656}%
\special{pa 2976 630}%
\special{pa 2996 604}%
\special{pa 3016 580}%
\special{pa 3036 556}%
\special{pa 3060 536}%
\special{pa 3084 516}%
\special{pa 3110 496}%
\special{pa 3140 480}%
\special{pa 3170 466}%
\special{pa 3202 452}%
\special{pa 3236 442}%
\special{pa 3270 432}%
\special{pa 3306 426}%
\special{pa 3340 422}%
\special{pa 3376 420}%
\special{pa 3410 420}%
\special{pa 3446 422}%
\special{pa 3478 428}%
\special{pa 3510 436}%
\special{pa 3542 446}%
\special{pa 3570 460}%
\special{pa 3596 478}%
\special{pa 3620 496}%
\special{pa 3642 520}%
\special{pa 3660 544}%
\special{pa 3676 572}%
\special{pa 3690 602}%
\special{pa 3700 634}%
\special{pa 3708 668}%
\special{pa 3712 702}%
\special{pa 3714 736}%
\special{pa 3714 770}%
\special{pa 3710 806}%
\special{pa 3704 840}%
\special{pa 3694 874}%
\special{pa 3684 908}%
\special{pa 3668 940}%
\special{pa 3652 970}%
\special{pa 3632 1000}%
\special{pa 3612 1028}%
\special{pa 3588 1054}%
\special{pa 3562 1078}%
\special{pa 3536 1102}%
\special{pa 3506 1124}%
\special{pa 3476 1144}%
\special{pa 3444 1162}%
\special{pa 3412 1180}%
\special{pa 3378 1196}%
\special{pa 3342 1208}%
\special{pa 3308 1220}%
\special{pa 3272 1230}%
\special{pa 3234 1240}%
\special{pa 3198 1246}%
\special{pa 3160 1250}%
\special{pa 3124 1252}%
\special{pa 3086 1254}%
\special{pa 3050 1252}%
\special{pa 3012 1248}%
\special{pa 2978 1244}%
\special{pa 2942 1236}%
\special{pa 2908 1226}%
\special{pa 2874 1214}%
\special{pa 2842 1200}%
\special{pa 2812 1184}%
\special{pa 2784 1166}%
\special{pa 2756 1146}%
\special{pa 2730 1122}%
\special{pa 2706 1098}%
\special{pa 2684 1070}%
\special{pa 2662 1042}%
\special{pa 2644 1012}%
\special{pa 2628 982}%
\special{pa 2614 950}%
\special{pa 2602 916}%
\special{pa 2592 884}%
\special{pa 2584 850}%
\special{pa 2578 816}%
\special{pa 2576 784}%
\special{pa 2576 750}%
\special{pa 2578 718}%
\special{pa 2582 686}%
\special{pa 2590 654}%
\special{pa 2600 624}%
\special{pa 2612 596}%
\special{pa 2626 570}%
\special{pa 2644 544}%
\special{pa 2664 518}%
\special{pa 2684 494}%
\special{pa 2706 470}%
\special{pa 2730 448}%
\special{pa 2754 424}%
\special{pa 2804 380}%
\special{pa 2828 360}%
\special{sp}%
}}%
%
{\color[named]{Black}{%
\special{pn 4}%
\special{sh 1}%
\special{ar 708 696 26 26 0  6.28318530717959E+0000}%
}}%
%
{\color[named]{Black}{%
\special{pn 13}%
\special{pa 908 504}%
\special{pa 908 616}%
\special{fp}%
\special{pa 908 504}%
\special{pa 796 504}%
\special{fp}%
}}%
%
{\color[named]{Black}{%
\special{pn 13}%
\special{pa 3004 600}%
\special{pa 3028 720}%
\special{fp}%
\special{pa 3004 600}%
\special{pa 2868 648}%
\special{fp}%
}}%
\put(5.4000,-11.9000){\makebox(0,0)[lb]{$*_1 \in {\rm Int}(S)$}}%
\put(28.4400,-11.6000){\makebox(0,0)[lb]{$*_1$}}%
\put(6.0400,-3.0400){\makebox(0,0)[lb]{$*_2$}}%
\put(28.4400,-3.0400){\makebox(0,0)[lb]{$*_2$}}%
\put(17.8000,-6.8000){\makebox(0,0)[lb]{$\longleftrightarrow$}}%
\put(5.7200,-7.2000){\makebox(0,0)[lb]{$p$}}%
\end{picture}%

\end{center}
\end{figure}
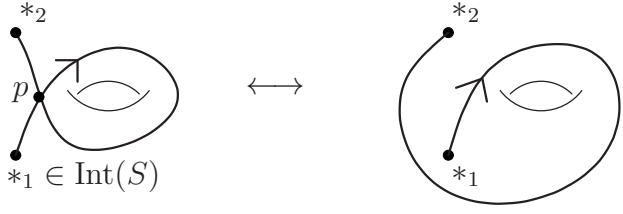

We show that $\sigma$ and $\mu$ satisfy the compatibility (\ref{A6compatible+})
and the involutivity (\ref{A6involutive}).

\begin{prop}
\label{GTbim}
The $\mathbb{Q}$-vector space $\mathbb{Q}\Pi S(*_0,*_1)$ is an involutive
right $\mathbb{Q}\hat{\pi}^{\prime}(S)$-bimodule
with respect to $\sigma\colon
\mathbb{Q}\hat{\pi}^{\prime}(S) \otimes \mathbb{Q}\Pi S(*_0,*_1)
\to \mathbb{Q}\Pi S(*_0,*_1)$ and
$\mu\colon \mathbb{Q}\Pi S(*_0,*_1) \to
\mathbb{Q}\Pi S(*_0,*_1) \otimes \mathbb{Q}\hat{\pi}^{\prime}(S)$.
\end{prop}

\begin{proof}
We first prove the involutivity.
Let $\gamma \colon [0,1]\to S$ be an immersed path from
$*_0$ to $*_1$ such that its self intersections consist of transverse
double points. Then
$$\overline{\sigma}\mu(\gamma)=
\sum_{p\in \Gamma} \varepsilon(\dot{\gamma}(t_1^p), \dot{\gamma}(t_2^p))
\sigma \left( |\gamma_{t_1^p t_2^p}|^{\prime} \otimes
(\gamma_{0 t_1^p} \gamma_{t_2^p 1}) \right).$$
Let $q$ be an intersection of the loop $\gamma_{t_1^p t_2^p}$ and
the path $\gamma_{0 t_1^p} \gamma_{t_2^p 1}$ and we denote
$\gamma^{-1}(q)=\{ t_1^q, t_2^q \}$ so that $t_1^q < t_2^q$.
There are two possibilities:
(i) $t_1^q < t_1^p < t_2^q < t_2^p$, (ii) $t_1^p < t_1^q < t_2^p < t_2^q$.
The contribution to $\overline{\sigma}\mu(\gamma)$ from $(p,q)$ are
\begin{equation}
\label{(i)cncl}
\varepsilon(\dot{\gamma}(t_1^p), \dot{\gamma}(t_2^p))
\varepsilon(\dot{\gamma}(t_2^q), \dot{\gamma}(t_1^q))
\gamma_{0 t_1^q} \gamma_{t_2^q t_2^p} \gamma_{t_1^p t_2^q}
\gamma_{t_1^q t_1^p} \gamma_{t_2^p 1}
\end{equation}
in case (i), and
\begin{equation}
\label{(ii)cncl}
\varepsilon(\dot{\gamma}(t_1^p), \dot{\gamma}(t_2^p))
\varepsilon(\dot{\gamma}(t_1^q), \dot{\gamma}(t_2^q))
\gamma_{0 t_1^p} \gamma_{t_2^p t_2^q} \gamma_{t_1^q t_2^p}
\gamma_{t_1^p t_1^q} \gamma_{t_2^q 1}
\end{equation}
in case (ii).
If we interchange $p$ and $q$ in (\ref{(ii)cncl}), we get
the minus of (\ref{(i)cncl}). Therefore the contributions from
$(p,q)$ in case (i) and those in case (ii) cancel and $\overline{\sigma}\mu(\gamma)=0$.

We next show the compatibility.
Let $\alpha$ be an immersed loop on $S$
and $\gamma \colon [0,1]\to S$ an immersed path from $*_0$ to $*_1$ such that
their intersections and self intersections consist of transverse double points.
The compatibility is equivalent to the following.
\begin{equation}
\label{mu-sigma}
\mu(\sigma(\alpha\otimes \gamma))=\sigma(\alpha)\mu(\gamma)
-(\overline{\sigma}\otimes 1_{\mathbb{Q}\hat{\pi}^{\prime}(S)})
(\gamma \otimes \delta(\alpha)).
\end{equation}
Here,
$\sigma(\alpha)\mu(\gamma)=(\sigma \otimes
1_{\mathbb{Q}\hat{\pi}^{\prime}(S)})(\alpha \otimes \mu(\gamma))
+(1_{\mathbb{Q}\Pi S(*_0,*_1)} \otimes {\rm ad}(\alpha))\mu(\gamma)$.
We compute the left hand side of (\ref{mu-sigma}). First of all, we have
$$\mu(\sigma(\alpha \otimes \gamma))=
\sum_{p\in \alpha \cap \gamma} \varepsilon(p;\alpha,\gamma)
\mu(\gamma_{*_0 p}\alpha_p\gamma_{p*_1}).$$
Let $q$ be a double point of $\gamma_{*_0 p}\alpha_p\gamma_{p*_1}$.
There are three possibilities: (i) $q$ comes from a double point
of $\alpha$, (ii) $q$ comes from a double point of $\gamma$,
(iii) $q$ comes from an intersection of $\alpha$ and $\gamma$,
which is different from $p$.

Suppose $q$ comes from a double point of $\alpha$.
We denote $\alpha^{-1}(q)=\{ t_1^q, t_2^q \}\subset S^1$, so that
$t_2^q$, $\alpha^{-1}(p)$, $t_1^q$ are arranged in this order
according to the orientation of $S^1$. (Since $p$ is a simple point
of $\alpha$, the preimage $\alpha^{-1}(p)$ consists of one point.
For simplicity, we write $\alpha^{-1}(p)$ for the unique point
in the preimage.)
The contribution
to $\mu(\gamma_{*_0p}\alpha_p\gamma_{p*_1})$ from such $q$ is
$$-\varepsilon(\dot{\alpha}(t_1^q),\dot{\alpha}(t_2^q))
\gamma_{*_0p}\alpha_{pq}\alpha_{qp}\gamma_{p*_1} \otimes |\alpha_{t_1^q t_2^q}|^{\prime}.$$
Here, $\alpha_{pq}$ (resp. $\alpha_{qp}$) means the restriction of $\alpha$ to the interval
$[\alpha^{-1}(p),t_1^q]$ (resp. $[t_2^q, \alpha^{-1}(p)]$).
Thus the contributions to $\mu(\sigma(\alpha \otimes \gamma))$
from $(p,q)$ such that $q$ is of type (i) is
\begin{equation}
\label{TYPE1}
-\sum_{p\in \alpha \cap \gamma} \sum_{(t_1^q,t_2^q)}
\varepsilon(p;\alpha,\gamma)
\varepsilon(\dot{\alpha}(t_1^q), \dot{\alpha}(t_2^q))
\gamma_{*_0p}\alpha_{pq}\alpha_{qp}\gamma_{p*_1} \otimes |\alpha_{t_1^q t_2^q}|^{\prime},
\end{equation}
where the second sum is taken over ordered pairs $(t_1^q, t_2^q)$ such that
$\alpha(t_1^q)=\alpha(t_2^q)$ and $\alpha^{-1}(p)\in [t_2^q, t_1^q]$.
On the other hand, we have
$$\delta(\alpha)=\sum_{(t_1^q,t_2^q)}
\varepsilon(\dot{\alpha}(t_2^q),\dot{\alpha}(t_1^q))|\alpha_{t_2^q t_1^q}|^{\prime}
\otimes |\alpha_{t_1^q t_2^q}|^{\prime},$$
where the sum is taken over ordered pairs $(t_1^q, t_2^q)$ such that
$\alpha(t_1^q)=\alpha(t_2^q)$, $t_1^q \neq t_2^q$, and
$$\sigma(|\alpha_{t_2^q t_1^q}| \otimes \gamma)
=\sum_p \varepsilon(p;\alpha,\gamma)
\gamma_{*_0p}\alpha_{pq}\alpha_{qp}\gamma_{p*_1},$$
where the sum is taken over $p\in \alpha \cap \gamma$ such that
$\alpha^{-1}(p)\in [t_2^q, t_1^q]$. Therefore, (\ref{TYPE1})
is equal to $-(\overline{\sigma}\otimes 1_{\mathbb{Q}\hat{\pi}^{\prime}(S)})
(\gamma \otimes \delta(\alpha))$.

Suppose $q$ comes from a double point of $\gamma$.
We denote $\gamma^{-1}(q)=\{ s_1^q, s_2^q \}$, so that $s_1^q < s_2^q$.
There are three possibilities: (ii-a)
$\gamma^{-1}(p) < s_1^q < s_2^q$, (ii-b) $s_1^q < \gamma^{-1}(p) < s_2^q$,
(ii-c) $s_1^q < s_2^q < \gamma^{-1}(p)$. The contributions
to $\mu(\sigma(\alpha\otimes \gamma))$ from $(p,q)$ of type (ii-a) are
$$\sum_{p\in \alpha \cap \gamma}\sum_{\substack{q \\ \gamma^{-1}(p)< s_1^q < s_2^q}}
\varepsilon(p;\alpha,\gamma)
\varepsilon(\dot{\gamma}(s_1^q),\dot{\gamma}(s_2^q))
\gamma_{*_0 p}\alpha_p\gamma_{pq}\gamma_{q*_1} \otimes \gamma_{s_1^q s_2^q},$$
and those from $(p,q)$ of type (ii-c) are
$$\sum_{p\in \alpha \cap \gamma}\sum_{\substack{q \\ s_1^q < s_2^q < \gamma^{-1}(p)}}
\varepsilon(p;\alpha,\gamma)
\varepsilon(\dot{\gamma}(s_1^q),\dot{\gamma}(s_2^q))
\gamma_{*_0 q}\gamma_{qp}\alpha_p \gamma_{p *_1} \otimes \gamma_{s_1^q s_2^q}.$$
The sum of these two is equal to
$(\sigma \otimes 1_{\mathbb{Q}\hat{\pi}^{\prime}(S)})(\alpha \otimes \mu(\gamma))$.
The contributions to $\mu(\sigma(\alpha\otimes \gamma))$ from $(p,q)$ of type (ii-b) are
$$\sum_{p\in \alpha \cap \gamma}\sum_{\substack{q \\ s_1^q < \gamma^{-1}(p) < s_2^q}}
\varepsilon(p;\alpha,\gamma)
\varepsilon(\dot{\gamma}(s_1^q),\dot{\gamma}(s_2^q))
\gamma_{*_0 q}\gamma_{q *_1} \otimes \gamma_{qp}\alpha_p \gamma_{pq}.$$
This is equal to $(1_{\mathbb{Q}\Pi S(*_0,*_1)}\otimes
{\rm ad}(\alpha))\mu(\gamma)$. Therefore,
the contributions from $(p,q)$ such that $q$ is of type (ii)
is $(\sigma \otimes 1_{\mathbb{Q}\hat{\pi}^{\prime}(S)})(\alpha \otimes \mu(\gamma))
+(1_{\mathbb{Q}\Pi S(*_0,*_1)}\otimes {\rm ad}(\alpha))
\mu(\gamma)=\sigma(\alpha)\mu(\gamma)$.

Suppose $q$ comes from an intersection of $\alpha$ and $\gamma$,
which is different from $p$. If $\gamma^{-1}(p)<\gamma^{-1}(q)$, the contribution is
\begin{equation}
\label{(iii)cl1}
-\varepsilon(p;\alpha,\gamma)\varepsilon(q;\alpha,\gamma)
\gamma_{*_0p}\alpha_{pq}\gamma_{q*_1} \otimes \alpha_{qp}\gamma_{pq},
\end{equation}
and if $\gamma^{-1}(q)< \gamma^{-1}(p)$, the contribution is
\begin{equation}
\label{(iii)cl2}
-\varepsilon(p;\alpha,\gamma)\varepsilon(q;\gamma,\alpha)
\gamma_{*_0q}\alpha_{qp}\gamma_{p*_1} \otimes \gamma_{qp} \alpha_{pq}.
\end{equation}
See Figure 6.
If we interchange $p$ and $q$ in (\ref{(iii)cl2}), we get the minus
of (\ref{(iii)cl1}). Therefore the sum of the contributions
from $(p,q)$ of type (iii) is zero.
We have established the formula (\ref{mu-sigma}).
\end{proof}

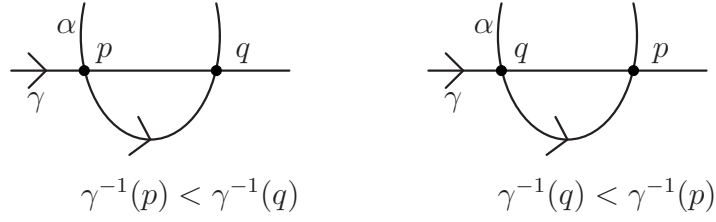
\begin{figure}
\begin{center}
\caption{the type (iii)}

\vspace{0.2cm}

\unitlength 0.1in
\begin{picture}( 36.0000,  9.5900)(  4.0000,-11.6900)
%
{\color[named]{Black}{%
\special{pn 13}%
\special{pa 400 562}%
\special{pa 1840 562}%
\special{fp}%
}}%
%
{\color[named]{Black}{%
\special{pn 13}%
\special{ar 1120 382 360 540  5.9614348 6.2831853}%
\special{ar 1120 382 360 540  0.0000000 3.4633432}%
}}%
%
{\color[named]{Black}{%
\special{pn 4}%
\special{sh 1}%
\special{ar 778 562 26 26 0  6.28318530717959E+0000}%
}}%
%
{\color[named]{Black}{%
\special{pn 4}%
\special{sh 1}%
\special{ar 1462 562 26 26 0  6.28318530717959E+0000}%
}}%
%
{\color[named]{Black}{%
\special{pn 13}%
\special{pa 580 562}%
\special{pa 490 472}%
\special{fp}%
\special{pa 580 562}%
\special{pa 490 652}%
\special{fp}%
}}%
%
{\color[named]{Black}{%
\special{pn 13}%
\special{pa 1120 922}%
\special{pa 1058 804}%
\special{fp}%
\special{pa 1120 922}%
\special{pa 1012 1002}%
\special{fp}%
}}%
\put(4.8100,-7.7700){\makebox(0,0)[lb]{$\gamma$}}%
\put(6.3400,-3.6300){\makebox(0,0)[lb]{$\alpha$}}%
\put(8.4100,-5.2500){\makebox(0,0)[lb]{$p$}}%
\put(15.6100,-5.2500){\makebox(0,0)[lb]{$q$}}%
\put(7.6000,-12.9900){\makebox(0,0)[lb]{$\gamma^{-1}(p)< \gamma^{-1}(q)$}}%
%
{\color[named]{Black}{%
\special{pn 13}%
\special{pa 2560 562}%
\special{pa 4000 562}%
\special{fp}%
}}%
%
{\color[named]{Black}{%
\special{pn 13}%
\special{ar 3280 382 360 540  5.9614348 6.2831853}%
\special{ar 3280 382 360 540  0.0000000 3.4633432}%
}}%
%
{\color[named]{Black}{%
\special{pn 4}%
\special{sh 1}%
\special{ar 2938 562 26 26 0  6.28318530717959E+0000}%
}}%
%
{\color[named]{Black}{%
\special{pn 4}%
\special{sh 1}%
\special{ar 3622 562 26 26 0  6.28318530717959E+0000}%
}}%
%
{\color[named]{Black}{%
\special{pn 13}%
\special{pa 2740 562}%
\special{pa 2650 472}%
\special{fp}%
\special{pa 2740 562}%
\special{pa 2650 652}%
\special{fp}%
}}%
%
{\color[named]{Black}{%
\special{pn 13}%
\special{pa 3280 922}%
\special{pa 3218 804}%
\special{fp}%
\special{pa 3280 922}%
\special{pa 3172 1002}%
\special{fp}%
}}%
\put(26.4100,-7.7700){\makebox(0,0)[lb]{$\gamma$}}%
\put(27.9400,-3.6300){\makebox(0,0)[lb]{$\alpha$}}%
\put(30.0100,-5.2500){\makebox(0,0)[lb]{$q$}}%
\put(37.2100,-5.2500){\makebox(0,0)[lb]{$p$}}%
\put(29.2000,-12.9900){\makebox(0,0)[lb]{$\gamma^{-1}(q)< \gamma^{-1}(p)$}}%
\end{picture}%
\end{center}
\end{figure}

\subsection{The case $*_0=*_1$}
\label{*_1=*_2}

Fix $* \in \partial S$. We shall give a structure of an
involutive right $\mathbb{Q}\hat{\pi}^{\prime}(S)$-bimodule
on the $\mathbb{Q}$-vector space $\mathbb{Q}\pi_1(S,*)=\mathbb{Q}\Pi S(*,*)$.

\textit{Definition of $\sigma$}.
The $\mathbb{Q}$-linear map $\sigma \colon \mathbb{Q}\hat{\pi}^{\prime}(S)
\otimes \mathbb{Q}\pi_1(S,*) \to \mathbb{Q}\pi_1(S,*)$
is defined by setting $*_0=*_1=*$ and applying the formula (\ref{sigma}).

\begin{figure}
\label{fig:secondbasepoint}
\begin{center}
\caption{$*$, $\bullet$, and $\nu$}

\vspace{0.2cm}

\unitlength 0.1in
\begin{picture}( 15.0400, 14.6200)( 14.5000,-18.6200)
%
{\color[named]{Black}{%
\special{pn 8}%
\special{ar 1530 840 80 80  6.2831853 6.2831853}%
\special{ar 1530 840 80 80  0.0000000 4.7123890}%
}}%
%
{\color[named]{Black}{%
\special{pn 8}%
\special{pa 1610 840}%
\special{pa 1562 864}%
\special{fp}%
\special{pa 1610 840}%
\special{pa 1638 884}%
\special{fp}%
}}%
%
{\color[named]{Black}{%
\special{pn 13}%
\special{ar 2954 1300 180 900  2.4668517 3.9269908}%
}}%
%
{\color[named]{Black}{%
\special{pn 8}%
\special{ar 2954 1300 180 900  3.9269908 3.9492130}%
\special{ar 2954 1300 180 900  4.0158797 4.0381019}%
\special{ar 2954 1300 180 900  4.1047686 4.1269908}%
\special{ar 2954 1300 180 900  4.1936575 4.2158797}%
\special{ar 2954 1300 180 900  4.2825464 4.3047686}%
\special{ar 2954 1300 180 900  4.3714353 4.3936575}%
\special{ar 2954 1300 180 900  4.4603242 4.4825464}%
\special{ar 2954 1300 180 900  4.5492130 4.5714353}%
\special{ar 2954 1300 180 900  4.6381019 4.6603242}%
}}%
\put(18.5000,-15.6000){\makebox(0,0)[lb]{$S$}}%
%
{\color[named]{Black}{%
\special{pn 4}%
\special{sh 1}%
\special{ar 2774 1120 26 26 0  6.28318530717959E+0000}%
}}%
%
{\color[named]{Black}{%
\special{pn 4}%
\special{sh 1}%
\special{ar 2774 1390 26 26 0  6.28318530717959E+0000}%
}}%
\put(28.2800,-11.6500){\makebox(0,0)[lb]{$*$}}%
\put(28.2800,-14.3500){\makebox(0,0)[lb]{$\bullet$}}%
%
{\color[named]{Black}{%
\special{pn 13}%
\special{pa 2774 1238}%
\special{pa 2676 1300}%
\special{fp}%
\special{pa 2774 1238}%
\special{pa 2856 1310}%
\special{fp}%
}}%
\put(25.4000,-13.4500){\makebox(0,0)[lb]{$\nu$}}%
\put(28.7300,-18.6700){\makebox(0,0)[lb]{$\partial S$}}%
\end{picture}%

\end{center}
\end{figure}
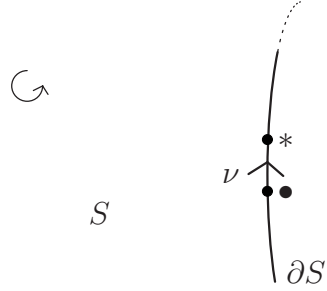

\textit{Definition of $\mu$}.
We regard that the orientation of $\partial S$ is induced from that
of $S$. Pick an orientation preserving embedding $\nu\colon [0,1]\to \partial S$
such that $\nu(1)=*$, and set $\nu(0)=\bullet$. See Figure 7.
Notice that we have an isomorphism $\nu \colon \mathbb{Q}\pi_1(S,*)=\mathbb{Q}\Pi S(*,*)
\cong \mathbb{Q}\Pi S(\bullet, *), u\mapsto \nu u$.
We define $\mu\colon \mathbb{Q}\pi_1(S,*)\to \mathbb{Q}\pi_1(S,*) \otimes
\mathbb{Q}\hat{\pi}^{\prime}(S)$ to be $(\nu^{-1}\otimes 1_{\mathbb{Q}\hat{\pi}^{\prime}(S)})
\circ \mu \circ \nu$. Namely, if $*_0=*_1=*$ we define $\mu$ so that
the following diagram commutes:
\begin{equation}
\label{mu_+,mu_-}
\begin{CD}
\mathbb{Q}\pi_1(S,*) @>{\mu}>> \mathbb{Q}\pi_1(S,*) \otimes
\mathbb{Q}\hat{\pi}^{\prime}(S)\\
@V{\nu}VV @VV{\nu \otimes 1_{\mathbb{Q}\hat{\pi}^{\prime}(S)}}V\\
\mathbb{Q}\Pi S(\bullet,*)
@>{\mu}
>> \mathbb{Q}\Pi S(\bullet,*) \otimes \mathbb{Q}\hat{\pi}^{\prime}(S).
\end{CD}
\end{equation}

Let $\gamma \colon [0,1]\to S$ be an immersed path with
$\gamma(0)=\gamma(1)=*$ such that its self intersections consist
of transverse double points and the velocity vectors
$\dot{\gamma}(0)$ and $\dot{\gamma}(1)$ are linearly
independent. Let $\Gamma \subset {\rm Int}(S)$ be the set of double points of $\gamma$ except $*$.
Then we have
\begin{equation}
\label{mu_+}
\mu(\gamma)=\begin{cases}
-\displaystyle\sum_{p\in \Gamma}\varepsilon(\dot{\gamma}(t_1^p),
\dot{\gamma}(t_2^p))
( \gamma_{0 t_1^p} \gamma_{t_2^p 1})
\otimes |\gamma_{t_1^p t_2^p}|^{\prime},
\quad {\rm if\ }\varepsilon(\dot{\gamma}(0),\dot{\gamma}(1))=+1 \\
1\otimes |\gamma|^{\prime}-\displaystyle\sum_{p\in \Gamma}\varepsilon(\dot{\gamma}(t_1^p),
\dot{\gamma}(t_2^p))
( \gamma_{0 t_1^p} \gamma_{t_2^p 1})
\otimes |\gamma_{t_1^p t_2^p}|^{\prime},
\quad {\rm if\ } \varepsilon(\dot{\gamma}(0),\dot{\gamma}(1))=-1.
\end{cases}
\end{equation}

\begin{rem}
\label{rem:Tu}
{\rm
Our construction is inspired by Turaev's self intersection
$\mu=\mu^T \colon \pi_1(S,*) \to \mathbb{Z}\pi_1(S,*)$ introduced
in \cite{T1} \S 1.4. Actually, for any $\gamma \in \pi_1(S,*)$,
we have
\begin{equation}
\mu^T(\gamma)\gamma=-(1_{\mathbb{Q}\pi_1(S,*)}
\otimes \varepsilon)\mu(\gamma).
\end{equation}
Here $\varepsilon \colon \mathbb{Q}\hat{\pi}^{\prime}(S)
\to \mathbb{Q}$ is the $\mathbb{Q}$-linear map given by
$\varepsilon(\alpha)=1$ for $\alpha \in \hat{\pi}^{\prime}(S)$.
}
\end{rem}

\begin{rem}{\rm
To define $\mu$ for the case $*_0=*_1=*$,
we have moved the start point of paths slightly
along the {\it negatively} oriented boundary of $S$. It is also possible
to move the start point slightly along the {\it positively} oriented
boundary of $S$, and to define another operation which is similar to but different from
what we have defined. If we denote by $\mu_-$ the former operation and
$\mu_+$ the latter, we have
$$\mu_-(\gamma)-\mu_+(\gamma)=-1\otimes |\gamma|^{\prime}$$
for any $\gamma\in \mathbb{Q}\pi_1(S,*)$.
Our choice of convention matches with that of the homotopy intersection form
in \cite{MT}. See also \S \ref{sec:GQT}.
}
\end{rem}

\begin{prop}
If $*_0=*_1=*$, the pair $(\sigma,\mu)$ defined above gives
an involutive right $\mathbb{Q}\hat{\pi}^{\prime}(S)$-bimodule
structure on the $\mathbb{Q}$-vector space $\mathbb{Q}\pi_1(S,*)$.
\end{prop}

\begin{proof}
By the commutativity of the diagram (\ref{mu_+,mu_-}) and
Proposition \ref{comod}, it follows that $\mu$ defines a right
$\mathbb{Q}\hat{\pi}^{\prime}(S)$-comodule
structure on $\mathbb{Q}\pi_1(S,*)$.
Since $\nu$ is compatible with $\sigma$, i.e., the action of $\mathbb{Q}\hat{\pi}^{\prime}(S)$,
Proposition \ref{GTbim} implies the compatibility and the involutivity
of $(\sigma,\mu)$ for $*_0=*_1=*$.
\end{proof}

\section{Completion of the Turaev cobracket}
\label{CTc}

The $\mathbb{Q}$-vector space $\mathbb{Q}\hat{\pi}(S)$ has a natural
decreasing filtration and we can consider the completion
$\widehat{\mathbb{Q}\hat{\pi}}(S)$. As is shown in \cite{KKp} \S 4,
the Goldman bracket induces a (complete) Lie bracket on $\widehat{\mathbb{Q}\hat{\pi}}(S)$. In this section
we show that $\mu$ is compatible with the filtrations of $\mathbb{Q}\Pi S(*_0,*_1)$ and
$\mathbb{Q}\hat{\pi}^{\prime}(S)$, and also show that the Turaev cobracket extends to a complete Lie cobracket
$\delta \colon \widehat{\mathbb{Q}\hat{\pi}}(S) \to
\widehat{\mathbb{Q}\hat{\pi}}(S) \widehat{\otimes} \widehat{\mathbb{Q}\hat{\pi}}(S)$.

\subsection{Completion of the Goldman Lie algebra}
\label{CGL}
We make a few remarks on filtered vector spaces.
Let $V=F_0V \supset F_1V \supset \cdots$ be a filtered $\mathbb{Q}$-vector space.
The projective limit $\widehat{V}:=\varprojlim_n V/F_nV$ is again
a filtered $\mathbb{Q}$-vector space with the filter $F_n\widehat{V}
:={\rm Ker}(\widehat{V} \to V/F_nV)$. We say $V$ is complete if the
natural map $V \to \widehat{V}$ is isomorphic.
If $V$ and $W$ are filtered $\mathbb{Q}$-vector
spaces, the tensor product $V\otimes W$ is naturally filtered by
$F_n(V\otimes W)=\sum_{p+q=n} F_pV \otimes F_qW$.
The complete tensor product $V\widehat{\otimes} W$ is defined as
$V\widehat{\otimes} W:=\widehat{V\otimes W}=\varprojlim_n V\otimes W/F_n(V\otimes W)$.
Note that we have a natural isomorphism
$\widehat{V} \widehat{\otimes} \widehat{W} \cong V\widehat{\otimes} W$.

\begin{dfn}
A complete Lie algebra is a pair $(V,\nabla)$, where
$V$ is a complete filtered $\mathbb{Q}$-vector space
and $\nabla\colon V\widehat{\otimes} V \to V$ is
a $\mathbb{Q}$-linear map continuous
with respect to the topologies coming from the filtrations,
and satisfies the skew condition
$\nabla T=-\nabla \colon V \widehat{\otimes} V\to V$
and the Jacobi identity $\nabla(\nabla \widehat{\otimes} 1)N=0
\colon V \widehat{\otimes} V \widehat{\otimes} V \to V$. Here,
$T\colon V\widehat{\otimes} V\to V\widehat{\otimes} V$ is
that induced from the switch map $T\colon V \otimes V \to V\otimes V$, etc.
We call $\nabla$ a complete Lie bracket.
Similarly, we define a complete Lie coalgebra, bialgebra,
and a complete $V$-module, comodule, and bimodule.
\end{dfn}

Let $S$ be a connected oriented surface.
Take some base point $*\in S$ and set
$$\mathbb{Q}\hat{\pi}(S)(n):=| \mathbb{Q}1 + (I\pi_1(S,*))^n |
\subset \mathbb{Q}\hat{\pi}(S), \quad {\rm for\ }n\ge 0.$$
Here, $I\pi_1(S,*):={\rm Ker}(\mathbb{Q}\pi_1(S,*) \to \mathbb{Q},
\pi \ni x \mapsto 1)$ is the augmentation ideal of the group ring $\mathbb{Q}\pi_1(S,*)$.
We regard $(I\pi_1(S,*))^0=\mathbb{Q}\pi_1(S,*)$.
The space $\mathbb{Q}\hat{\pi}(S)(n)$ is independent of the choice
of $*$. Moreover, the Goldman bracket satisfies
\begin{equation}
\label{Gbrafil}
[\mathbb{Q}\hat{\pi}(S)(n_1),\mathbb{Q}\hat{\pi}(S)(n_2)]
\subset \mathbb{Q}\hat{\pi}(S)(n_1+n_2-2), \quad {\rm for\ }
n_1,n_2\ge 1
\end{equation}
(see \cite{KKp} \S 4.1). This implies that
the Goldman bracket induces a complete Lie bracket on the projective limit
$$\widehat{\mathbb{Q}\hat{\pi}}(S):=\varprojlim_n
\mathbb{Q}\hat{\pi}(S)/\mathbb{Q}\hat{\pi}(S)(n).$$
We call $\widehat{\mathbb{Q}\hat{\pi}}(S)$ the {\it completed
Goldman Lie algebra} of $S$. We denote
$$\widehat{\mathbb{Q}\hat{\pi}}(S)(n):=F_n\widehat{\mathbb{Q}\hat{\pi}}(S)
={\rm Ker}(\widehat{\mathbb{Q}\hat{\pi}}(S)
\to \mathbb{Q}\hat{\pi}(S)/\mathbb{Q}\hat{\pi}(S)(n)), \quad n\ge 0.$$
For $n\ge 0$, let
$$\mathbb{Q}\hat{\pi}^{\prime}(S)(n):=\varpi(\mathbb{Q}\hat{\pi}(S)(n))
=|(I\pi_1(S,*))^n|^{\prime} \subset \mathbb{Q}\hat{\pi}^{\prime}(S).$$
Since $|1|^{\prime}=0$, $\mathbb{Q}\hat{\pi}^{\prime}(S)(0)=
\mathbb{Q}\hat{\pi}^{\prime}(S)(1)=\mathbb{Q}\hat{\pi}^{\prime}(S)$.
The natural map
$\mathbb{Q}\hat{\pi}(S)/\mathbb{Q}\hat{\pi}(S)(n)
\to \mathbb{Q}\hat{\pi}^{\prime}(S)/\mathbb{Q}\hat{\pi}^{\prime}(S)(n)$
is a $\mathbb{Q}$-linear isomorphism for any $n$. Hence
$\widehat{\mathbb{Q}\hat{\pi}}(S)$ is also written as
\begin{equation}
\label{betsudes}
\widehat{\mathbb{Q}\hat{\pi}}(S)=\varprojlim_n
\mathbb{Q}\hat{\pi}^{\prime}(S)/\mathbb{Q}\hat{\pi}^{\prime}(S)(n).
\end{equation}

Let $*_0,*_1\in \partial S$. We make
$\mathbb{Q}\Pi S(*_0,*_1)$ filtered by taking some path $\gamma\in \Pi S(*_0,*_1)$
and setting
\begin{equation}
\label{F_nPi}
F_n\mathbb{Q}\Pi S(*_0,*_1):=\gamma (I\pi_1(S,*_1))^n,\quad {\rm for\ }n\ge 0.
\end{equation}
Then this is independent of the choice of $\gamma$ (see \cite{KKp} Proposition 2.1.1).
In particular we can consider the completion $\widehat{\mathbb{Q}\Pi S}(*_0,*_1)$.
By \cite{KKp} \S 4.1, we see that $\sigma$ induces a $\mathbb{Q}$-linear map
$\sigma \colon \widehat{\mathbb{Q}\hat{\pi}}(S) \widehat{\otimes}
\widehat{\mathbb{Q}\Pi S}(*_0,*_1) \to \widehat{\mathbb{Q}\Pi S}(*_0,*_1)$,
and the complete vector space $\widehat{\mathbb{Q}\Pi S}(*_0,*_1)$
has a structure of a complete right $\widehat{\mathbb{Q}\hat{\pi}}(S)$-module.
As a special case, the completed group ring
$\widehat{\mathbb{Q}\pi_1(S,*)}:=\varprojlim_n \mathbb{Q}\pi_1(S,*)/(I\pi_1(S,*))^n$,
where $*\in \partial S$, has a structure of a complete right
$\widehat{\mathbb{Q}\hat{\pi}}(S)$-module.

\subsection{Intersection of paths}
\label{Intpath}
Take points $*_1,*_2,*_3,*_4\in \partial S$. We introduce a $\mathbb{Q}$-linear map
$$\kappa \colon \mathbb{Q}\Pi S(*_1,*_2)\otimes \mathbb{Q}\Pi S(*_3,*_4)
\to \mathbb{Q}\Pi S(*_1,*_4)\otimes \mathbb{Q}\Pi S(*_3,*_2)$$
using the intersections of two based paths.

{\it The generic case}.
First we consider the case $\{ *_1,*_2\}\cap \{*_3,*_4\}=\emptyset$.
Let $x,y \colon [0,1] \to S$ be immersed paths such that
$x(0)=*_1$, $x(1)=*_2$, $y(0)=*_3$, $y(1)=*_4$ and their
intersections consist of transverse double points. Set
\begin{equation}
\label{k(x,y)}
\kappa(x,y):=-\sum_{p\in x\cap y} \varepsilon(p;x,y)
(x_{*_1 p}y_{p *_4}) \otimes (y_{*_3 p}x_{p *_2}) \in
\mathbb{Q}\Pi S(*_1,*_4) \otimes \mathbb{Q}\Pi S(*_3,*_2).
\end{equation}
By an argument similar to the proof of Proposition \ref{comod},
we see that (\ref{k(x,y)}) gives rise to a well-defined $\mathbb{Q}$-linear map
$\kappa\colon
\mathbb{Q}\Pi S(*_1,*_2) \otimes \mathbb{Q}\Pi S(*_3,*_4)
\to \mathbb{Q}\Pi S(*_1,*_4) \otimes \mathbb{Q}\Pi S(*_3,*_2)$.

{\it The degenerate case}.
Next we consider the case $\{ *_1,*_2\} \cap \{*_3,*_4\}\neq \emptyset$.
The idea is, as was the case of $\mu$, to move the points $*_1,*_2$ slightly along the negatively oriented
boundary of $S$ to achieve $\{*_1,*_2\} \cap \{*_3,*_4\}=\emptyset$, then to
apply the formula (\ref{k(x,y)}). Using two examples we explain this procedure more precisely.

For the first example, suppose $*_1,*_2,*_4$ are distinct and $*_3=*_2$. Set $*=*_2$ and
let $\bullet$ and $\nu$ be as in Figure 7. We assume that the image of $\nu$ is
so small that it does not contain $*_1$ and $*_4$. Notice that we have isomorphisms
$\nu\colon \mathbb{Q}\Pi S(*_1,*) \to \mathbb{Q}\Pi S(*_1,\bullet), x\mapsto x\overline{\nu}$,
and $\nu\colon \mathbb{Q}\pi_1(S,*)=\mathbb{Q}\Pi S(*,*)\cong \mathbb{Q}\Pi S(*,\bullet),
x\mapsto x\overline{\nu}$. Here and throughout this paper
$\overline{\nu}\in \Pi S(*,\bullet)$ is the inverse of $\nu$.
We define $\kappa\colon \mathbb{Q}\Pi (*_1,*)\otimes \mathbb{Q}\Pi (*,*_4)
\to \mathbb{Q}\Pi (*_1,*_4)\otimes \mathbb{Q}\pi_1(S,*)$ so that the following diagram commetes:
$$
\begin{CD}
\mathbb{Q}\Pi S(*_1,*) \otimes \mathbb{Q}\Pi S(*,*_4) @>{\kappa}>>
\mathbb{Q}\Pi S(*_1,*_4) \otimes \mathbb{Q}\pi_1(S,*) \\
@V{\nu \otimes 1_{\mathbb{Q}\Pi S(*,*_4)}}VV @VV{1_{\mathbb{Q}\Pi S(*_1,*_4)} \otimes \nu}V\\
\mathbb{Q}\Pi S(*_1,\bullet) \otimes \mathbb{Q}\Pi S(*,*_4)
@>{\kappa}
>> \mathbb{Q}\Pi S(*_1,*_4) \otimes \mathbb{Q}\Pi S(*,\bullet).
\end{CD}
$$

For the second example, suppose $*_1=*_2=*_3=*_4$. (This is the most extreme case.)
Set $*=*_1$ and let $\bullet$ and $\nu$ be as in Figure 7. Notice that we have
three isomorphisms $\mathbb{Q}\pi_1(S,*)\cong \mathbb{Q}\pi_1(S,\bullet),x\mapsto \nu x \overline{\nu}$,
$\mathbb{Q}\pi_1(S,*)\cong \mathbb{Q}\Pi S(\bullet,*), x\mapsto \nu x$, and
$\mathbb{Q}\pi_1(S,*)\cong \mathbb{Q}\Pi S(*,\bullet), x\mapsto x\overline{\nu}$,
for all of which we use the letter $\nu$. We define $\kappa\colon
\mathbb{Q}\pi_1(S,*)\otimes \mathbb{Q}\pi_1(S,*) \to
\mathbb{Q}\pi_1(S,*)\otimes \mathbb{Q}\pi_1(S,*)$ so that the following diagram commutes:
$$
\begin{CD}
\mathbb{Q}\pi_1(S,*) \otimes \mathbb{Q}\pi_1(S,*) @>{\kappa}>>
\mathbb{Q}\pi_1(S,*) \otimes \mathbb{Q}\pi_1(S,*) \\
@V{\nu \otimes 1_{\mathbb{Q}\pi_1(S,*)}}VV @VV{\nu \otimes \nu}V\\
\mathbb{Q}\pi_1(S,\bullet) \otimes \mathbb{Q}\pi_1(S,*)
@>{\kappa}
>> \mathbb{Q}\Pi S(\bullet,*) \otimes \mathbb{Q}\Pi S(*,\bullet).
\end{CD}
$$

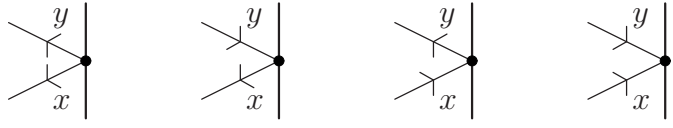
\begin{figure}
\begin{center}
\caption{$\varepsilon(\dot{x}(\delta_1),\dot{y}(\delta_2))=(-1)^{\delta_1+\delta_2+1}$}

\vspace{0.2cm}

\unitlength 0.1in
\begin{picture}( 34.0000,  6.0000)(  6.0000,-10.0000)
%
{\color[named]{Black}{%
\special{pn 13}%
\special{pa 1000 400}%
\special{pa 1000 1000}%
\special{fp}%
}}%
%
{\color[named]{Black}{%
\special{pn 8}%
\special{pa 1000 700}%
\special{pa 600 900}%
\special{fp}%
}}%
%
{\color[named]{Black}{%
\special{pn 8}%
\special{pa 1000 700}%
\special{pa 600 500}%
\special{fp}%
}}%
%
{\color[named]{Black}{%
\special{pn 4}%
\special{sh 1}%
\special{ar 1000 700 26 26 0  6.28318530717959E+0000}%
}}%
%
{\color[named]{Black}{%
\special{pn 8}%
\special{pa 800 600}%
\special{pa 800 680}%
\special{fp}%
\special{pa 800 600}%
\special{pa 870 560}%
\special{fp}%
}}%
%
{\color[named]{Black}{%
\special{pn 8}%
\special{pa 800 800}%
\special{pa 800 720}%
\special{fp}%
\special{pa 800 800}%
\special{pa 870 840}%
\special{fp}%
}}%
%
{\color[named]{Black}{%
\special{pn 8}%
\special{pa 1800 600}%
\special{pa 1800 520}%
\special{fp}%
\special{pa 1800 600}%
\special{pa 1730 640}%
\special{fp}%
}}%
\put(8.3000,-9.5000){\makebox(0,0)[lb]{$x$}}%
\put(8.3000,-5.3000){\makebox(0,0)[lb]{$y$}}%
%
{\color[named]{Black}{%
\special{pn 13}%
\special{pa 2000 400}%
\special{pa 2000 1000}%
\special{fp}%
}}%
%
{\color[named]{Black}{%
\special{pn 8}%
\special{pa 2000 700}%
\special{pa 1600 900}%
\special{fp}%
}}%
%
{\color[named]{Black}{%
\special{pn 8}%
\special{pa 2000 700}%
\special{pa 1600 500}%
\special{fp}%
}}%
%
{\color[named]{Black}{%
\special{pn 4}%
\special{sh 1}%
\special{ar 2000 700 26 26 0  6.28318530717959E+0000}%
}}%
\put(18.3000,-9.5000){\makebox(0,0)[lb]{$x$}}%
\put(18.3000,-5.3000){\makebox(0,0)[lb]{$y$}}%
%
{\color[named]{Black}{%
\special{pn 13}%
\special{pa 3000 400}%
\special{pa 3000 1000}%
\special{fp}%
}}%
%
{\color[named]{Black}{%
\special{pn 8}%
\special{pa 3000 700}%
\special{pa 2600 900}%
\special{fp}%
}}%
%
{\color[named]{Black}{%
\special{pn 8}%
\special{pa 3000 700}%
\special{pa 2600 500}%
\special{fp}%
}}%
%
{\color[named]{Black}{%
\special{pn 4}%
\special{sh 1}%
\special{ar 3000 700 26 26 0  6.28318530717959E+0000}%
}}%
%
{\color[named]{Black}{%
\special{pn 8}%
\special{pa 2800 800}%
\special{pa 2800 880}%
\special{fp}%
\special{pa 2800 800}%
\special{pa 2730 760}%
\special{fp}%
}}%
%
{\color[named]{Black}{%
\special{pn 8}%
\special{pa 2800 600}%
\special{pa 2800 680}%
\special{fp}%
\special{pa 2800 600}%
\special{pa 2870 560}%
\special{fp}%
}}%
\put(28.3000,-9.5000){\makebox(0,0)[lb]{$x$}}%
\put(28.3000,-5.3000){\makebox(0,0)[lb]{$y$}}%
%
{\color[named]{Black}{%
\special{pn 13}%
\special{pa 4000 400}%
\special{pa 4000 1000}%
\special{fp}%
}}%
%
{\color[named]{Black}{%
\special{pn 8}%
\special{pa 4000 700}%
\special{pa 3600 900}%
\special{fp}%
}}%
%
{\color[named]{Black}{%
\special{pn 8}%
\special{pa 4000 700}%
\special{pa 3600 500}%
\special{fp}%
}}%
%
{\color[named]{Black}{%
\special{pn 4}%
\special{sh 1}%
\special{ar 4000 700 26 26 0  6.28318530717959E+0000}%
}}%
%
{\color[named]{Black}{%
\special{pn 8}%
\special{pa 3800 800}%
\special{pa 3800 880}%
\special{fp}%
\special{pa 3800 800}%
\special{pa 3730 760}%
\special{fp}%
}}%
\put(38.3000,-9.5000){\makebox(0,0)[lb]{$x$}}%
\put(38.3000,-5.3000){\makebox(0,0)[lb]{$y$}}%
%
{\color[named]{Black}{%
\special{pn 8}%
\special{pa 1800 800}%
\special{pa 1800 720}%
\special{fp}%
\special{pa 1800 800}%
\special{pa 1870 840}%
\special{fp}%
}}%
%
{\color[named]{Black}{%
\special{pn 8}%
\special{pa 3800 600}%
\special{pa 3800 520}%
\special{fp}%
\special{pa 3800 600}%
\special{pa 3730 640}%
\special{fp}%
}}%
\end{picture}%

\end{center}
\end{figure}

An explicit formula for $\kappa$ is given as follows.
Let $x,y\colon [0,1]\to S$ be immersed paths with $x(0)=*_1$, $x(1)=*_2$,
$y(0)=*_3$, $y(1)=*_4$, such that their intersections in the interior of $S$ consist of transverse double points.
Furthermore, if $x(\delta_1)=y(\delta_2)$ for some $\delta_1,\delta_2\in \{ 0,1\}$,
we assume that $\varepsilon(\dot{x}(\delta_1),\dot{y}(\delta_2))=(-1)^{\delta_1+\delta_2+1}$.
See Figure 8. Note that this condition is always satisfied by a suitable homotopy.
Then we have

\begin{equation}
\label{k_+}
\kappa(x,y)=
-\displaystyle\sum_{p\in x\cap y \setminus \partial S}\varepsilon(p;x,y)
(x_{*_1 p}y_{p *_4}) \otimes (y_{*_3 p}x_{p *_2}).
\end{equation}

\begin{rem}
{\rm To define $\kappa$, we have taken $*_i$ from $\partial S$.
By the same reason as is explained in Remark \ref{rem-*},
the formula (\ref{k(x,y)}) does not work if at least one of $*_i$ lies in ${\rm Int}(S)$.
}
\end{rem}

\begin{rem}
{\rm The map $\kappa$ is inspired by
Turaev's intersection $\lambda \colon \mathbb{Z}\pi_1(S,*)
\otimes \mathbb{Z}\pi_1(S,*)\to \mathbb{Z}\pi_1(S,*)$ introduced
in \cite{T1} \S 1.4. Fix $*\in \partial S$ and consider the $\mathbb{Q}$-linear map
$\kappa\colon \mathbb{Q}\pi_1(S,*) \otimes \mathbb{Q}\pi_1(S,*)
\to \mathbb{Q}\pi_1(S,*) \otimes \mathbb{Q}\pi_1(S,*)$
introduced above. Then for any $x,y\in \pi_1(S,*)$, we have
\begin{equation}
\lambda(x,y)=-(1_{\mathbb{Q}\pi_1(S,*)} \otimes \varepsilon)\kappa(x,y^{-1}).
\end{equation}
Here, $\varepsilon$ is the augmentation map of the group ring $\mathbb{Q}\pi_1(S,*)$.
}
\end{rem}

\subsection{Product formulas}

We give formulas for how $\kappa$ and $\mu$ behaves under the conjunction of paths.

\begin{lem}
\label{k(xy,z)}
Let $*_1,*_2,*_2^{\prime},*_3,*_4,*_4^{\prime} \in \partial S$.

\begin{enumerate}
\item For any $x\in \mathbb{Q}\Pi S(*_1,*_2)$, $y\in \mathbb{Q}\Pi S(*_2,*_2^{\prime})$,
and $z\in \mathbb{Q}\Pi S(*_3,*_4)$, we have
$$\kappa(xy,z) = \kappa(x,z)(1\otimes y)+(x\otimes 1)\kappa(y,z)
\in \mathbb{Q}\Pi S(*_1,*_4) \otimes \mathbb{Q}\Pi S(*_3,*_2^{\prime}).$$
\item For any $x\in \mathbb{Q}\Pi S(*_1,*_2)$, $y\in \mathbb{Q}\Pi S(*_3,*_4)$,
and $z\in \mathbb{Q}\Pi S(*_4,*_4^{\prime})$, we have
$$\kappa(x,yz) = \kappa(x,y)(z\otimes 1)+(1\otimes y)\kappa(x,z)
\in \mathbb{Q}\Pi S(*_1,*_4^{\prime}) \otimes \mathbb{Q}\Pi S(*_3,*_2).$$
\end{enumerate}
Here, $\kappa(x,z)(1\otimes y)$ means the image of $\kappa(x,z)\otimes y$
by the map $\mathbb{Q}\Pi S(*_1,*_4)\otimes \mathbb{Q}\Pi S(*_3,*_2)
\otimes \mathbb{Q}\Pi S(*_2,*_2^{\prime}) \to \mathbb{Q}\Pi S(*_1,*_4)
\otimes \mathbb{Q}\Pi S(*_3,*_2^{\prime})$, $u\otimes v \otimes w \mapsto u\otimes vw$, etc.
\end{lem}

\begin{proof}
We only prove the formula $\kappa(xy,z)=\kappa(x,z)(1\otimes y)
+(x\otimes 1)\kappa(y,z)$. The other formula
is proved similarly. Let $x,y,z \colon [0,1] \to S$ be immersed paths with
$x(0)=*_1$, $x(1)=*_2$, $y(0)=*_2$, $y(1)=*_2^{\prime}$, $z(0)=*_3$, $z(1)=*_4$,
such that their intersections in the interior of $S$ consist
of transverse double points. Moreover, we assume that
$\varepsilon(\dot{x}(\delta_1),\dot{z}(\delta_2))=(-1)^{\delta_1+\delta_2+1}$
(resp.\ $\varepsilon(\dot{y}(\delta_1),\dot{z}(\delta_2))=(-1)^{\delta_1+\delta_2+1}$)
for any $\delta_1,\delta_2\in \{ 0,1\}$ with $x(\delta_1)=z(\delta_2)$ (resp.\ $y(\delta_1)=z(\delta_2)$).

Applying the formula (\ref{k_+}), we compute
\begin{eqnarray*}
\kappa(xy,z) &= &-\sum_{p\in (xy)\cap z \setminus \partial S}
\varepsilon(p;xy,z) (xy)_{*_1 p}z_{p *_4} \otimes z_{*_3 p}(xy)_{p *_3} \\
&=& -\sum_{p\in x\cap z \setminus \partial S} \varepsilon(p;x,z) (x_{*_1 p}z_{p *_4}) \otimes
(z_{*_3 p}x_{p *_2} y) \\
& & -\sum_{p\in y\cap z \setminus \partial S}\varepsilon(p;y,z)
(x y_{*_2 p}z_{p *_4}) \otimes (z_{*_3 p}y_{p *_3}).
\end{eqnarray*}
The first and the second terms are equal to $\kappa(x,z)(1\otimes y)$ and
$(x\otimes 1)\kappa(y,z)$, respectively.
\end{proof}

\begin{cor}
\label{k(x_1-x_n)}
Let $n\ge 2$ and let $*_1,\ldots,*_{n+2}\in \partial S$ be $n+2$ points.
For any $x_1,\ldots,x_{n+1}$, $x_i \in \mathbb{Q}\Pi S(*_i,*_{i+1})$,
we have
\begin{eqnarray*}
\kappa(x_1\cdot \cdots \cdot x_n,x_{n+1}) &=&
\sum_{i=1}^n ((x_1\cdots x_{i-1}) \otimes 1)\kappa(x_i,x_{n+1})(1\otimes
(x_{i+1}\cdots x_n)) \\
\kappa(x_1,x_2\cdot \cdots \cdot x_nx_{n+1}) &=&
\sum_{i=2}^{n+1} (1\otimes (x_2\cdots x_{i-1}))\kappa(x_1,x_i)
((x_{i+1}\cdots x_{n+1})\otimes 1).
\end{eqnarray*}
\end{cor}

\begin{proof}
Note that $\kappa(x_i,x_{n+2})(1\otimes (x_{i+1}\cdots x_n))(1\otimes x_{n+1})
=\kappa(x_i,x_{n+2})(1\otimes (x_{i+1}\cdots x_nx_{n+1}))$, etc.
The assertion follows from Lemma \ref{k(xy,z)} and induction on $n$.
\end{proof}

\begin{lem}
\label{mu(xy)}
Let $*_1,*_2,*_3 \in \partial S$ be three points. For
any $x\in \mathbb{Q}\Pi S(*_1,*_2)$ and $y\in \mathbb{Q}\Pi S(*_2,*_3)$, we have
$$\mu(xy) = \mu(x)(y\otimes 1)+(x\otimes 1)\mu(y)
+(1_{\mathbb{Q}\Pi S(*_1,*_3)} \otimes |\ |^{\prime})\kappa(x,y)
\in \mathbb{Q}\Pi S(*_1,*_3) \otimes \mathbb{Q}\hat{\pi}^{\prime}(S).$$
Here, $\mu(x)(y\otimes 1)$ means the image of $\mu(x)\otimes y$
by the map $\mathbb{Q}\Pi S(*_1,*_2) \otimes \mathbb{Q}\hat{\pi}^{\prime}(S)
\otimes \mathbb{Q}\Pi S(*_2,*_3) \to \mathbb{Q}\Pi S(*_1,*_3) \otimes
\mathbb{Q}\hat{\pi}^{\prime}(S)$, $u\otimes v\otimes w \mapsto uw \otimes v$, etc.
\end{lem}

\begin{proof}
Let $x,y \colon [0,1]\to S$ be immersed paths with $x(0)=*_1$,
$x(1)=y(0)=*_2$, $y(1)=*_3$, such that their intersections and self intersections
in the interior of $S$ consist of transverse double points. Moreover, we assume that
$\varepsilon(\dot{x}(\delta_1),\dot{y}(\delta_2))=(-1)^{\delta_1+\delta_2+1}$
for any $\delta_1,\delta_2\in \{ 0,1\}$ with $x(\delta_1)=y(\delta_2)$,
and $\varepsilon(\dot{x}(0),\dot{x}(1))=+1$ (resp.\ $\varepsilon(\dot{y}(0),\dot{y}(1))=+1$)
if $x(0)=x(1)$ (resp.\ $y(0)=y(1)$). Note that it is always possible to achieve this by a suitable homotopy.
Let $\Gamma_x$ and $\Gamma_y$ be the set of double points of $x$ and $y$, respectively.
Then the set of double points of $xy$ is
$\Gamma_x \cup \Gamma_y \cup (x\cap y \setminus \partial S)$.
We have
\begin{eqnarray*}
\mu(xy) &=& -\sum_{p\in \Gamma_x} \varepsilon(x^{\prime}(t_1^p),x^{\prime}(t_2^p))
(x_{*_1 p}x_{p *_2}y) \otimes |x_{t_1^p t_2^p}|^{\prime} \\
& & -\sum_{p\in \Gamma_y} \varepsilon(y^{\prime}(t_1^p),y^{\prime}(t_2^p))
(xy_{*_2 p}y_{p *_3}) \otimes |y_{t_1^p t_2^p}|^{\prime} \\
& & -\sum_{p\in x\cap y \setminus \partial S}
\varepsilon(x^{\prime}(p),y^{\prime}(p))(x_{*_1 p}y_{p *_3}) \otimes |x_{p *_2}y_{*_2 p}|^{\prime}.
\end{eqnarray*}
The first and the second terms are equal to $\mu(x)(y\otimes 1)$
and $(x\otimes 1)\mu(y)$, respectively. Since $|x_{p *_2}y_{*_2 p}|=|y_{*_2 p}x_{p *_2}|$,
the third term is equal to $(1_{\mathbb{Q}\Pi S(*_1,*_3)} \otimes |\ |^{\prime})\kappa(x,y)$.
Hence $\mu(xy)=\mu(x)(y\otimes 1)+(x\otimes 1)\mu(y)
+(1_{\mathbb{Q}\Pi S(*_1,*_3)} \otimes |\ |^{\prime})\kappa(x,y)$.
\end{proof}

By Corollary \ref{k(x_1-x_n)}, Lemma \ref{mu(xy)}, and
induction on $n$, we have the following.

\begin{cor}
\label{mu(1-n)}
Let $n\ge 2$ and let $*_1,\ldots, *_{n+1} \in \partial S$ be $n+1$ points. For any $x_1,\ldots,x_n$,
$x_i \in \mathbb{Q}\Pi S(*_i,*_{i+1})$, we have
\begin{eqnarray*}
\mu(x_1\cdots x_n) &=&
\sum_{i=1}^n ((x_1\cdots x_{i-1})\otimes 1)\mu(x_i)((x_{i+1}\cdots x_n)
\otimes 1) \\
& & +\sum_{i<j} ((x_1\cdots x_{i-1})\otimes 1)
K_{i,j}((x_{j+1}\cdots x_n)\otimes 1),
\end{eqnarray*}
where $K_{i,j}=(1_{\mathbb{Q}\Pi S(*_i,*_{j+1})}\otimes |\ |^{\prime})
(\kappa(x_i,x_j)(1\otimes (x_{i+1}\cdots x_{j-1})))$.
\end{cor}

\subsection{$\mu$ and the filtrations of $\mathbb{Q}\Pi S(*_0,*_1)$,
$\mathbb{Q}\hat{\pi}^{\prime}(S)$}

We assume that the boundary of $S$ is not empty. Take $*\in \partial S$.

\begin{lem}
\label{mu-delta}
The following diagram is commutative:
$$
\begin{CD}
\mathbb{Q}\pi_1(S,*) @>{\mu}>>
\mathbb{Q}\pi_1(S,*) \otimes \mathbb{Q}\hat{\pi}^{\prime}(S) \\
@V{|\ |^{\prime}}VV 
@VV{(1-T)(|\ |^{\prime} \otimes 1_{\mathbb{Q}\hat{\pi}^{\prime}(S)})}V\\
\mathbb{Q}\hat{\pi}^{\prime}(S)
@>{\delta}
>> \mathbb{Q}\hat{\pi}^{\prime}(S) \otimes \mathbb{Q}\hat{\pi}^{\prime}(S).
\end{CD}
$$
\end{lem}

\begin{proof}
Let $\gamma \colon [0,1]\to S$ be an immersed path with $\gamma(0)=\gamma(1)$
such that its self intersections consist of transverse double points
and $\varepsilon(\dot{\gamma}(0),\dot{\gamma}(1))=+1$. By
(\ref{mu_+}),
$$\mu(\gamma)=-\displaystyle\sum_{p\in \Gamma}
\varepsilon(\dot{\gamma}(t_1^p),\dot{\gamma}(t_2^p))
( \gamma_{0 t_1^p} \gamma_{t_2^p 1})
\otimes |\gamma_{t_1^p t_2^p}|^{\prime}.$$
Using $|\gamma_{0 t_1^p} \gamma_{t_2^p 1}|=|\gamma_{t_2^p 1} \gamma_{0 t_1^p}|$,
we obtain
\begin{eqnarray*}
& & (1-T)(|\ |^{\prime} \otimes 1_{\mathbb{Q}\hat{\pi}^{\prime}(S)})\mu(\gamma) \\
&=& \sum_{p\in \Gamma}\varepsilon(\dot{\gamma}(t_1^p),\dot{\gamma}(t_2^p))
(|\gamma_{t_1^p t_2^p}|^{\prime} \otimes |\gamma_{t_2^p 1} \gamma_{0 t_1^p}|^{\prime}
-|\gamma_{t_2^p 1} \gamma_{0 t_1^p}|^{\prime} \otimes |\gamma_{t_1^p t_2^p}|^{\prime}).
\end{eqnarray*}
This coincides with $\delta(|\gamma|^{\prime})$.
\end{proof}

\begin{prop}
\label{mu-fil}
Let $*_0,*_1\in \partial S$. For $n\ge 2$, we have
$$\mu(F_n\mathbb{Q}\Pi S(*_0,*_1)) \subset
\sum_{p+q=n-2} F_p\mathbb{Q}\Pi S(*_0,*_1) \otimes \mathbb{Q}\hat{\pi}^{\prime}(S)(q)
\subset \mathbb{Q}\Pi S(*_0,*_1) \otimes \mathbb{Q}\hat{\pi}^{\prime}(S).$$
In particular, for $*\in \partial S$ and $n\ge 2$, we have
$$\mu(I\pi_1(S,*)^n)\subset
\sum_{p+q=n-2} I\pi_1(S,*)^p \otimes \mathbb{Q}\hat{\pi}^{\prime}(S)(q).$$
\end{prop}

\begin{proof}
If $*_0=*_1=*$, the assertion follows from Corollary \ref{mu(1-n)} by
setting $*_i=*$ and $x_i\in I\pi_1(S,*)$. In the general case, take
$\gamma\in \Pi S(*_0,*_1)$ and $x_i\in I\pi_1(S,*_1)$, $1\le i\le n$.
Again by Corollary \ref{mu(1-n)}, we have $\mu(\gamma x_1\cdots x_n)
\in \sum_{p+q=n-2} F_p\mathbb{Q}\Pi S(*_0,*_1) \otimes \mathbb{Q}\hat{\pi}^{\prime}(S)(q)$.
\end{proof}

By Lemma \ref{mu-delta} and (\ref{betsudes}), we obtain

\begin{cor}
The Turaev cobracket $\delta \colon \mathbb{Q}\hat{\pi}^{\prime}(S)
\to \mathbb{Q}\hat{\pi}^{\prime}(S)\otimes \mathbb{Q}\hat{\pi}^{\prime}(S)$
satisfies
$$\delta(\mathbb{Q}\hat{\pi}^{\prime}(S)(n))\subset
\sum_{p+q=n-2} \mathbb{Q}\hat{\pi}^{\prime}(S)(p) \otimes
\mathbb{Q}\hat{\pi}^{\prime}(S)(q)$$
for any $n\ge 2$. Moreover, $\delta$ induces a complete Lie cobracket
$\delta \colon \widehat{\mathbb{Q}\hat{\pi}}(S) \to
\widehat{\mathbb{Q}\hat{\pi}}(S) \hat{\otimes}
\widehat{\mathbb{Q}\hat{\pi}}(S)$.
\end{cor}

The $\mathbb{Q}$-vector space $\widehat{\mathbb{Q}\hat{\pi}}(S)$ has a
structure of a involutive complete Lie bialgebra with respect to the complete Lie
bracket in \S \ref{CGL} and the complete Lie cobracket above. We call this
{\it the completed Goldman-Turaev Lie bialgebra}. We define 
$$
\widehat{\mathbb{Q}\Pi S}(*_0,*_1) := \varprojlim_{n\to\infty}
\mathbb{Q}\Pi S(*_0,*_1)/F_n\mathbb{Q}\Pi S(*_0,*_1),
$$
which is a $\widehat{\mathbb{Q}{\hat\pi}}(S)$-module by means of
$\sigma$, as was stated in \S \ref{CGL}.
As another consequence of Proposition \ref{mu-fil}, $\mu$ induces a $\mathbb{Q}$-linear map 
\begin{equation}
\label{mucpleted}
\mu\colon \widehat{\mathbb{Q}\Pi S}(*_0,*_1)\to 
\widehat{\mathbb{Q}\Pi S}(*_0,*_1)
\widehat{\otimes}\widehat{\mathbb{Q}{\hat\pi}}(S),
\end{equation}
which makes $\widehat{\mathbb{Q}\Pi S}(*_0,*_1)$ a complete involutive right 
$\widehat{\mathbb{Q}{\hat\pi}}(S)$-bimodule. In \S \ref{sec:App} we will use this
bimodule structure to prove that some generalized Dehn twists are 
not realized by diffeomorphisms. 

\section{Application to generalized Dehn twists}
\label{sec:App}

In this and the next section we discuss applications of our consideration
of the (self) intersections of curves to the study of the mapping class groups.
In this section we study {\it generalized Dehn twists}, which was introduced in \cite{KKp} \cite{Ku2}. 

\subsection{Generalized Dehn twists}
\label{GDT}

Generalized Dehn twists are associated with not necessarily
simple loops on a surface, and are defined as elements of
a certain enlargement of the mapping class group of the surface.
We recall generalized Dehn twists following \cite{KKp} \S 5.
For another treatment, see \cite{MT}.

Let $S$ be a compact connected oriented surface with non-empty boundary,
or a surface obtained from such a surface by removing finitely many points
in the interior. We denote by $\mathcal{M}(S,\partial S)$ the mapping class group
of the pair $(S,\partial S)$, i.e., the group of orientation preserving
diffeomorphisms of $S$ fixing $\partial S$ pointwise, modulo
isotopies relative to $\partial S$. The group
$\mathcal{M}(S,\partial S)$ naturally acts on each
$\Pi S(p_0,p_1)$, $p_0,p_1\in \partial S$.

Let $E\subset \partial S$ be a subset which
contains at least one point of any connected component of $\partial S$.
Then we can construct a small additive category $\mathbb{Q}\mathcal{C}(S,E)$,
whose set of objects is $E$, and whose set of morphisms from $p_0\in E$
to $p_1 \in E$ is $\mathbb{Q}\Pi S(p_0,p_1)$. As we mentioned in \S \ref{CGL},
$\mathbb{Q}\Pi S(p_0,p_1)$ is filtered and its completion
$\widehat{\mathbb{Q}\Pi S}(p_0,p_1)$ is defined.
Let $\widehat{\mathbb{Q}\mathcal{C}(S,E)}$ be a small additive category
whose set of objects is $E$, and whose set of morphisms from
$p_0\in E$ to $p_1 \in E$ is $\widehat{\mathbb{Q}\Pi S}(p_0,p_1)$.
In \cite{KKp},
$\widehat{\mathbb{Q}\mathcal{C}(S,E)}$ is called the completion of
$\mathbb{Q}\mathcal{C}(S,E)$.

The action of $\mathcal{M}(S,\partial S)$ on $\Pi S(p_0,p_1)$
naturally induces a $\mathbb{Q}$-linear automorphism of $\mathbb{Q}\Pi S(p_0,p_1)$,
as well as a $\mathbb{Q}$-linear automorphism of $\widehat{\mathbb{Q}\Pi S}(p_0,p_1)$.
In this way we obtain a group homomorphism of Dehn-Nielsen type
\begin{equation}
\label{DN}
\widehat{{\sf DN}}\colon \mathcal{M}(S,\partial S)
\to {\rm Aut}(\widehat{\mathbb{Q}\mathcal{C}(S,E)}),
\end{equation}
where ${\rm Aut}(\widehat{\mathbb{Q}\mathcal{C}(S,E)})$ is the group of covariant
functors from $\widehat{\mathbb{Q}\mathcal{C}(S,E)}$ to itself, which act
on the set of objects as the identity, and act on each set of morphisms
as $\mathbb{Q}$-linear automorphisms. This group homomorphism is
injective (see \cite{KKp} Theorem 3.1.1).

Let $C\subset S\setminus \partial S$ be an unoriented loop.
Take $q\in S$ and let $x\in \pi_1(S,q)$ be a based loop which
is homotopic to $C$ as an unoriented loop. The quantity
$$L(C):=\left| \frac{1}{2}(\log x)^2 \right|\in
\widehat{\mathbb{Q}\hat{\pi}}(S)(2),$$
where $|\ |\colon \widehat{\mathbb{Q}\pi_1(S,q)}\to
\widehat{\mathbb{Q}\hat{\pi}}(S)$ is the map induced by $|\ |\colon
\mathbb{Q}\pi_1(S,q) \to \mathbb{Q}\hat{\pi}(S)$, is independent
of the choice of $q$ and $x$.

A family of $\mathbb{Q}$-linear homomorphisms $D=D^{(p_0,p_1)}\colon
\widehat{\mathbb{Q}\Pi S}(p_0,p_1)\to \widehat{\mathbb{Q}\Pi S}(p_0,p_1)$,
$p_0,p_1\in E$, is called a derivation of $\widehat{\mathbb{Q}\mathcal{C}(S,E)}$,
if it satisfies the Leibniz rule
$$D(uv)=(Du)v+u(Dv)$$
for any $p_0,p_1,p_2\in E$, $u\in \widehat{\mathbb{Q}\Pi S}(p_0,p_1)$,
and $v\in \widehat{\mathbb{Q}\Pi S}(p_1,p_2)$. The set of derivations
of $\widehat{\mathbb{Q}\mathcal{C}(S,E)}$ naturally has a structure
of a Lie algebra, which we denote by ${\rm Der}(\widehat{\mathbb{Q}\mathcal{C}(S,E)})$.
Then we obtain a Lie algebra homomorphism
$$\sigma \colon \widehat{\mathbb{Q}\hat{\pi}}(S)
\to {\rm Der}(\widehat{\mathbb{Q}\mathcal{C}(S,E)}),$$
by collecting the structure morphisms
$\sigma \colon \widehat{\mathbb{Q}\hat{\pi}}(S)
\widehat{\otimes} \widehat{\mathbb{Q}\Pi S}(p_0,p_1) \to
\widehat{\mathbb{Q}\Pi S}(p_0,p_1)$, $p_0,p_1 \in E$ (see \cite{KKp} \S 4.1).
For $p_0,p_1\in E$, the exponential of the derivation
$\sigma(L(C)) \in {\rm End}(\widehat{\mathbb{Q}\Pi S}(p_0,p_1))$
converges and we obtain an automorphism
$$\exp(\sigma(L(C))) \in {\rm Aut}(\widehat{\mathbb{Q}\mathcal{C}(S,E)}),$$
which we call the {\it generalized Dehn twist along $C$}
(\cite{KKp} Lemma 5.1.1, Definition 5.3.1).
If $C$ is simple, then this is (the $\widehat{\sf DN}$-image of)
the usual right handed Dehn twist along $C$ (\cite{KKp} Theorem 5.2.1).

\begin{rem}
{\rm 
Actually $\exp(\sigma(L(C)))$ lies in a subgroup
$A(S,E)\subset {\rm Aut}(\widehat{\mathbb{Q}\mathcal{C}(S,E)})$,
which was introduced in \cite{KKp} Definition 3.3.1.
}
\end{rem}

\subsection{A criterion of the realizability}

A natural question is whether $\exp(\sigma(L(C)))$ is realized by a diffeomorphism,
i.e., is in the $\widehat{\sf DN}$-image. In \cite{KKp} \cite{Ku2}
we showed that if $C$ is a figure eight, then $\exp(\sigma(L(C)))$ is not in the
$\widehat{\sf DN}$-image. To extend this result for curves in wider classes,
we consider the self intersections of curves.

Let $C \subset S\setminus\partial S$ be an unoriented free loop, 
and $N \subset S\setminus\partial S$ a connected compact subsurface which
is a neighborhood of $C$, and not diffeomorphic to $D^2$.
If the generalized Dehn twist $\exp(\sigma(L(C)))$ is
the $\widehat{\sf DN}$-image of a mapping class
$\varphi \in \mathcal{M}(S,\partial S)$, the support of
(a representative of) $\varphi$ is included in the subsurface $N$,
by the localization theorem \cite{KKp} Theorem 5.3.3.

Using the fact that $\mu$ maps simple paths to zero and
a diffeomorphism preserves the simplicity of curves,
together with cut and paste techniques developed in \cite{KKp},
we have the following.

\begin{prop}
\label{criterion}
Suppose the inclusion homomorphism $\pi_1(N) \to \pi_1(S)$ is
injective. Assume the generalized Dehn twist $\exp(\sigma(L(C)))$
is realized by a diffeomorphism. Then we have
$$
\mu(\sigma(L(C))(\gamma)) = 0\in \widehat{\mathbb{Q}\Pi N}(*_0,*_1) \widehat{\otimes}
\widehat{\mathbb{Q}\hat{\pi}}(N)
$$
for any distinct points $*_0, *_1 \in \partial N$ and any simple path $\gamma \in
\Pi N(*_0.*_1)$. When $*_0=*_1=*$, the same conclusion holds if
$\varepsilon(\dot{\gamma}(0),\dot{\gamma}(1))=+1$.
\end{prop}
\begin{proof} Let $\varphi \in {\rm Diff}(S,\partial S)$ be a representative
of $\exp(\sigma(L(C)))$. By the remark above, we may assume that the support of $\varphi$
is included in $N$. We denote by the same letter $\varphi$ the restriction
of $\varphi$ to $N$, which we can regard as an element of the mapping
class group $\mathcal{M}(N,\partial N)$. Also we regard $C$ as an
unoriented free loop on $N$ and $L(C)$ as an element of $\widehat{\mathbb{Q}\hat{\pi}}(N)$.

Let $\partial N=\coprod_i \partial_iN$ be the decomposition into
connected components. Then, by \cite{KKp} Proposition 3.3.4,
there exist some $a_i \in \mathbb{Q}$ such that 
$$\varphi\exp(-\sigma(L(C))) = \exp \left(\sigma \left(\sum_i a_iL(\partial_iN)\right) \right)
\in {\rm Aut}(\widehat{\mathbb{Q}\mathcal{C}(N, \partial N)})$$
(see also the proof of \cite{KKp} Theorem 5.4.1).
Since $C$ and $\partial_i N$ are disjoint, the derivations
$L(C)$ and $L(\partial_i N)$ commute with each other.
This implies $\varphi^n = \exp(n\sigma(L(C) + \sum_i a_iL(\partial_i N)))$
for any $n \in \mathbb{Z}$. Since 
$\varphi^n(\gamma)$ is a simple path, we have $\mu(\varphi^n(\gamma)) =0$.
Hence we obtain
$$\mu \left(\sigma \left(L(C)+\sum_i a_iL(\partial_i N) \right)(\gamma)\right) = 0.$$
On the other hand, by \cite{KKp} Theorem 5.2.1,
$\exp(\sigma(L(\partial_iN)))$ is realized by the Dehn twist along the
simple closed curve $\partial_iN$. This implies
$\mu(\sigma(L(\partial_iN))(\gamma)) = 0$. Hence we obtain 
$\mu(\sigma(L(C))(\gamma)) = 0$. This completes the proof. 
\end{proof}

In the case $S$ is compact, i.e., has no punctures,
we have another criterion for the realizability of generalized Dehn twists.

\begin{prop}\label{deltaDT}
Assume $S$ is compact, and let $C \subset S\setminus \partial S$
be an unoriented loop whose generalized 
Dehn twist $\exp(\sigma(L(C)))$ is realized by a diffeomorphism. 
Then we have
$$
\delta L(C) = 0 \in
\widehat{\mathbb{Q}{\hat\pi}}(S) \widehat{\otimes}
\widehat{\mathbb{Q}{\hat\pi}}(S).
$$ 
\end{prop}
\begin{proof}
Take $*_0,*_1\in E$.
Any orientation-preserving diffeomorphism $\varphi$ of $S$
fixing the boundary pointwise preserves the comodule structure map
$\mu \colon
\widehat{\mathbb{Q}\Pi S}(*_0,*_1)\to
\widehat{\mathbb{Q}\Pi S}(*_0,*_1) \widehat{\otimes}
\widehat{\mathbb{Q}{\hat\pi}}(S)$.
Hence, for any $n \in \mathbb{Z}$, we have 
$$\mu \exp(n\sigma(L(C))) = \exp(n\sigma(L(C)))\mu,$$
and so 
$$(\sigma(L(C))\widehat{\otimes}1 + 1\widehat{\otimes}\sigma(L(C)))\mu
= \mu \sigma(L(C))\colon \widehat{\mathbb{Q}\Pi S}(*_0,*_1)\to
\widehat{\mathbb{Q}\Pi S}(*_0,*_1)\widehat{\otimes}
\widehat{\mathbb{Q}{\hat\pi}}(S).
$$
From (\ref{A6compatible+}) this is equivalent to 
$$
(\overline{\sigma}\widehat{\otimes}1_{\widehat{\mathbb{Q}{\hat\pi}}(S)})
(1_{\widehat{\mathbb{Q}\Pi S}(*_0,*_1)}\widehat{\otimes}\delta)(v\widehat{\otimes}
L(C)) = 0 \in \widehat{\mathbb{Q}\Pi S}(*_0,*_1) \widehat{\otimes}
\widehat{\mathbb{Q}{\hat\pi}}(S)
$$
for any $v \in \widehat{\mathbb{Q}\Pi S}(*_0,*_1)$.
By \cite{KKp} Theorem 6.2.1, the intersection of the kernels of the
structure map 
$\sigma\colon \widehat{\mathbb{Q}{\hat\pi}}(S) \to
{\rm End}({\widehat{\mathbb{Q}\Pi S}}(*_0,*_1))$ for all $*_0,*_1\in E$,
is zero. Hence we have $\delta L(C) = 0$. This proves the
proposition. 
\end{proof}

In \cite{Ku2} the second-named author posed the following question. 

\begin{quest}[\cite{Ku2} Question 5.3.4]
Let $C$ be an unoriented loop on $\Sigma_{g,1}$, a surface of genus $g$
with one boundary component, and suppose the
generalized Dehn twist along $C$ is realized by a diffeormorphism. 
Is $C$ homotopic to a power of a simple closed curve?
\end{quest}

In view of Proposition \ref{deltaDT} we come to the following conjecture. 

\begin{conj}\label{simple-d}
Suppose an unoriented loop $C$ satisfies $\delta L(C) = 0$. Then
$C$ would be homotopic to a power of a simple closed curve.
\end{conj}

If the conjecture is true, then the question is also
affirmative. But the conjecture looks like the question which was posed
by Turaev \cite{T2} and whose counter-examples Chas gave in \cite{Chas}. 

\subsection{New examples not realized by a diffeormorphism}

In this subsection we prove the following.

\begin{thm}
\label{main.thm}
Let $S$ and $E\subset \partial S$ be as in {\rm \S \ref{GDT}}
and $C\subset S\setminus \partial S$ an unoriented immersed loop
whose self intersections consist of transverse double points.
Assume $C$ is non-simple and the inclusion homomorphism 
$\pi_1(N(C))\to \pi_1(S)$ is injective, where $N(C)$ is a closed regular neighborhood of $C$.
Then the generalized Dehn twist $\exp(\sigma(L(C)))$ is not in the image of
$\widehat{\sf DN}\colon \mathcal{M}(S,\partial S) \to
{\rm Aut}(\widehat{\mathbb{Q}\mathcal{C}(S,E)})$.
\end{thm}

The rest of this subsection is devoted to the proof of Theorem \ref{main.thm}.

Let $S$ be an oriented surface and $*_0,*_1\in \partial S$
distinct points. Using $\mu$ and the augmentation
$\mathbb{Q}\Pi S(*_0,*_1)\to \mathbb{Q}, \Pi S(*_0,*_1) \ni x
\mapsto 1$, we define a $\mathbb{Q}$-linear map
$\hat{\mu}\colon \mathbb{Q}\Pi S(*_0,*_1) \to
\mathbb{Q}\hat{\pi}^{\prime}(S)$ as the composite
\begin{eqnarray*}
\hat{\mu}\colon \mathbb{Q}\Pi S(*_0,*_1) &\stackrel{\mu}{\to}&
\mathbb{Q}\Pi S(*_0,*_1) \otimes \mathbb{Q}\hat{\pi}^{\prime}(S) \\
&\to & \mathbb{Q}\otimes \mathbb{Q}\hat{\pi}^{\prime}(S)=
\mathbb{Q}\hat{\pi}^{\prime}(S).
\end{eqnarray*}
By Proposition \ref{mu-fil}, $\hat{\mu}$ extends to
a $\mathbb{Q}$-linear map $\hat{\mu}\colon \widehat{\mathbb{Q}\Pi S}(*_0,*_1)
\to \widehat{\mathbb{Q}\hat{\pi}}(S)$.
We denote by $\widehat{\mathbb{Q}H_1(S;\mathbb{Z})}$ the completed
group ring of the integral first homology group $H_1(S;\mathbb{Z})$.
There is a natural projection $\hat{\pi}(S)\to H_1(S;\mathbb{Z})$,
which induces a $\mathbb{Q}$-linear map
$\varpi\colon \widehat{\mathbb{Q}\hat{\pi}}(S)
\to \widehat{\mathbb{Q}H_1(S;\mathbb{Z})}/\mathbb{Q}1$.
Here $\mathbb{Q}1$ is the 1-dimensional subspace spanned
by the identity element of $H_1(S;\mathbb{Z})$.

Let $C\subset S\setminus \partial S$ be an unoriented immersed
loop such that its self intersections consist of transverse double points,
and let $\gamma \in \Pi S(*_0,*_1)$ be a simple path meeting $C$
transversally in a single point. In this situation, we shall
compute the quantity $\varpi \hat{\mu}(\sigma(L(C))\gamma)$.

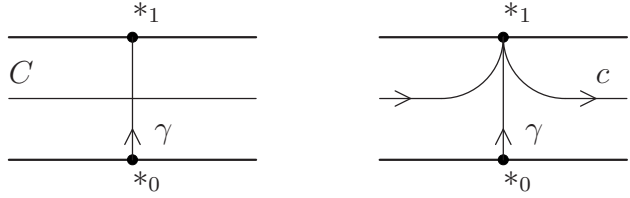
\begin{figure}
\begin{center}
\caption{$\pi_1(S,*_1)$-representative of $C$}

\vspace{0.2cm}

\unitlength 0.1in
\begin{picture}( 32.0000,  8.8000)(  4.0000,-10.2400)
%
{\color[named]{Black}{%
\special{pn 13}%
\special{pa 400 354}%
\special{pa 1680 354}%
\special{fp}%
}}%
%
{\color[named]{Black}{%
\special{pn 8}%
\special{pa 400 674}%
\special{pa 1680 674}%
\special{fp}%
}}%
%
{\color[named]{Black}{%
\special{pn 13}%
\special{pa 400 994}%
\special{pa 1680 994}%
\special{fp}%
}}%
%
{\color[named]{Black}{%
\special{pn 4}%
\special{sh 1}%
\special{ar 1040 994 26 26 0  6.28318530717959E+0000}%
}}%
%
{\color[named]{Black}{%
\special{pn 4}%
\special{sh 1}%
\special{ar 1040 354 26 26 0  6.28318530717959E+0000}%
}}%
%
{\color[named]{Black}{%
\special{pn 8}%
\special{pa 1040 994}%
\special{pa 1040 354}%
\special{fp}%
}}%
%
{\color[named]{Black}{%
\special{pn 8}%
\special{pa 1040 834}%
\special{pa 1080 914}%
\special{fp}%
}}%
%
{\color[named]{Black}{%
\special{pn 8}%
\special{pa 1040 834}%
\special{pa 1000 914}%
\special{fp}%
}}%
%
{\color[named]{Black}{%
\special{pn 13}%
\special{pa 2320 354}%
\special{pa 3600 354}%
\special{fp}%
}}%
%
{\color[named]{Black}{%
\special{pn 13}%
\special{pa 2320 994}%
\special{pa 3600 994}%
\special{fp}%
}}%
%
{\color[named]{Black}{%
\special{pn 4}%
\special{sh 1}%
\special{ar 2960 994 26 26 0  6.28318530717959E+0000}%
}}%
%
{\color[named]{Black}{%
\special{pn 4}%
\special{sh 1}%
\special{ar 2960 354 26 26 0  6.28318530717959E+0000}%
}}%
%
{\color[named]{Black}{%
\special{pn 8}%
\special{pa 2960 994}%
\special{pa 2960 354}%
\special{fp}%
}}%
%
{\color[named]{Black}{%
\special{pn 8}%
\special{pa 2960 834}%
\special{pa 3000 914}%
\special{fp}%
}}%
%
{\color[named]{Black}{%
\special{pn 8}%
\special{pa 2960 834}%
\special{pa 2920 914}%
\special{fp}%
}}%
%
{\color[named]{Black}{%
\special{pn 8}%
\special{ar 3280 354 320 320  1.5707963 3.1415927}%
}}%
%
{\color[named]{Black}{%
\special{pn 8}%
\special{ar 2640 354 320 320  6.2831853 6.2831853}%
\special{ar 2640 354 320 320  0.0000000 1.5707963}%
}}%
%
{\color[named]{Black}{%
\special{pn 8}%
\special{pa 2320 674}%
\special{pa 2640 674}%
\special{fp}%
}}%
%
{\color[named]{Black}{%
\special{pn 8}%
\special{pa 3280 674}%
\special{pa 3600 674}%
\special{fp}%
}}%
%
{\color[named]{Black}{%
\special{pn 8}%
\special{pa 3440 674}%
\special{pa 3360 634}%
\special{fp}%
}}%
%
{\color[named]{Black}{%
\special{pn 8}%
\special{pa 3440 674}%
\special{pa 3360 714}%
\special{fp}%
}}%
%
{\color[named]{Black}{%
\special{pn 8}%
\special{pa 2480 674}%
\special{pa 2400 714}%
\special{fp}%
}}%
%
{\color[named]{Black}{%
\special{pn 8}%
\special{pa 2480 674}%
\special{pa 2400 634}%
\special{fp}%
}}%
\put(10.4000,-11.5400){\makebox(0,0)[lb]{$*_0$}}%
\put(29.6000,-11.5400){\makebox(0,0)[lb]{$*_0$}}%
\put(29.6000,-2.7400){\makebox(0,0)[lb]{$*_1$}}%
\put(10.4000,-2.7400){\makebox(0,0)[lb]{$*_1$}}%
\put(4.0000,-5.9400){\makebox(0,0)[lb]{$C$}}%
\put(11.5200,-9.1400){\makebox(0,0)[lb]{$\gamma$}}%
\put(30.7200,-9.1400){\makebox(0,0)[lb]{$\gamma$}}%
\put(34.4000,-5.9400){\makebox(0,0)[lb]{$c$}}%
\end{picture}%
\end{center}
\end{figure}

Let $c$ be a $\pi_1(S,*_1)$-representative of $C$, as in Figure 9.
Then we have
$$\sigma(L(C))\gamma =\gamma \log c \in \widehat{\mathbb{Q}\Pi S}(*_0,*_1),$$
since $\sigma(|c^n|)\gamma=n\gamma c^n$ for $n\ge 0$. Now fix a
parametrization $c\colon ([0,1],\{ 0,1\}) \to (S,*_1)$.
When $p\in S$ is a double point of $C$, we denote $c^{-1}(p)=\{ t_1^p,t_2^p \}$
so that $t_1^p <t_2^p$. Set $x_p:=c_{t_2^p1}c_{0t_1^p}$, and $y_p:=c_{t_1^pt_2^p}$.
By abuse of notation, we use the same letter $x_p$ and $y_p$ for
the homology classes represented by these loops. Finally,
let $h(x)$ be the formal power series defined by
$$h(x):=\sum_{n=0}^{\infty}\frac{(-1)^n}{n+2}(x-1)^n.$$

\begin{prop}
\label{pimu}
Keep the notations as above. then
\begin{eqnarray*}
& & \varpi \hat{\mu}(\sigma(L(C))\gamma) \\
&=& -\sum_p \varepsilon(\dot{c}(t_1^p),\dot{c}(t_2^p))
(y_p+x_p(y_p^2-1)h(c)) \in \widehat{\mathbb{Q}H_1(S;\mathbb{Z})}/\mathbb{Q}1.
\end{eqnarray*}
Here we write the product of the group ring $\mathbb{Q}H_1(S;\mathbb{Z})$
multiplicatively.
\end{prop}

To prove this proposition, we need a lemma.

\begin{lem}
\label{f_n(x)}
In the polynomial ring $\mathbb{Q}[x]$,
the following equalities hold.
\begin{enumerate}
\item For $n\ge 1$,
$$\sum_{k=0}^n \binom{n}{k}(-1)^{n-k}\sum_{j=0}^{k-1}x^j=(x-1)^{n-1}.$$
\item For $n\ge 1$, set
$$f_n(x):=\sum_{k=0}^n \binom{n}{k}(-1)^{n-k}
\sum_{j=0}^{k-1} (k-j)x^j.$$
Then $f_1(x)=1$ and $f_n(x)=x(x-1)^{n-2}$ for $n\ge 2$.
\end{enumerate}
\end{lem}

\begin{proof}
1. Since $\sum_{j=0}^{k-1}x^j=(x^k-1)/(x-1)$, the left hand side is
equal to
\begin{eqnarray*}
\frac{1}{x-1}\sum_{k=0}^n \binom{n}{k}(-1)^{n-k}(x^k-1)
&=& \frac{1}{x-1}\sum_{k=0}^n \binom{n}{k}(-1)^{n-k}x^k \\
&=& \frac{1}{x-1}(x-1)^n=(x-1)^{n-1}.
\end{eqnarray*}

2. The case $n=1$ is clear. Let $n\ge 2$. By the first part, we compute
\begin{eqnarray*}
f_n(x)-(x-1)^{n-1}
&=& \sum_{k=0}^n \binom{n}{k}(-1)^{n-k}\sum_{j=0}^{k-1}
(k-1-j)x^j \\
&=& \sum_{k=0}^n \left( \binom{n-1}{k-1}+\binom{n-1}{k} \right)
(-1)^{n-k}\sum_{j=0}^{k-1}(k-1-j)x^j \\
&=& \sum_{k=0}^{n-1} \binom{n-1}{k}(-1)^{n-k-1}
\sum_{j=0}^k(k-j)x^j \\
& &+\sum_{k=0}^n \binom{n-1}{k}(-1)^{n-k}
\sum_{j=0}^{k-1} (k-1-j)x^j \\
&=& \sum_{k=0}^{n-1} \binom{n-1}{k} (-1)^{n-1-k}
\sum_{j=0}^{k-1} x^j=(x-1)^{n-2}.
\end{eqnarray*}
Therefore $f_n(x)=(x-1)^{n-1}+(x-1)^{n-2}=x(x-1)^{n-2}$.
\end{proof}

\begin{proof}[Proof of Proposition \ref{pimu}]
We first compute $\hat{\mu}(\gamma c^k)$ for $k\ge 0$.
We choose a representative of $\gamma c^k$ by sliding
$c$ into the left. See Figure 10. Each self intersection
$p$ of $C$ creates $k^2$ self intersections of $\gamma c^k$.
These $k^2$ points are classified into $k+(k-1)$ classes,
according to their contributions to $\hat{\mu}(\gamma c^k)$.
Namely, if $\varepsilon(\dot{c}(t_1^p),\dot{c}(t_2^p))=1$,
for $0\le j \le k-1$, there are $k-j$ self intersections
whose contributions are $-|y_p(x_py_p)^j|^{\prime}$,
and for $1\le j\le k-1$, there are $k-j$ self intersections
whose contributions are $+|x_p(y_px_p)^{j-1}|^{\prime}$.
This is illustrated in Figure 11. The points in the box $j^+$
($0\le j \le k-1$) contribute as $-|y_p(x_py_p)^j|^{\prime}$,
and the points in the box $j^-$ ($1\le j\le k-1$) contribute
as $+|x_p(y_px_p)^{j-1}|^{\prime}$.
If $\varepsilon(\dot{c}(t_1^p),\dot{c}(t_2^p))=-1$,
the contributions are the minus of the case of
$\varepsilon(\dot{c}(t_1^p),\dot{c}(t_2^p))=1$.
Therefore, we obtain
\begin{equation}
\label{gc^k}
\hat{\mu}(\gamma c^k)=-\sum_p \varepsilon(\dot{c}(t_1^p),\dot{c}(t_2^p))
\left( \sum_{j=0}^{k-1}(k-j)|y_p(x_py_p)^j|^{\prime}
-\sum_{j=1}^{k-1}(k-j)|x_p(y_px_p)^{j-1}|^{\prime} \right).
\end{equation}
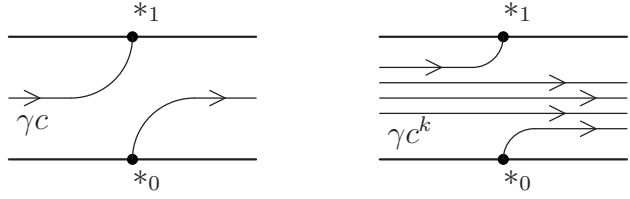
\begin{figure}
\begin{center}
\caption{a representative of $\gamma c^k$ ($k=4$)}

\vspace{0.2cm}

\unitlength 0.1in
\begin{picture}( 32.0000,  8.8000)(  4.0000,-10.2400)
%
{\color[named]{Black}{%
\special{pn 13}%
\special{pa 400 354}%
\special{pa 1680 354}%
\special{fp}%
}}%
%
{\color[named]{Black}{%
\special{pn 13}%
\special{pa 400 994}%
\special{pa 1680 994}%
\special{fp}%
}}%
%
{\color[named]{Black}{%
\special{pn 4}%
\special{sh 1}%
\special{ar 1040 994 26 26 0  6.28318530717959E+0000}%
}}%
%
{\color[named]{Black}{%
\special{pn 4}%
\special{sh 1}%
\special{ar 1040 354 26 26 0  6.28318530717959E+0000}%
}}%
%
{\color[named]{Black}{%
\special{pn 13}%
\special{pa 2320 354}%
\special{pa 3600 354}%
\special{fp}%
}}%
%
{\color[named]{Black}{%
\special{pn 13}%
\special{pa 2320 994}%
\special{pa 3600 994}%
\special{fp}%
}}%
%
{\color[named]{Black}{%
\special{pn 4}%
\special{sh 1}%
\special{ar 2960 994 26 26 0  6.28318530717959E+0000}%
}}%
%
{\color[named]{Black}{%
\special{pn 4}%
\special{sh 1}%
\special{ar 2960 354 26 26 0  6.28318530717959E+0000}%
}}%
\put(10.4000,-11.5400){\makebox(0,0)[lb]{$*_0$}}%
\put(29.6000,-11.5400){\makebox(0,0)[lb]{$*_0$}}%
\put(29.6000,-2.7400){\makebox(0,0)[lb]{$*_1$}}%
\put(10.4000,-2.7400){\makebox(0,0)[lb]{$*_1$}}%
%
{\color[named]{Black}{%
\special{pn 8}%
\special{ar 1360 994 320 320  3.1415927 4.7123890}%
}}%
%
{\color[named]{Black}{%
\special{pn 8}%
\special{ar 720 354 320 320  6.2831853 6.2831853}%
\special{ar 720 354 320 320  0.0000000 1.5707963}%
}}%
%
{\color[named]{Black}{%
\special{pn 8}%
\special{pa 400 674}%
\special{pa 720 674}%
\special{fp}%
}}%
%
{\color[named]{Black}{%
\special{pn 8}%
\special{pa 1360 674}%
\special{pa 1680 674}%
\special{fp}%
}}%
%
{\color[named]{Black}{%
\special{pn 8}%
\special{pa 1520 674}%
\special{pa 1440 714}%
\special{fp}%
}}%
%
{\color[named]{Black}{%
\special{pn 8}%
\special{pa 1520 674}%
\special{pa 1440 634}%
\special{fp}%
}}%
%
{\color[named]{Black}{%
\special{pn 8}%
\special{pa 560 674}%
\special{pa 480 714}%
\special{fp}%
}}%
%
{\color[named]{Black}{%
\special{pn 8}%
\special{pa 560 674}%
\special{pa 480 634}%
\special{fp}%
}}%
\put(4.4000,-8.6600){\makebox(0,0)[lb]{$\gamma c$}}%
%
{\color[named]{Black}{%
\special{pn 8}%
\special{ar 3120 994 160 160  3.1415927 4.7123890}%
}}%
%
{\color[named]{Black}{%
\special{pn 8}%
\special{ar 2800 354 160 160  6.2831853 6.2831853}%
\special{ar 2800 354 160 160  0.0000000 1.5707963}%
}}%
%
{\color[named]{Black}{%
\special{pn 8}%
\special{pa 3120 834}%
\special{pa 3600 834}%
\special{fp}%
}}%
%
{\color[named]{Black}{%
\special{pn 8}%
\special{pa 2320 514}%
\special{pa 2800 514}%
\special{fp}%
}}%
%
{\color[named]{Black}{%
\special{pn 8}%
\special{pa 2320 754}%
\special{pa 3600 754}%
\special{fp}%
}}%
%
{\color[named]{Black}{%
\special{pn 8}%
\special{pa 2320 594}%
\special{pa 3600 594}%
\special{fp}%
}}%
%
{\color[named]{Black}{%
\special{pn 8}%
\special{pa 2320 674}%
\special{pa 3600 674}%
\special{fp}%
}}%
%
{\color[named]{Black}{%
\special{pn 8}%
\special{pa 3440 834}%
\special{pa 3360 874}%
\special{fp}%
\special{pa 3440 834}%
\special{pa 3360 794}%
\special{fp}%
}}%
%
{\color[named]{Black}{%
\special{pn 8}%
\special{pa 3440 674}%
\special{pa 3360 714}%
\special{fp}%
\special{pa 3440 674}%
\special{pa 3360 634}%
\special{fp}%
}}%
%
{\color[named]{Black}{%
\special{pn 8}%
\special{pa 3280 754}%
\special{pa 3200 794}%
\special{fp}%
\special{pa 3280 754}%
\special{pa 3200 714}%
\special{fp}%
}}%
%
{\color[named]{Black}{%
\special{pn 8}%
\special{pa 3280 594}%
\special{pa 3200 634}%
\special{fp}%
\special{pa 3280 594}%
\special{pa 3200 554}%
\special{fp}%
}}%
%
{\color[named]{Black}{%
\special{pn 8}%
\special{pa 2640 514}%
\special{pa 2560 554}%
\special{fp}%
\special{pa 2640 514}%
\special{pa 2560 474}%
\special{fp}%
}}%
\put(23.6000,-9.4000){\makebox(0,0)[lb]{$\gamma c^k$}}%
\end{picture}%
\end{center}
\end{figure}

\begin{figure}
\begin{center}
\caption{a picture near $p$
($\varepsilon(\dot{c}(t_1^p),\dot{c}(t_2^p))=1$, $k=4$)}

\vspace{0.2cm}

\input{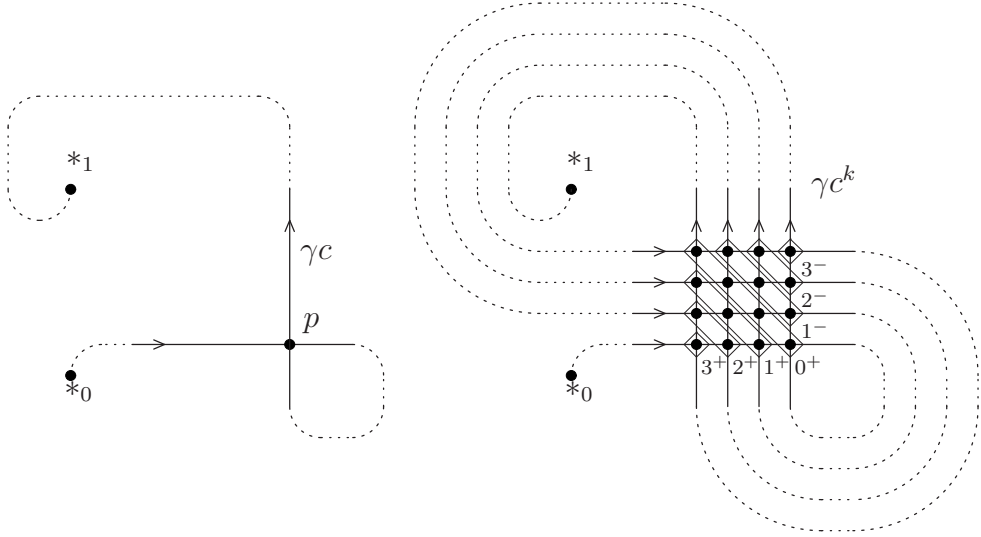}
\end{center}
\end{figure}

We next compute $\hat{\mu}(\gamma (c-1)^n)$ for $n\ge 1$.
We claim that the contribution from $p$ to $\hat{\mu}(\gamma (c-1)^n)$
is $-\varepsilon(\dot{c}(t_1^p),\dot{c}(t_2^p))|y_p|^{\prime}$ if
$n=1$, and
$$-\varepsilon(\dot{c}(t_1^p),\dot{c}(t_2^p))
\left( |y_px_py_p(x_py_p-1)^{n-2}|^{\prime}
-|x_p(y_px_p-1)^{n-2}|^{\prime} \right)$$
if $n\ge 2$. The case $n=1$ is clear. If $n\ge 2$,
by (\ref{gc^k}) the contribution
is $-\varepsilon(\dot{c}(t_1^p),\dot{c}(t_2^p))$ times
\begin{equation}
\label{g(c-1)^n}
\sum_{k=0}^n \binom{n}{k} (-1)^{n-k}
\left( \sum_{j=0}^{k-1}(k-j)|y_p(x_py_p)^j|^{\prime}
-\sum_{j=1}^{k-1}(k-j)|x_p(y_px_p)^{j-1}|^{\prime} \right).
\end{equation}

By
$$\sum_{j=1}^{k-j}(k-j)|x_p(y_px_p)^{j-1}|^{\prime}
=\sum_{j=0}^{k-1}(k-j)|x_p(y_px_p)^j|^{\prime}
-\sum_{j=0}^{k-1}|x_p(y_px_p)^j|^{\prime}$$
and Lemma \ref{f_n(x)}, (\ref{g(c-1)^n}) is equal to
\begin{eqnarray*}
& &
|y_px_py_p(x_py_p-1)^{n-2}|^{\prime}-
|x_py_px_p(y_px_p-1)^{n-2}|^{\prime}+|x_p(y_px_p-1)^{n-1}|^{\prime} \\
&=&
|y_px_py_p(x_py_p-1)^{n-2}|^{\prime}-|x_p(y_px_p-1)^{n-2}|^{\prime}.
\end{eqnarray*}
The claim is proved. Now we conclude
$$\hat{\mu}(\gamma \log c)=
-\sum_p \varepsilon(\dot{c}(t_1^p),\dot{c}(t_2^p))
\left( |y_p|^{\prime}+|y_px_py_ph(x_py_p)|^{\prime}
-|x_ph(y_px_p)|^{\prime} \right).$$
Applying $\varpi$ and using $x_py_p=c=y_px_p \in H_1(S;\mathbb{Z})$,
we obtain the desired formula. This completes the proof.
\end{proof}

\begin{proof}[Proof of Theorem \ref{main.thm}]
Assume the generalized Dehn twist $\exp(\sigma(L(C)))$ is
realized by a diffeomorphism. Let $N=N(C)$ be a closed regular neighborhood of $C$.
Take a simple point $a\in S$ of $C$
and let $\gamma\colon ([0,1],\{ 0,1\}) \to (N,\partial N)$
be a simple path in $N$ meeting $C$ transversally only at $a$.
We denote $\gamma(0)=*_0$ and $\gamma(1)=*_1$.
By Proposition \ref{criterion}, we have $\mu(\sigma(L(C))\gamma)=0
\in \widehat{\mathbb{Q}\Pi N}(*_0,*_1)\widehat{\otimes}
\widehat{\mathbb{Q}\hat{\pi}}(N)$.
In particular, we have $\varpi \hat{\mu}(\sigma(L(C))\gamma)=0
\in \widehat{\mathbb{Q}H_1(N;\mathbb{Z})}/\mathbb{Q}1$.

We claim: 1) $\{ x_p\}_p \cup \{ c\}$ constitute a $\mathbb{Z}$-basis
of $H_1(N;\mathbb{Z})=H_1(C;\mathbb{Z})$, and 2) by an appropriate choice of $a$,
we can arrange that $\sum_p \varepsilon(\dot{c}(t_1^p),\dot{c}(t_2^p))
\neq 0$.

To prove the first claim, note that only
the underlyng 4-valent graph structure of $C$, together
with its (unoriented) parametrization matters.
We proceed by induction on the number
of double points of $C$. If $C$ is simple, the claim is trivial.
Suppose $C$ has at least one self intersection and let $q$ be a double
point of $C$. Let $f\colon \widetilde{C}\to C$ be a resolution of $q$.
Namely, $\widetilde{C}$ is a 4-valent graph with a
surjective map $S^1 \to \widetilde{C}$,
such that the composition $S^1\to \widetilde{C} \to C$ gives
a parametrization of $C$, $f^{-1}(x)$ consist of a single
point if $x\neq q$, and $f^{-1}(q)$ consist of two points,
say $q_+$ and $q_-$.

By the excision isomorphism, we have $H_1(C;\mathbb{Z})
= H_1(C,\{ q\};\mathbb{Z})\cong H_1(\widetilde{C},\{ q_+,q_-\};\mathbb{Z})$.
Consider the homology exact sequence of the pair
$$0\to H_1(\widetilde{C};\mathbb{Z}) \to
H_1(\widetilde{C},\{ q_+,q_- \};\mathbb{Z}) \stackrel{\partial}{\to}
\widetilde{H}_0(\{ q_+,q_- \};\mathbb{Z}) \to 0.$$
Then the $\partial$-image of
$x_q\in H_1(C;\mathbb{Z})=H_1(\widetilde{C},\{ q_+,q_-\};\mathbb{Z})$
is $\pm (q_+-q_-)$, which is a generator of
$\widetilde{H}_0(\{ q_+,q_- \};\mathbb{Z})\cong \mathbb{Z}$.
By the inductive assumption, the lifts of $\{x_p\}_{p\neq q}\cup \{ c\}$ to
$\widetilde{C}$ constitute a $\mathbb{Z}$-basis of
$H_1(\widetilde{C};\mathbb{Z})$. Therefore the lifts of $\{ x_p \}_p\cup \{ c\}$
to $\widetilde{C}$ constitute a $\mathbb{Z}$-basis of
$H_1(\widetilde{C},\{ q_+,q_- \};\mathbb{Z})$, which completes
the proof of the first claim.

To prove the second claim, let $\ell$ be the component of the
set of simple points of $C$ containing $a$, and $\ell^{\prime}$
a component next to $\ell$. Take a simple point
$a^{\prime}\in \ell^{\prime}$ and let $\gamma^{\prime}$ be
a simple path meeting $C$ transversally only at $a^{\prime}$.
We arrange that $\varepsilon(\dot{c}(a),\dot{\gamma}(a))
=\varepsilon(\dot{c}(a^{\prime}),\dot{\gamma}^{\prime}(a^{\prime}))$.
Let $q$ be the double point of $C$ between $\ell$ and $\ell^{\prime}$.
We compare $\varepsilon(\dot{c}(t_1^p),\dot{c}(t_2^p))$'s
with respect to $a$ and $a^{\prime}$. If $p\neq q$, then
they are the same. If $p=q$, they are minus of each other.
Hence the difference of the sums $\sum_p \varepsilon(\dot{c}(t_1^p),\dot{c}(t_2^p))$
for $a$ and $a^{\prime}$ is two, in particular
at least one of them is not zero. This proves the second claim.

Now choose $a$ such that $\sum_p \varepsilon(\dot{c}(t_1^p),\dot{c}(t_2^p))
\neq 0$. By the first claim, we can define a group homomorphism
$\Phi\colon H_1(N;\mathbb{Z})\to \langle t \rangle$
to an infinite cyclic group generated by $t$, by
$\Phi(x_p)=1$ and $\Phi(c)=t$. This group homomorphism induces a $\mathbb{Q}$-linear
map $\Phi\colon \widehat{\mathbb{Q}H_1(N;\mathbb{Z})}/\mathbb{Q}1
\to \widehat{\mathbb{Q}\langle t \rangle}/\mathbb{Q}1$.
Since $x_py_p=c$, we have $\Phi(y_p)=t$. By Proposition \ref{pimu},
$$\Phi(\varpi \mu(\sigma(L(C))\gamma))=
-\left( \sum_p \varepsilon(\dot{c}(t_1^p),\dot{c}(t_2^p)) \right)
(t+(t^2-1)h(t)) \in \widehat{\mathbb{Q}\langle t \rangle}/\mathbb{Q}1.$$
Finally, we claim $t+(t^2-1)h(t) \neq 0$. To prove this,
consider an algebra homomorphism from $\widehat{\mathbb{Q}\langle t \rangle}$
to $\mathbb{Q}[[s]]$, the ring of formal power series in $s$, given by
$t\mapsto 1+s$. This is a filter-preserving isomorphism, and the image of $t+(t^2-1)h(t)$ is
$1+2s+({\rm higher\ term\ })$, which is not zero in $\mathbb{Q}[[s]]/\mathbb{Q}1$.
This shows $t+(t^2-1)h(t)\neq 0$, which contradicts to
$\varpi \mu(\sigma(L(C))\gamma)=0$. This completes the proof.
\end{proof}

\section{A geometric approach to the Johnson homomorphisms}
\label{sec:App2}

In this section we study the Johnson homomorphisms following the treatments in \cite{KKp},
which is briefly recalled in \S5.1. In \S5.2, we show that the Turaev cobracket gives a geometric constraint on the Johnson image.
As was shown in \cite{KK1}, in the case $S = \Sigma_{g,1}$, a once bordered surface
of genus $g\ge 1$,
the completed Goldman Lie algebra is isomorphic to the Lie algebra of symplectic derivations
of the completed tensor algebra generated by the first rational homology group of the surface, $\mathfrak{a}_g^-$,
through a symplectic expansion $\theta$ introduced by Massuyeau \cite{Mas}.
This isomorphism induces a complete Lie cobracket $\delta^\theta$ on the Lie algebra $\mathfrak{a}_g^-$.
We can consider the ``Laurent expansion" of the cobracket $\delta^\theta$ with respect to the natural degree on $\mathfrak{a}_g^-$.
In \S5.3, based on a theorem of Massuyeau and Turaev \cite{MT}, we prove that the principal term, which is of degree $-2$,
equals Schedler's cobracket \cite{Sch}. Moreover we prove that the $(-1)$-st and the $0$-th terms vanish.
The latter term is computed in \S 5.5. 
These results are obtained independently by Massuyeau and Turaev \cite{MT3}.
In \S5.4, we prove that all the Morita traces \cite{MoAQ} factor through Schedler's cobracket.
As a corollary, all the Morita traces are outside of our geometric constraint.

\subsection{The Johnson homomorphisms}

The higher Johnson homomorphisms on the higher Torelli groups for a once
bordered surface are important tools to study the algebraic structure of the mapping class group.
In \cite{KKp}, we gave a generalization of the classical construction of the Johnson
homomorphisms to arbitrary compact oriented surfaces with non-empty boundaries.
In this subsection we briefly recall this construction.

Let $S$ be a compact connected oriented surface of with non-empty boundary,
and $E \subset \partial S$ a subset such that each
connected component of $\partial S$ has a unique point of $E$.
We define the {\it Torelli group} $\ISE$ to be the kernel of the action of the
mapping class group $\mathcal{M}(S, E)$ on the first homology group
$H_1(S,E; \mathbb{Z})$, which is the smallest Torelli group in the sense
of Putman \cite{P}. On the other hand, we define 
$$
L^+(S,E) := \{u \in \widehat{\mathbb{Q}\hat\pi}(S)(3);
(\sigma(u)\widehat{\otimes}\sigma(u))\circ \Delta^{(*_0,*_1)} = \Delta^{(*_0,*_1)}\circ\sigma(u) \,\,\mbox{for any $*_0,*_1 \in E$}\}, 
$$
where $\Delta$ is the coproduct $\Delta = \Delta^{(*_0,*_1)}\colon
\widehat{\mathbb{Q}\Pi S}(*_0,*_1) \to \widehat{\mathbb{Q}\Pi
S}(*_0,*_1)\widehat{\otimes}\widehat{\mathbb{Q}\Pi S}(*_0,*_1)$ given by
$\Delta x:= x\widehat{\otimes}x$ for any $x \in \Pi S(*_0,*_1)$, $*_0,*_1
\in E$. 

Using the Hausdorff series, we can regard
$L^+(S,E)$ as a pro-nilpotent group. In other words, using the injectivity of
$\sigma$ (\cite{KKp} Theorem 6.2.1) and the exponential map, we have a bijection $L^+(S,E)\stackrel{\cong}{\to}
\exp (\sigma(L^+(S,E)))\subset {\rm Aut}(\widehat{\mathbb{Q}\mathcal{C}(S,E)})$,
which endows $L^+(S,E)$ with a group structure.
In \cite{KKp} \S6.3 we showed the inclusion
$$\widehat{\sf DN}(\ISE )\subset \exp (\sigma(L^+(S,E))),$$
using a result of Putman \cite{P} about generators of $\ISE$ and our formula
for Dehn twists \cite{KK1}. Hence we obtain a unique injective group homomorphism
\begin{equation}
\tau\colon \ISE \to L^+(S,E)
\end{equation}
satisfying $\widehat{\sf DN}|_{\ISE}=\exp \circ \sigma \circ \tau$.
We call it {\it the geometric Johnson homomorphism} of the Torelli group $\ISE$. 
\par
In the case $S= \Sigma_{g,1}$ and $E$ consists of a single point $*\in \partial S$,
the Torelli group $\ISE$ is just the classical Torelli group $\mathcal{I}_{g,1}$.
As was shown in \S 6.3 \cite{KKp}, 
the map $\tau$ is essentially the same as Massuyeau's improvement \cite{Mas} of
the total Jonson map \cite{Ka}. 
In particular, the graded quotients of the geometric Johnson homomorphism
$\tau\colon \mathcal{I}_{g,1} \to L^+(\Sigma_{g,1}, \{*\})$ with respect
to the Johnson filtration on $\mathcal{I}_{g,1}$ and the filtration
$\{L^+(\Sigma_{g,1}, \{*\})\cap \widehat{\mathbb{Q}\hat\pi}(\Sigma_{g,1})(p)\}_{p \geq 3}$
are exactly the Johnson homomorphims of all degrees. 

\subsection{A constraint on the Johnson image}

Now we show that the Turaev cobracket gives an obstruction of the surjectivity of $\tau$.

\begin{thm}\label{zero} 
$$
\delta\circ\tau = 0\colon \ISE
\overset\tau\to L^+(S,E) \subset \widehat{\mathbb{Q}{\hat\pi}}(S)
\overset\delta\to \widehat{\mathbb{Q}{\hat\pi}}(S)\widehat{\otimes}
\widehat{\mathbb{Q}{\hat\pi}}(S).
$$
\end{thm}

\begin{proof} 
The proof is similar to that of Proposition \ref{deltaDT}. 
From the definition of $\tau$, for any $\varphi \in
\ISE$, there exists a unique $u\in L^+(S,E)$ such that 
$\varphi = \exp\sigma(u)$ on $\widehat{\mathbb{Q}\Pi S}(*_0,*_1)$ for any
$*_0$ and $*_1 \in E$. Then we have 
$\tau(\varphi) = u$ by definition. 
Let $\mu\colon \widehat{\mathbb{Q}\Pi S}(*_0,*_1) \to
\widehat{\mathbb{Q}\Pi
S}(*_0,*_1)\widehat{\otimes}\widehat{\mathbb{Q}{\hat\pi}}(S)$ be the
structure map of the comodule $\widehat{\mathbb{Q}\Pi S}(*_0,*_1)$. 
It is clear that $\mu$ is preserved by $\varphi^n$ for any $n \in
\mathbb{Z}$, namely, we have
$$
(\exp\sigma(nu)\widehat{\otimes}\exp\sigma(nu))\mu(v)  =
\mu(\exp\sigma(nu)(v))
$$ 
for any $n \in \mathbb{Z}$ and $v \in \widehat{\mathbb{Q}\Pi
S}(*_0,*_1)$.  Hence we have 
$$
(\sigma(u)\widehat{\otimes}1 + 1\widehat{\otimes}\sigma(u))\mu(v) 
= \mu(\sigma(u)(v))
$$
which is equivalent to 
$$
(\overline{\sigma}\otimes 1)(v\widehat{\otimes}\delta u) = 0
\in \widehat{\mathbb{Q}\Pi
S}(*_0,*_1)\widehat{\otimes}\widehat{\mathbb{Q}{\hat\pi}}(S)
$$
for any $*_0$ and $*_1 \in E$, 
from (\ref{A6compatible+}). Again by \cite{KKp} Theorem 6.2.1,
we conclude $\delta u = 0$. This proves the theorem. 
\end{proof}

\begin{figure}
\begin{center}
\caption{the case $g=3$, $r=2$}

\vspace{0.2cm}

\input{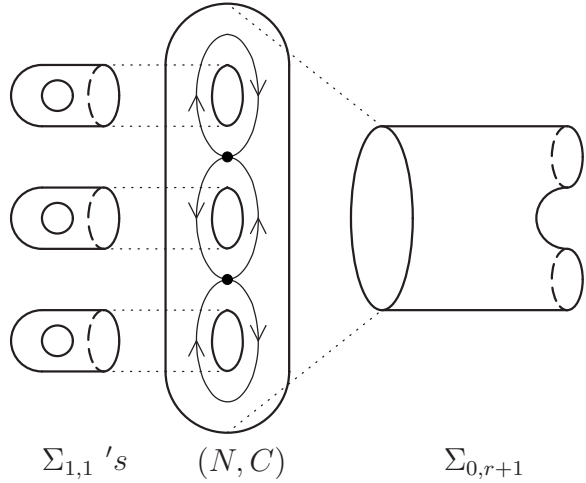}

\end{center}
\end{figure}

This constraint is non-trivial if the genus of the surface $S$ is greater than $1$.
\begin{prop} 
If $g \geq 2$, we have $\delta\vert_{L^+(S,E)} \neq 0$.
\end{prop}
\begin{proof} We denote by $\Sigma_{g,r}$ a compact connected oriented
surface of genus $g$ with $r$ boundary components.
Consider a spine $C$ of the surface $N := \Sigma_{0, g+1}$ as in Figure 12.
If $g \geq 2$, $C$ has a self-intersection. We cap each of the $g$
boundaries the curve $C$ surrounds by a surface diffeomorphic to
$\Sigma_{1,1}$ to obtain a compact surface $S_0$ diffeomorphic to
$\Sigma_{g,1}$, and glue $\Sigma_{0, r+1}$ to the boundary of $S_0$
to get a compact surface $S$ diffeomorphic to $\Sigma_{g,r}$.
See Figure 12. Choose one point in each boundary component of $S$.
We define $E$ by the set of all these points. 
We may regard $N$ as a regular neighborhood of $C$. \par
Consider the invariant $L(C) \in \widehat{\mathbb{Q}\hat\pi}(S)$.
As was proved in \cite{KKp} Lemma 5.1.2, the action of $L(C)$ stabilizes
the coproduct $\Delta$. Since $[C] = 0 \in H_1(S; \mathbb{Q})$,
we have $L(C) \in L^+(S, E)$. 
From the proof of Theorem \ref{main.thm} and the compatibility of
the comodule structure map and the cobracket (\ref{A6compatible+}),
we have $\delta L(C) \neq 0 \in 
\widehat{\mathbb{Q}\hat\pi}(N)\widehat{\otimes}\widehat{\mathbb{Q}\hat\pi}(N)$. 
As was proved in \cite{KKp} Proposition 6.2.3, the inclusion homomorphism 
$\widehat{\mathbb{Q}\hat{\pi}}(N) \to \widehat{\mathbb{Q}\hat{\pi}}(S_0)$ is injective. 
Since the inclusion homomorphism $\pi_1(S_0) \to \pi_1(S)$ has a right
inverse coming from capping all the boundaries except one by $r-1$ discs,
the inclusion homomorphism 
$\widehat{\mathbb{Q}\hat{\pi}}(S_0) \to \widehat{\mathbb{Q}\hat{\pi}}(S)$ is injective. 
Hence we have $\delta L(C) \neq 0 \in 
\widehat{\mathbb{Q}\hat\pi}(S)\widehat{\otimes}\widehat{\mathbb{Q}\hat\pi}(S)$.
This proves the proposition.
\end{proof}

From Theorem \ref{zero} the Zariski closure of the subgroup $\tau(\ISE)$ is
included in the closed Lie subalgebra $\Ker(\delta\vert_{L^+(S,E)})$. 
In view of this theorem we raise the following conjecture.
\begin{conj}\label{Zariski}
The Zariski closure of the subgroup $\tau(\ISE)$ equals the
closed Lie subalgebra $\Ker\,(\delta\vert_{L^+(S,E)})$:
$$
\overline{\tau(\ISE)} = \Ker\,(\delta\vert_{L^+(S,E)}).
$$
\end{conj}
This conjecture questions the extensionality of the Johnson image.
Hain \cite{Hain} already described its comprehension.
By Turaev's theorem \cite{T1}, p.234, Corollary 2,  $\mu$ captures the
simplicity of a based loop on a surface. This conjecture is its analogue
in the mapping class group. It is closely related to Conjecture
\ref{simple-d}. But it seems quite optimistic even in the simplest case $S = \Sigma_{g,1}$.
The cokernel of the Johnson homomorphisms in the case $S=\Sigma_{g,1}$ is known to have plenty of
$Sp$-irreducible components including the Morita traces \cite{MoAQ}. For
detalis, see \cite{ES} and references therein. In the succeeding subsections 
we will prove that all the Morita traces are outside of our constraint.
Very recently Enomoto \cite{eno} proved that the Enomoto-Satoh traces \cite{ES} are {\it inside}
of the leading term $\delta^{\rm alg}$ in the ``Laurent expansion" of the Turaev cobracket $\delta$
in \S \ref{sec:GQT}. But we do not know whether they are inside of $\delta$ itself or not.

\subsection{The graded quotient of the Turaev cobracket}
\label{sec:GQT}
For the rest of this section we suppose
$S=\Sigma_{g,1}$, a compact connected oriented surface of genus $g$ with one boundary component.
Then the Torelli group $\mathcal{I}(S,E) = \mathcal{I}(\Sigma_{g,1}, \{*\})$ 
is classically denoted by $\mathcal{I}_{g,1}$. 
Moreover as we will briefly recall below (for details, see \cite{KKp} \S 6.3),
the Lie algebra $L^+(\Sigma_{g,1},\{ *\})$ is identified with
(the completion of) the positive part of Kontsevich's Lie, $\ell_g^+$ \cite{Kon}.
Preceding Kontsevich, Morita \cite{MoPJA} \cite{MoICM} introduced the
Lie algebra $\mathcal{H}_+=\ell_g^+$ as a target of the higher Johnson
homomorphisms.\par

Let $H:=H_1(\Sigma_{g,1};\mathbb{Q})$ be the first homology group of $\Sigma_{g,1}$, and
consider $\widehat{T}:=\prod_{m=0}^{\infty}H^{\otimes m}$, the completed tensor algebra
generated by $H$, which has a natural filtration $\{\widehat{T}_p\}_{p \geq 0}$ defined by $\widehat{T}_p := \prod_{m=p}^{\infty}H^{\otimes m}$. 
Via the intersetion pairing $(\ \cdot \ )\colon H\times H \to \mathbb{Q}$,
which is skew-symmetric and non-degenerate,
we identify $H$ and its dual $H^*={\rm Hom}(H,\mathbb{Q})$:
$H\stackrel{\cong}{\to} H^*, X\mapsto (Y\mapsto (Y\cdot X))$. Let $\omega \in H^{\otimes 2}$
be the two tensor corresponding to $-1_H\in {\rm Hom}(H,H)=H^* \otimes H=H^{\otimes 2}$. In other words, if $\{A_i, B_i\}^g_{i=1} \subset H$ is a symplectic basis of the first homology group $H$, then we have $\omega = \sum^g_{i=1}
A_iB_i - B_iA_i \in H^{\otimes 2}$. Here and for the rest of this paper, we omit the symbol $\otimes$ because we regard it as the product operation on the completed tensor algebra $\widehat{T}$. 
By definition, the {\it Lie algebra of symplectic derivations} is
$\mathfrak{a}_g^-={\rm Der}_{\omega}(\widehat{T})$, i.e., the Lie algebra of (continuous)
derivations on the algebra $\widehat{T}$ annihilating $\omega$. The restriction
$$
\mathfrak{a}_g^- \to {\rm Hom}(H,\widehat{T})=H^* \otimes \widehat{T}=H\otimes \widehat{T}
=\widehat{T}_1, \quad D\mapsto D|_H$$
is injective. The image is described as follows. Define a $\mathbb{Q}$-linear map $N\colon \widehat{T} \to \widehat{T}$ by 
$$
N\vert_{H^{\otimes m}} = \sum^{m-1}_{p=0} \nu^p, \quad
\mbox{for $m \geq 1$,}
$$
where $\nu$ is the cyclic permutation given by $X_1X_2\cdots X_m \mapsto X_2X_3\cdots X_1$ ($X_i \in H$), and $N\vert_{H^{\otimes 0}} = 0$. As was shown in \cite{KK1}, we can identify
$$
\mathfrak{a}_g^- = \prod^\infty_{m=1}N(H^{\otimes m}) \subset H\otimes \widehat{T}
$$
by the restriction map stated above. 
In particular, the Lie algebra $\mathfrak{a}_g^-$ is naturally graded. We say
$D\in \mathfrak{a}_g^-$ is of degree $m$ if $D\in N(H^{\otimes m})$. 
Then $D$ is of degree $m-2$ as a derivation of the graded algebra $\widehat{T}$. Moreover the Lie bracket on $\mathfrak{a}_g^-$ is homogeneous of degree $-2$. 
Now the algebra $\widehat{T}$ has the complete coproduct $\Delta$ given by
$\Delta(X)=X\widehat{\otimes}1+1\widehat{\otimes}X$, $X\in H$. Let $\mathfrak{l}_g^+$
be the Lie subalgebra of $\mathfrak{a}_g^-$ consisting of the derivations of 
degree $\geq 3$ and stabilizing the coproduct on $\widehat{T}$. The Lie algebra $\mathfrak{l}_g^+$
is an ideal of (the completion of) Kontsevich's Lie $\ell_g$ \cite{Kon}.
\par
We denote $\pi := \pi_1(\Sigma_{g,1}, *)$. Let  $\theta\colon \pi \to \widehat{T}$ be a {\it symplectic expansion}. By definition, $\theta$ is a group homomorphism of $\pi$ into the group-like elements of $\widehat{T}$, $\theta(x) \equiv 1 + [x] \pmod{\widehat{T}_2}$ for any $x \in \pi$, and $\theta(\zeta) = \exp(\omega)$
\cite{Mas}. Here $[x]\in H$ is the homology class of $x$ and
$\zeta \in \pi$ is a boundary loop on $\partial \Sigma_{g,1}$ in the opposite direction. See Figure 13.

\begin{figure}
\begin{center}
\caption{the boundary loop and symplectic generators ($g=2$)}

\vspace{0.2cm}

\unitlength 0.1in
\begin{picture}( 30.0000, 14.5000)(  2.0000,-15.1400)
%
{\color[named]{Black}{%
\special{pn 13}%
\special{ar 3040 874 160 640  0.0000000 6.2831853}%
}}%
%
{\color[named]{Black}{%
\special{pn 13}%
\special{ar 840 874 640 640  1.5707963 4.7123890}%
}}%
%
{\color[named]{Black}{%
\special{pn 13}%
\special{ar 840 874 240 240  0.0000000 6.2831853}%
}}%
%
{\color[named]{Black}{%
\special{pn 13}%
\special{ar 1960 874 240 240  0.0000000 6.2831853}%
}}%
%
{\color[named]{Black}{%
\special{pn 8}%
\special{ar 840 874 400 400  1.5707963 6.2831853}%
}}%
%
{\color[named]{Black}{%
\special{pn 8}%
\special{ar 1960 874 400 400  1.5707963 6.2831853}%
}}%
%
{\color[named]{Black}{%
\special{pn 8}%
\special{ar 1640 874 400 400  1.5707963 3.1415927}%
}}%
%
{\color[named]{Black}{%
\special{pn 8}%
\special{ar 2760 874 400 400  1.5707963 3.1415927}%
}}%
%
{\color[named]{Black}{%
\special{pn 8}%
\special{ar 1480 874 400 400  1.5707963 3.1415927}%
}}%
%
{\color[named]{Black}{%
\special{pn 8}%
\special{ar 2600 874 400 400  1.5707963 3.1415927}%
}}%
%
{\color[named]{Black}{%
\special{pn 8}%
\special{ar 1800 874 400 400  1.5707963 3.1415927}%
}}%
%
{\color[named]{Black}{%
\special{pn 8}%
\special{ar 2920 874 400 400  1.5707963 3.1415927}%
}}%
%
{\color[named]{Black}{%
\special{pn 8}%
\special{pa 2200 874}%
\special{pa 2200 714}%
\special{pa 2200 682}%
\special{pa 2200 650}%
\special{pa 2198 618}%
\special{pa 2198 586}%
\special{pa 2200 554}%
\special{pa 2204 522}%
\special{pa 2208 490}%
\special{pa 2216 458}%
\special{pa 2224 426}%
\special{pa 2234 396}%
\special{pa 2248 366}%
\special{pa 2264 338}%
\special{pa 2282 312}%
\special{pa 2304 288}%
\special{pa 2326 266}%
\special{pa 2350 244}%
\special{pa 2360 234}%
\special{sp 0.070}%
}}%
%
{\color[named]{Black}{%
\special{pn 8}%
\special{pa 2520 874}%
\special{pa 2520 714}%
\special{pa 2522 682}%
\special{pa 2522 650}%
\special{pa 2522 618}%
\special{pa 2522 586}%
\special{pa 2520 554}%
\special{pa 2518 522}%
\special{pa 2512 490}%
\special{pa 2506 458}%
\special{pa 2498 426}%
\special{pa 2486 396}%
\special{pa 2472 366}%
\special{pa 2456 338}%
\special{pa 2438 312}%
\special{pa 2418 288}%
\special{pa 2396 266}%
\special{pa 2372 244}%
\special{pa 2360 234}%
\special{sp}%
}}%
%
{\color[named]{Black}{%
\special{pn 8}%
\special{pa 1080 874}%
\special{pa 1080 714}%
\special{pa 1080 682}%
\special{pa 1080 650}%
\special{pa 1078 618}%
\special{pa 1078 586}%
\special{pa 1080 554}%
\special{pa 1084 522}%
\special{pa 1088 490}%
\special{pa 1096 458}%
\special{pa 1104 426}%
\special{pa 1114 396}%
\special{pa 1128 366}%
\special{pa 1144 338}%
\special{pa 1162 312}%
\special{pa 1184 288}%
\special{pa 1206 266}%
\special{pa 1230 244}%
\special{pa 1240 234}%
\special{sp 0.070}%
}}%
%
{\color[named]{Black}{%
\special{pn 8}%
\special{pa 1400 874}%
\special{pa 1400 714}%
\special{pa 1402 682}%
\special{pa 1402 650}%
\special{pa 1402 618}%
\special{pa 1402 586}%
\special{pa 1400 554}%
\special{pa 1398 522}%
\special{pa 1392 490}%
\special{pa 1386 458}%
\special{pa 1378 426}%
\special{pa 1366 396}%
\special{pa 1352 366}%
\special{pa 1336 338}%
\special{pa 1318 312}%
\special{pa 1298 288}%
\special{pa 1276 266}%
\special{pa 1252 244}%
\special{pa 1240 234}%
\special{sp}%
}}%
%
{\color[named]{Black}{%
\special{pn 8}%
\special{pa 840 1274}%
\special{pa 2920 1274}%
\special{fp}%
}}%
%
{\color[named]{Black}{%
\special{pn 13}%
\special{pa 840 1514}%
\special{pa 3036 1514}%
\special{fp}%
}}%
%
{\color[named]{Black}{%
\special{pn 13}%
\special{pa 840 234}%
\special{pa 3036 234}%
\special{fp}%
}}%
%
{\color[named]{Black}{%
\special{pn 4}%
\special{sh 1}%
\special{ar 2920 1274 26 26 0  6.28318530717959E+0000}%
\special{sh 1}%
\special{ar 2920 1274 26 26 0  6.28318530717959E+0000}%
}}%
%
{\color[named]{Black}{%
\special{pn 8}%
\special{pa 1380 634}%
\special{pa 1400 554}%
\special{fp}%
\special{pa 1400 554}%
\special{pa 1420 634}%
\special{fp}%
}}%
%
{\color[named]{Black}{%
\special{pn 8}%
\special{pa 2500 634}%
\special{pa 2520 554}%
\special{fp}%
\special{pa 2520 554}%
\special{pa 2540 634}%
\special{fp}%
}}%
%
{\color[named]{Black}{%
\special{pn 8}%
\special{pa 780 454}%
\special{pa 860 474}%
\special{fp}%
\special{pa 860 474}%
\special{pa 780 494}%
\special{fp}%
}}%
%
{\color[named]{Black}{%
\special{pn 8}%
\special{pa 1900 454}%
\special{pa 1980 474}%
\special{fp}%
\special{pa 1980 474}%
\special{pa 1900 494}%
\special{fp}%
}}%
%
{\color[named]{Black}{%
\special{pn 8}%
\special{pa 2920 1378}%
\special{pa 2968 1444}%
\special{fp}%
\special{pa 2968 1444}%
\special{pa 2958 1364}%
\special{fp}%
}}%
\put(4.8800,-5.0200){\makebox(0,0)[lb]{$\alpha_1$}}%
\put(11.9600,-2.0200){\makebox(0,0)[lb]{$\beta_1$}}%
\put(17.2000,-4.4200){\makebox(0,0)[lb]{$\alpha_2$}}%
\put(23.0000,-1.9400){\makebox(0,0)[lb]{$\beta_2$}}%
\put(29.4800,-11.9800){\makebox(0,0)[lb]{$*$}}%
\put(27.8000,-14.4600){\makebox(0,0)[lb]{$\zeta$}}%
%
{\color[named]{Black}{%
\special{pn 8}%
\special{ar 2720 386 80 80  6.2831853 6.2831853}%
\special{ar 2720 386 80 80  0.0000000 4.7123890}%
}}%
%
{\color[named]{Black}{%
\special{pn 8}%
\special{pa 2800 386}%
\special{pa 2752 410}%
\special{fp}%
\special{pa 2800 386}%
\special{pa 2828 430}%
\special{fp}%
}}%
\end{picture}%

\end{center}
\end{figure}
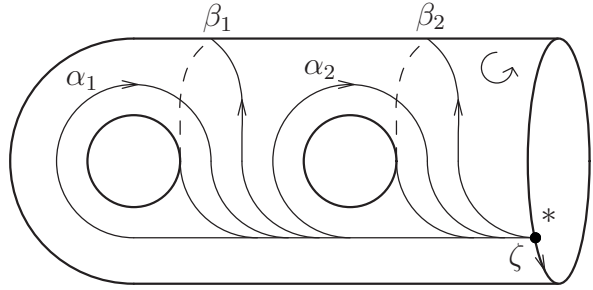

Any symplectic expansion induces an isomorphism $\theta\colon \widehat{\mathbb{Q}\pi} \overset\cong\to \widehat{T}$ of complete Hopf algebras. Here $\widehat{\mathbb{Q}\pi}$ is the completed group ring of $\pi$
(see \S \ref{CGL}). Moreover, in \cite{KKp} Theorem 6.1.4 and 6.1.5, we showed that the map 
\begin{equation}
\label{GLA-ASS}
-\lambda_{\theta}\colon \widehat{\mathbb{Q}\hat{\pi}}(\Sigma_{g,1})
\stackrel{\cong}{\to} \mathfrak{a}_g^-,
\quad \vert x\vert \mapsto -N\theta(x)
\end{equation}
is a filter preserving isomorphism of Lie algebras, and induces 
a filter preserving isomorphism of Lie algebras
$$
-\lambda_\theta\colon L^+(\Sigma_{g,1},\{ *\})\overset\cong\to \mathfrak{l}^+_g.
$$
As was mentioned in \cite{KKp} \S 6.3,
the composite $-\lambda_\theta\circ \tau$ is exactly Massuyeau's improvement 
$\rho^\theta$ \cite{Mas} of the
Johnson map introduced by Kawazumi \cite{Ka}. Its graded quotients
with respect to the Johnson filtration on $\mathcal{I}_{g,1}$ and the degree filtration on $\mathfrak{l}^+_g$ are the Johnson homomorphisms
of all degrees introduced by Johnson \cite{JS} and improved by Morita \cite{MoAQ}.
Indeed, it is this context in which the Lie algebra $\ell_g$ was introduced by Morita \cite{MoPJA} \cite{MoICM}.
\par
Through the isomorphism $-\lambda_\theta$, the Turaev cobracket $\delta$ on 
$\widehat{\mathbb{Q}\hat\pi}(\Sigma_{g,1})$ induces the complete Lie cobracket
$\delta^\theta$ on the Lie algebra $\mathfrak{a}_g^-$.
Namely, $\delta^{\theta}$ is defined so that the following diagram commutes:
$$
\begin{CD}
\widehat{\mathbb{Q}\hat\pi}(\Sigma_{g,1}) @>{\delta}>> \widehat{\mathbb{Q}\hat\pi}(\Sigma_{g,1})\widehat{\otimes}\widehat{\mathbb{Q}\hat\pi}(\Sigma_{g,1})\\
@V{-\lambda_\theta}VV @V{(-\lambda_\theta)^{\widehat{\otimes}2}}VV\\
\mathfrak{a}_g^- @>{\delta^\theta}>> \mathfrak{a}_g^-\widehat{\otimes}\mathfrak{a}_g^-. 
\end{CD}
$$
The grading on $\mathfrak{a}_g^-$ defines the Laurent expansion of the cobracket $\delta^\theta$.
Namely, for any $u\in H^{\otimes m}$ we can write
\begin{eqnarray*}
&&\delta^\theta(N(u)) = \sum^\infty_{p= -\infty}\delta^\theta_{(p)}
(N(u)), \quad {\rm where} \\
&&\delta^\theta_{(p)}
(N(u)) \in (\mathfrak{a}_g^-\widehat{\otimes}\mathfrak{a}_g^-)_{(m+p)} := \bigoplus_{k+l = m+p} N(H^{\otimes k}) \otimes N(H^{\otimes l}).
\end{eqnarray*}
In this subsection and \S5.5 we prove the following
\begin{thm}\label{laurent}
 For any symplectic expansion $\theta$ we have 
\begin{enumerate}
\item $\delta^\theta_{(p)} = 0$ for $p = 0, -1$, and $p \leq -3$.
\item $\delta^\theta_{(-2)}$ is the same as Schedler's cobracket {\rm \cite{Sch}}, i.e., 
\begin{eqnarray*}
&& \delta^\theta_{(-2)}(N(X_1X_2\cdots X_m)) \\
&=& -\sum_{i<j}(X_i\cdot X_j)
\{N(X_{i+1}\cdots X_{j-1})\widehat{\otimes}
N(X_{j+1}\cdots X_mX_1\cdots X_{i-1})\\
&& \quad\quad - 
N(X_{j+1}\cdots X_mX_1\cdots X_{i-1})
\widehat{\otimes}
N(X_{i+1}\cdots X_{j-1})\}
\end{eqnarray*}
for any $X_i \in H$ and $m \geq 1$. 
\end{enumerate}
\end{thm}
In particular, the $(-2)$-nd term $\delta^\theta_{(-2)}$, or equivalently Schedler's cobracket, is independent of the choice of a symplectic expansion $\theta$, so that we denote it by $\delta^{\rm alg}$.
Since $\delta^{\rm alg}$ is the graded quotient of the cobracket $\delta^\theta$ with respect to the degree filtration, it induces a structure of a complete involutive Lie bialgebra on the Lie algebra $\mathfrak{a}_g^-$. \par
Recall that the graded quotients of the geometric Johnson map are the classical Johnson homomorphisms. Hence, 
as a corollary of Theorem \ref{laurent}, we obtain
\begin{cor}\label{53johnson}
 For any $k \geq 1$, we have 
$$
\delta^{\rm alg} \circ \tau_k = 0\colon \mathcal{I}_{g,1}(k)/\mathcal{I}_{g,1}(k+1) \to (\mathfrak{a}_g^-\widehat{\otimes}\mathfrak{a}_g^-)_{(k)}.
$$
Here $\mathcal{I}_{g,1}(k)$ is the $k$-th term of the Johnson filtration, and $\tau_k$ is the $k$-th Johnson homomorphism. 
\end{cor}
\par
The proof of Theorem \ref{laurent} is based on a theorem of Massuyeau and Turaev \cite{MT} Theorem 10.4, 
which gives a tensorial description of the homotopy intersection form.
As is announced in \cite{MT2} Remark 7.4.3,
Massuyeau and Turaev \cite{MT3} also prove Theorem \ref{laurent} in a similar way to ours.
In this subsection we prove Theorem \ref{laurent} except for the case $p = 0$, which will be proved in \S 5.5. \par

Recall the definition of the homotopy intersection form \cite{MT}.
Taking a path $\nu$ and a second base point $\bullet$ as in Figure 7, we identify
the fundamental group $\pi_1(\Sigma_{g,1}, \bullet)$ with $\pi=\pi_1(\Sigma_{g,1},*)$
by the isomorphism $\alpha \mapsto \overline{\nu} \alpha\nu$.
(Note that in \cite{MT} \S7, $\nu$ and $\overline{\nu}$ are denoted by
$\nu_{\bullet\ast}$ and $\overline{\nu}_{\ast\bullet}$, respectively.)
Let $\alpha$ be an oriented based immersed loop on $\Sigma_{g,1}$ with base point $\bullet$, and $\beta$ an oriented based immersed loop on $\Sigma_{g,1}$ with base point $\ast$ such that their intersections consists of transverse double points. Then the formula 
\begin{equation}\label{53eta}
\eta(\alpha, \beta)  := \sum_{p\in \alpha\cap\beta}\varepsilon(p;\alpha, \beta) \overline{\nu}
\alpha_{\bullet p}\beta_{p\ast} \in \mathbb{Q}\pi
\end{equation}
gives rise to a well-defined $\mathbb{Q}$-bilinear map $\eta\colon \mathbb{Q}\pi\times\mathbb{Q}\pi \to \mathbb{Q}\pi$ via the identification $\pi_1(\Sigma_{g,1}, \bullet) = \pi$ stated above. The map $\eta$ is called the {\it homotopy intersection form} of $\Sigma_{g,1}$ \cite{MT}, 
which is essentially the same as what Papakyriakopoulos \cite{Pa} and Turaev \cite{T1} independently introduced. \par
Massuyeau and Turaev \cite{MT} proves that this map $\eta$ naturally extends to a $\mathbb{Q}$-bilinear map $\hat\eta\colon \widehat{\mathbb{Q}\pi}\times\widehat{\mathbb{Q}\pi} \to \widehat{\mathbb{Q}\pi}$, and gives its tensorial description through any symlectic expansion
$\theta\colon  \widehat{\mathbb{Q}\pi}\overset\cong\to \widehat{T}$. 
Let $\varepsilon\colon \widehat{T} \to \widehat{T}/\widehat{T}_1 = \mathbb{Q}$ be the augmentation map. Define a $\mathbb{Q}$-bilinear map $\overset\bullet\leadsto\colon \widehat{T}_1 \times \widehat{T}_1 \to \widehat{T}$ by 
$$
\left(X_1\cdots X_m\overset\bullet\leadsto Y_1\cdots Y_n\right) := (X_m\cdot Y_1) X_1\cdots X_{m-1}Y_2\cdots Y_n \in H^{\otimes m+n-2}
$$
for any $n$, $m \geq 1$, and $X_i$, $Y_j \in H$. Here $(X_m\cdot Y_1) \in \mathbb{Q}$ is the intersection pairing of $X_m$ and $Y_1 \in H$. 
A $\mathbb{Q}$-bilinear map $\rho\colon \widehat{T}\times\widehat{T} \to \widehat{T}$ is defined by 
\begin{equation}
\rho(a, b) := (a -\varepsilon(a))\overset\bullet\leadsto(b -\varepsilon(b)) 
+ (a -\varepsilon(a))s(\omega)(b -\varepsilon(b))
\label{53rho}
\end{equation}
for any $a$ and $b \in \widehat{T}$, where $s(z)$ is the formal power series
$$
s(z) = \frac{1}{e^{-z}-1} + \frac{1}{z}
= -\frac{1}{2} - \sum_{k\geq 1}\frac{B_{2k}}{(2k)!}z^{2k-1} = 
 -\frac{1}{2}  -\frac{z}{12}  +\frac{z^3}{720}  -\frac{z^5}{30240} + \cdots.
$$
Massuyeau and Turaev \cite{MT} proved the following. 
\begin{thm}[\cite{MT} Theorem 10.4]\label{53MT}
Let $\theta\colon \pi \to \widehat{T}$ be a symplectic expansion. Then the following diagram commutes: 
$$
\begin{CD}
\widehat{\mathbb{Q}\pi}\times \widehat{\mathbb{Q}\pi}
@>{\hat\eta}>> \widehat{\mathbb{Q}\pi}\\
@V{\theta\times\theta}VV @V{\theta}VV\\
\widehat{T}\times\widehat{T} @>{\rho}>> \widehat{T}.
\end{CD}
$$
\end{thm}
\par 
Recall the $\mathbb{Q}$-linear map $\kappa\colon \mathbb{Q}\pi \otimes \mathbb{Q}\pi \to \mathbb{Q}\pi \otimes \mathbb{Q}\pi$ introduced in \S \ref{Intpath}.
Let $\Delta\colon \mathbb{Q}\pi\to \mathbb{Q}\pi\otimes \mathbb{Q}\pi$ be the coproduct defined by $\alpha \in \pi \mapsto \alpha\otimes\alpha$, and $\iota\colon \mathbb{Q}\pi \to \mathbb{Q}\pi$ the antipode defined by 
$\alpha \in \pi \mapsto \alpha^{-1}$. Then, 
for any $\alpha$ and $\beta \in \pi$ in general position as stated above, we have 
\begin{eqnarray}
&&-(1\otimes\beta)((1\otimes\iota)\Delta\eta(\alpha,\beta))(1\otimes\alpha)\nonumber\\
&=&-\sum_{p \in\alpha\cap\beta}\varepsilon(p; \alpha, \beta)(1\otimes\beta)(\overline{\nu}
\alpha_{\bullet p}\beta_{p\ast}\otimes (\beta_{p\ast})^{-1}(\alpha_{\bullet p})^{-1}\nu)
(1\otimes \overline{\nu}\alpha\nu)\nonumber\\
&=& -\sum_{p \in\alpha\cap\beta}\varepsilon(p; \alpha, \beta)(\overline{\nu}
\alpha_{\bullet p}\beta_{p\ast}\otimes \beta_{\ast p}\alpha_{p\bullet}{\nu})\nonumber\\
&=&\kappa(\alpha,\beta).
\nonumber
\end{eqnarray}
Hence we obtain
\begin{equation}
\kappa(u,v) = - \sum(1\otimes v'')((1\otimes\iota)(\Delta\eta(u', v')))(1\otimes u'')
\label{53kappa}
\end{equation}
for any $u$ and $v \in \mathbb{Q}\pi$. Here we denote $\Delta u = \sum u'\otimes u''$ and $\Delta v = \sum v'\otimes v''$. \par
Let $\theta\colon \pi \to \widehat{T}$ be a symplectic expansion. 
Then we have a unique $\mathbb{Q}$-linear map $\kappa^\theta\colon \widehat{T}\widehat{\otimes}\widehat{T} \to \widehat{T}\widehat{\otimes}\widehat{T}$ such that the diagram 
$$
\begin{CD}
\widehat{\mathbb{Q}\pi}\widehat{\otimes} \widehat{\mathbb{Q}\pi}
@>{\kappa}>> \widehat{\mathbb{Q}\pi}\widehat{\otimes} \widehat{\mathbb{Q}\pi}\\
@V{\theta\widehat{\otimes}\theta}VV @V{\theta\widehat{\otimes}\theta}VV\\
\widehat{T}\widehat{\otimes}\widehat{T} @>{\kappa^\theta}>> \widehat{T}\widehat{\otimes}\widehat{T}
\end{CD}
$$
commutes. By Theorem \ref{53MT} and (\ref{53kappa}), the map $\kappa^\theta$ does not depend on the expansion $\theta$. From (\ref{53kappa})
we have
\begin{equation}
\kappa^\theta(X,Y) = -(1\widehat{\otimes}1)((1\widehat{\otimes}\iota)\Delta\rho(X, Y))(1\widehat{\otimes}1) = -(X\cdot Y)(1\widehat{\otimes}1) - (1\widehat{\otimes}\iota)\Delta(Xs(\omega)Y)
\label{53xy}
\end{equation}
for any $X$ and $Y \in H$. \par
Now we define a $\mathbb{Q}$-linear map $\mu^\theta\colon \widehat{T} \to \widehat{T}\widehat{\otimes}\mathfrak{a}_g^-$ by the commutative diagram
$$
\begin{CD}
\widehat{\mathbb{Q}\pi}
@>{\mu}>> \widehat{\mathbb{Q}\pi}\widehat{\otimes} \widehat{\mathbb{Q}\hat\pi}\\
@V{\theta}VV @V{-\theta\widehat{\otimes}\lambda_\theta}VV\\
\widehat{T} @>{\mu^\theta}>> \widehat{T}\widehat{\otimes}\mathfrak{a}_g^-.
\end{CD}
$$
From Corollary \ref{mu(1-n)} we have 
\begin{eqnarray}
&&\mu^\theta(X_1\cdots X_m)\nonumber\\
&=& (1\widehat{\otimes}(-N))\sum_{1\leq i<j \leq m}(X_1\cdots X_{i-1}\widehat{\otimes}1)\kappa^\theta(X_i, X_j)(X_{j+1}\cdots X_m\widehat{\otimes}X_{i+1}\cdots X_{j-1})\nonumber\\
&& + \sum^m_{i=1}(X_1\cdots X_{i-1}\widehat{\otimes}1)\mu^\theta(X_i)(X_{i+1}\cdots X_m\widehat{\otimes}1)\nonumber\\
&=& \sum_{1\leq i<j \leq m}(X_i\cdot X_j)X_1\cdots X_{i-1}X_{j+1}\cdots X_m\widehat{\otimes}N(X_{i+1}\cdots X_{j-1}) \nonumber\\
&& + (1\widehat{\otimes}N)\sum_{1\leq i<j \leq m}(X_1\cdots X_{i-1}\widehat{\otimes}1)((1\widehat{\otimes}\iota)\Delta(X_is(\omega)X_j))(X_{j+1}\cdots X_m\widehat{\otimes}X_{i+1}\cdots X_{j-1})\nonumber\\
&& + \sum^m_{i=1}(X_1\cdots X_{i-1}\widehat{\otimes}1)\mu^\theta(X_i)(X_{i+1}\cdots X_m\widehat{\otimes}1) \label{53mu}
\end{eqnarray}
for any $m \geq 0$ and $X_i \in H$. We remark that the first term in (\ref{53mu}) is of degree $m-2$, and the second and the third terms are of degree $\geq m$. 
In fact, $\mu^\theta(X_i) \in \widehat{T}\widehat{\otimes}\mathfrak{a}_g^-$ and  $ \widehat{T}\widehat{\otimes}\mathfrak{a}_g^-$ starts from degree $1$.
On the other hand, by Lemma \ref{mu-delta} the maps $\delta^{\theta}$ and $\mu^{\theta}$ are related by the formula
\begin{equation}
\label{mu-del-th}
\delta^{\theta}\circ N=(1-T)(N\widehat{\otimes} 1_{\mathfrak{a}_g^-})\mu^{\theta},
\end{equation}
where $T\colon \mathfrak{a}_g^- \widehat{\otimes} \mathfrak{a}_g^-
\to \mathfrak{a}_g^- \widehat{\otimes} \mathfrak{a}_g^-$ is the switch map:
$T(u\widehat{\otimes}v)=v\widehat{\otimes}u$.
Theorem \ref{laurent} except for the case $p=0$ follows from the above observation and (\ref{mu-del-th}).
\qed

\subsection{The Morita traces}

In this subsection we prove that Schedler's cobracket $\delta^{\rm alg}$ restricted to $\mathfrak{l}_g^+$ covers the {\it Morita traces} of all degrees ${\rm Tr}_k\colon (\mathfrak{l}_g^+)_{(k+1)} \to {\rm Sym}^{k-1}H$, $k \geq 4$ \cite{MoAQ}.
Here $(\mathfrak{l}_g^+)_{(n)}$ is the degree $n$ part of $\mathfrak{l}_g^+ \subset \mathfrak{a}_g^-$, and ${\rm Sym}^{n}H$ is the $n$-th symmetric power of the homology group $H = H_1(\Sigma_{g,1}; \mathbb{Q})$. \par
To state our result precisely, we need some notations. Let $p_1\colon \mathfrak{a}_g^- = \prod^\infty_{m=1} N(H^{\otimes m}) \to N(H^{\otimes 1}) = H$ be the first projection, $i\colon \mathfrak{a}_g^- = \prod^\infty_{m=1} N(H^{\otimes m}) \hookrightarrow \prod^\infty_{m=1} H^{\otimes m} = \widehat{T}_1$ the inclusion map, and $\varpi\colon \widehat{T} \to \widehat{\rm Sym}(H) := \prod^\infty_{m=0}{\rm Sym}^m(H)$ the natural projection. We define 
$$
\mathfrak{s} := \varpi\circ (p_1\widehat{\otimes}i)\colon \mathfrak{a}_g^-\widehat{\otimes}\mathfrak{a}_g^- \to H\otimes \widehat{T}_1 = \widehat{T}_2 \to \widehat{\rm Sym}(H).
$$
Then we have 
\begin{thm}\label{54trace}
$$
\mathfrak{s}\circ\delta^{\rm alg}\vert_{(\mathfrak{l}_g^+)_{(m+2)}} = (-1)^mm\times{\rm Tr}_{m+1}\colon(\mathfrak{l}_g^+)_{(m+2)} \to {\rm Sym}^mH
$$
for any $m \geq 3$. 
\end{thm}

From Corollary \ref{53johnson} $\mathfrak{s}\circ\delta^{\rm lag}\circ \tau_k = 0\colon \mathcal{I}_{g,1}(k)/\mathcal{I}_{g,1}(k+1) \to {\rm Sym}^kH$. Hence Theorem \ref{54trace} means that all the Morita traces are derived from the fundamental geometric fact that any diffeomorphism preserves the self-intersection of any curve on the surface. \par

The rest of this subsection is devoted to the proof of Theorem \ref{54trace}. 
Let $\widehat{\mathcal{L}} \subset \widehat{T}_1$ be the completed free Lie algebra over the vector space $H$. It is known that (a multiple of) the Dynkin idempotent $\Phi\colon \widehat{T}_1 \to \widehat{\mathcal{L}}$ defined by 
$\Phi(X_1X_2\cdots X_m) := [X_1, [X_2, [\cdots [X_{m-1}, X_m]]]]$ for $X_i \in H$ and $m \geq 1$ satisfies $\Phi\vert_{\widehat{\mathcal{L}}\cap H^{\otimes m}} = m1_{\widehat{\mathcal{L}}\cap H^{\otimes m}}$. On the other hand,  
as was shown in \S2.7 \cite{KK1},
we have $\mathfrak{l}_g^+\rtimes \mathfrak{sp}(H) = N(H\otimes \widehat{\mathcal{L}})
={\rm Ker}([\, , \,]) \subset H\otimes \widehat{L}$.
Here $[\, , \,]: H\otimes \widehat{\mathcal{L}} \to \widehat{\mathcal{L}}$,
$X\otimes u \mapsto [X, u]$, is the bracket map.
Since $N(Y[X_1, \Phi(X_2\cdots X_m)]) = N([Y, X_1]\Phi(X_2\cdots X_m))$ for any $Y \in H$, we have
\begin{equation}
N([H, H]\otimes \widehat{\mathcal{L}}) = \mathfrak{l}^+_g.
\label{54HH}
\end{equation}
Hence it suffices to prove Theorem \ref{54trace} at $N([Y,Z]\Phi(X_1\cdots X_m))$ for any $Y, Z, X_i \in H$. \par

Originally the Morita traces were defined as the trace of certain matrix representation
of ${\rm Hom}(H_{\mathbb{Z}},\Gamma_k/\Gamma_{k+1})$ using the Fox free derivative.
Here $H_{\mathbb{Z}}=\pi^{\rm ab}=H_1(\Sigma_{g,1};\mathbb{Z})$ and $\{ \Gamma_k\}_k$ is the
the lower central series of $\pi$ (with $\Gamma_1=\pi$). Actually it is known
that the $k$-th Morita trace ${\rm Tr}_k: H\otimes \mathcal{L}_k \to
{\rm Sym}^{k-1}(H)$ coincides with $(-1)^k$ times the map $H\otimes 
\mathcal{L}_k \subset H\otimes H^{\otimes k} \overset{C_{12}}\to 
H^{\otimes (k-1)} \to  {\rm Sym}^{k-1}(H)$, where $C_{12}: H\otimes 
H^{\otimes k} \to H^{\otimes (k-1)}$ is defined by $X_0\otimes X_1X_2\cdots X_k
\mapsto (X_0\cdot X_1)X_2\cdots X_k$ for any $X_i \in H$, cf. \cite{MoProblems} Remark 22
(note that the convention about indices in \cite{MoAQ} \cite{MoProblems} is different from ours).
We use this description of ${\rm Tr}_k$. To compute ${\rm Tr}_{m+1}(N([Y,Z]\Phi(X_1\cdots X_m))) \in {\rm Sym}^m(H)$,
note that any tensor including $[Y, Z]$ or $\Phi(X_1\cdots X_m)$ must vanish
in the symmetric power ${\rm Sym}^m(H)$. Hence, if we fix $Y$ and $Z$ and define a map
$\beta = \beta_{Y,Z}\colon H^{\otimes m} \to {\rm Sym}^m(H)$ by 
\begin{eqnarray*}
\beta(X_1\cdots X_m) &:=& -(Y\cdot X_1)ZX_2\cdots X_m 
-(Y\cdot X_m)ZX_1\cdots X_{m-1} \\
&&+ (Z\cdot X_1)YX_2\cdots X_m 
+ (Z\cdot X_m) YX_1\cdots X_{m-1},
\end{eqnarray*}
then we have
$$
{\rm Tr}_{m+1}(N([Y,Z]\Phi(X_1\cdots X_m))) =(-1)^{m+1} \beta(\Phi(X_1\cdots X_m))
$$
for any $X_i \in H$.

\begin{lem}\label{54beta}
\begin{enumerate}
\item If $m$ is even, then $\beta(\Phi(X_1\cdots X_m)) = 0$. 
\item If $m \geq 3$, then $\beta(\Phi(X_1X_2X_3\cdots X_m)) = X_1X_2\beta(\Phi(X_3\cdots X_m))$. 
\end{enumerate}
\end{lem}
The first assertion was already known by Morita \cite{MoAQ} Theorem 6.1 (ii). 
\begin{proof}
(1) Let $\iota\colon \widehat{T} \to \widehat{T}$ denote the antipode, i.e., we have 
$\iota(X_1\cdots X_m) = (-1)^mX_m\cdots X_1$. The map $\beta$ is `symmetric'. In other words, we have $\beta\iota(X_1\cdots X_m) = (-1)^m\beta(X_m\cdots X_1) = (-1)^m\beta(X_1\cdots X_m) \in {\rm Sym}^m(H)$, while we have $\iota\vert_{\widehat{\mathcal{L}}} = -1$. Hence $-\beta\Phi(X_1\cdots X_m) = \beta\iota\Phi(X_1\cdots X_m) = (-1)^m\beta\Phi(X_1\cdots X_m)$, which implies $\beta\Phi(X_1\cdots X_m) = 0$ if $m$ is even.
\par
(2) We remark
\begin{eqnarray}
\Phi(X_1X_2X_3\cdots X_m) &=& X_1X_2\Phi(X_3\cdots X_m) + \Phi(X_3\cdots X_m)X_2X_1\nonumber\\
&& - X_1\Phi(X_3\cdots X_m)X_2 - X_2\Phi(X_3\cdots X_m)X_1.
\label{54Phi}
\end{eqnarray}
Then $\beta(X_1\Phi(X_3\cdots X_m)X_2) = \beta(X_1\Phi(X_3\cdots X_m)X_2) = 0$, since $\Phi(X_3\cdots X_m)$ remains in ${\rm Sym}^m(H)$. 
When we compute the $\beta$-image of the first and the second terms in the right hand side, the terms coming from the contraction of $X_1$ and $Y$ or $Z$ must vanish since they include $\Phi(X_3\cdots X_m)$. The rest terms equal $X_1X_2\beta\Phi(X_3\cdots X_m)$. This proves the lemma. 
\end{proof}

From the definition, we have
$$
\mathfrak{s}\delta^{\rm alg}(N(X_1X_2\cdots X_m)) = 
\sum_{i=1}^m (X_{i-1}\cdot X_{i+1})X_i X_{i+2}\cdots X_m X_1\cdots X_{i-2}
\in {\rm Sym}^{m-2}(H),
$$
where we denote $X_{m+1} = X_1$, $X_{0} = X_m$ and so on. 
Similarly $[Y,Z]$ and $\Phi(X_1\cdots X_m)$ vanish in the symmetric power ${\rm Sym}^m(H)$. Hence, if we fix $Y$ and $Z$ and define maps $\alpha = \alpha_{Y,Z}\colon H^{\otimes m} \to {\rm Sym}^m(H)$ 
and $\gamma = \gamma_{Y,Z}\colon H^{\otimes m} \to {\rm Sym}^m(H)$
by 
\begin{eqnarray*}
\alpha(X_1\cdots X_m) &:=& - \beta(X_1\cdots X_m) + \gamma(X_1\cdots X_m), 
\quad\mbox{and}\\
\gamma(X_1\cdots X_m) &:=&
-(Y\cdot X_2)ZX_1X_3\cdots X_{m} 
-(Y\cdot X_{m-1})ZX_1\cdots X_{m-2}X_m \\
&&+ (Z\cdot X_2)YX_1X_3\cdots X_m 
+ (Z\cdot X_{m-1}) YX_1\cdots X_{m-2}X_m,
\end{eqnarray*}
respectively, then we have
$$
\mathfrak{s}\delta^{\rm alg}(N([Y,Z]\Phi(X_1\cdots X_m))) 
= \alpha(\Phi(X_1\cdots X_m))
$$
for any $X_i \in H$. 

\begin{lem}\label{54alpha}
\begin{enumerate}
\item If $m \geq 4$, then $\gamma(X_1X_2\cdots X_{m-1}X_m) = X_1X_m\beta(X_2\cdots X_{m-1})$. 
\item If $m$ is even, then $\alpha(\Phi(X_1\cdots X_m)) =\gamma(\Phi(X_1\cdots X_m)) = 0$. 
\item If $m \geq 5$, then $$\gamma(\Phi(X_1X_2X_3\cdots X_m)) = X_1X_2\gamma(\Phi(X_3\cdots X_m))
- 2X_1X_2\beta(\Phi(X_3\cdots X_m)).$$ 
\end{enumerate}
\end{lem}
\begin{proof} (1) is clear from the definition. \par
(2) is proved in a similar way to Lemma \ref{54beta} (1), since the map $\gamma$ is also `symmetric' with respect to the antipode $\iota$. \par
(3) Recall the equation (\ref{54Phi})
\begin{eqnarray*}
\Phi(X_1X_2X_3\cdots X_m) &=& X_1X_2\Phi(X_3\cdots X_m) + \Phi(X_3\cdots X_m)X_2X_1\nonumber\\
&& - X_1\Phi(X_3\cdots X_m)X_2 - X_2\Phi(X_3\cdots X_m)X_1.
\end{eqnarray*}
By (1) the $\gamma$-image of each of the third and the fourth terms in the right hand side equals $-X_1X_2\beta(\Phi(X_3\cdots X_m))$. When we compute the $\gamma$-image of the first and the second terms, the terms coming from the contraction of $X_2$ and $Y$ or $Z$ must vanish, since they include $\Phi(X_3\cdots X_m)$. Hence the $\gamma$-image of the sum of the first and the second terms equals $X_1X_2\gamma(\Phi(X_3\cdots X_m))$. This proves the lemma. 
\end{proof}

\begin{lem}\label{54gamma}
If $m \geq 3$, then we have
$$
\gamma(\Phi(X_1X_2\cdots X_m)) = -(m-1)\beta(\Phi(X_1X_2\cdots X_m)).
$$
\end{lem}
\begin{proof} If $m$ is even, 
$\gamma(\Phi(X_1X_2\cdots X_m)) = \beta(\Phi(X_1X_2\cdots X_m)) = 0$ 
from Lemma \ref{54Phi} (1) and Lemma \ref{54alpha} (2). 
Hence it suffices to prove the lemma in the case $m$ is odd $\geq 3$ by induction.\par
For $m=3$ we compute the both sides in the lemma explicitly. 
Then we have 
\begin{eqnarray*}
&&\gamma(\Phi(X_1X_2X_3))\\
 &=& -4(Y\cdot X_2)ZX_1X_3 + 4(Z\cdot X_2)YX_1X_3 + 4(Y\cdot X_3)ZX_1X_2 - 4(Z\cdot X_3)YX_1X_2\\
&=& -2\beta(\Phi(X_1X_2X_3)),
\end{eqnarray*}
as was to be shown.\par
Next suppose $m \geq 5$. From Lemma \ref{54alpha} (1) we have 
$$
\gamma(\Phi(X_1X_2X_3\cdots X_m)) = X_1X_2\gamma(\Phi(X_2\cdots X_m)) - 2X_1X_2\beta(\Phi(X_3\cdots X_m)),
$$ 
which equals 
$$
-(m-3)X_1X_2\beta(\Phi(X_3\cdots X_m)) - 2X_1X_2\beta(\Phi(X_3\cdots X_m)) = -(m-1)X_1X_2\beta(\Phi(X_3\cdots X_m))
$$ 
by the inductive assumption. By Lemma \ref{54beta} (2) we have $\beta(\Phi(X_1X_2X_3\cdots X_m)) = X_1X_2\beta(\Phi(X_3\cdots X_m))$. Hence we obtain 
$$
\gamma(\Phi(X_1X_2X_3\cdots X_m)) = -(m-1)\beta(\Phi(X_1X_2X_3\cdots X_m)).
$$
This completes the induction and the proof of the lemma. 
\end{proof}

As a corollary of Lemma \ref{54gamma}, we have 
$$\alpha(\Phi(X_1\cdots X_m)) = -\beta(\Phi(X_1\cdots X_m)) + \gamma(\Phi(X_1\cdots X_m)) = -m\beta(\Phi(X_1\cdots X_m)).
$$ This completes the proof of Theorem \ref{54trace}. \qed

\subsection{The $0$-th term of the Laurent expansion of the Turaev cobracket}

In this subsection we prove Theorem \ref{laurent} for $p = 0$.

Using the grading on $\mathfrak{a}_g^-$, for any $u\in H^{\otimes m}$ we can write
\begin{eqnarray*}
&& \mu^{\theta}(u)=\sum_{p=-\infty}^{\infty} \mu^{\theta}_{(p)}(u), \quad {\rm where} \\
&& \mu^{\theta}_{(p)}(u)\in (\widehat{T} \widehat{\otimes} \mathfrak{a}_g^-)_{(m+p)}
:= \bigoplus_{k+l=m+p}H^{\otimes k} \otimes N(H^{\otimes l}).
\end{eqnarray*}
From (\ref{53mu}) we have $\mu^{\theta}_{(p)}=0$ for $p\le -3$ and $p=-1$,
and $\mu^{\theta}_{(-2)}$ is given by the first term of (\ref{53mu}). Since it does not
depend on the choice of $\theta$, we denote it by $\mu^{\rm alg}$. Namely,
\begin{eqnarray*}
\mu^{\theta}_{(-2)}(X_1\cdots X_m) &=& \mu^{\rm alg}(X_1\cdots X_m) \\
&=& \sum_{1\le i<j \le m}(X_i\cdot X_j) X_1\cdots X_{i-1}X_{j+1}\cdots X_m
\otimes N(X_{i+1}\cdots X_{j-1}).
\end{eqnarray*}
Therefore, by (\ref{mu-del-th}), Theorem \ref{laurent} for $p=0$ is deduced from the following

\begin{thm}
\label{mu-zero}
For any symplectic expansion $\theta$ and $u\in H^{\otimes m}$, we have
$$\mu^{\theta}_{(0)}(u)=-\frac{1}{2}(1\otimes N(u)).$$
\end{thm}

The rest of this subsection is devoted to the proof of Theorem \ref{mu-zero}.

First we prove Theorem \ref{mu-zero} for $m=1$, i.e.,
\begin{equation}
\label{mu-zero-1}
\mu^{\theta}_{(0)}(X)=-\frac{1}{2}(1\otimes X),\quad X\in H.
\end{equation}
Let $\alpha_1,\ldots,\alpha_g,\beta_1,\ldots,\beta_g\in \pi$ be symplectic generators of
$\pi=\pi_1(\Sigma_{g,1},*)$. See Figure 13. The homology classes $A_i=[\alpha_i], B_i=[\beta_i]\in H$
($1\le i\le g$) constitute a symplectic basis for $H$. By (\ref{mu_+}) we have
$\mu(\alpha_i)=0$ and $\mu(\beta_i)=1\otimes |\beta_i|^{\prime}$. From the definition of $\mu^{\theta}$,
we obtain
$$\mu^{\theta}(\theta(\alpha_i))=0 \quad {\rm and} \quad
\mu^{\theta}(\theta(\beta_i))=-(\theta \widehat{\otimes} \lambda_{\theta})(1\otimes |\beta_i|^{\prime})
=-1\otimes N\theta(\beta_i).$$
We denote by $\theta_k$ the degree $k$ part of $\theta$. Looking at the degree 1 part of the
above equation, for $1\le i\le g$, we have
\begin{eqnarray}
&& \mu^{\theta}_{(0)}(A_i)=-\mu^{\rm alg}(\theta_3(\alpha_i)), \nonumber \\
&& \mu^{\theta}_{(0)}(B_i)=-\mu^{\rm alg}(\theta_3(\beta_i))-1\otimes B_i.
\label{mu-A-B}
\end{eqnarray}

\begin{lem}
\label{LH3}
Let $\widehat{\mathcal{L}}\subset \widehat{T}$ be the set of primitive elements of $\widehat{T}$.
In other words, $\widehat{\mathcal{L}}$ is the completed free Lie algebra generated by $H$. Then
for any $u\in \widehat{\mathcal{L}}\cap H^{\otimes 3}$, we have $\mu^{\rm alg}(u)=0$.
\end{lem}

\begin{proof}
It is sufficient to prove the formula for $u=[X,[Y,Z]]$ where $X,Y,Z\in H$. Since
$[X,[Y,Z]]=[X,YZ-ZY]=XYZ-XZY-YZX+ZYX$, we compute
$$\mu^{\rm alg}([X,[Y,Z]])=(X\cdot Z)1\otimes Y-(X\cdot Y)1\otimes Z-(Y\cdot X)1\otimes Z
+(Z\cdot X)1\otimes Y=0.$$
\end{proof}

Now we prove (\ref{mu-zero-1}). We denote $\ell^{\theta}:=\log \theta\colon \pi \to \widehat{\mathcal{L}}$.
Note that the logarithm $\log\colon 1+\widehat{T}_1 \to \widehat{T}_1, u\mapsto
\sum_{n=1}^{\infty}((-1)^{n-1}/n)(u-1)^n$ gives a bijection from the set of group-like elements of $\widehat{T}$
to $\widehat{\mathcal{L}}$. For $x\in \pi$ let $\ell^{\theta}_k(x)$ be the degree $k$ part of $\ell^{\theta}(x)\in \widehat{\mathcal{L}}$.
We remark that $\ell^{\theta}_1(x)=[x]\in H$.

{\it Step 1.} Let $\theta^0$ be a symplectic expansion satisfying
$$\ell^{\theta^0}_2(\alpha_i)=\frac{1}{2}[A_i,B_i], \quad \ell^{\theta^0}_2(\beta_i)=-\frac{1}{2}[A_i,B_i].$$
Such a symplectic expansion does exist. For example, see Massuyeau \cite{Mas} and Kuno \cite{Ku1}.
Since $\theta^0=\exp(\ell^{\theta^0})$,
\begin{eqnarray*}
\theta^0_3(\alpha_i) &=&
\ell^{\theta^0}_3(\alpha_i)+\frac{1}{2}(A_i \ell^{\theta^0}_2(\alpha_i)+\ell^{\theta^0}_2(\alpha_i)A_i)+\frac{1}{6}A_iA_iA_i \\
&=& \ell^{\theta^0}_3(\alpha_i)+\frac{1}{4}(A_i[A_i,B_i]+[A_i,B_i]A_i)+\frac{1}{6}A_iA_iA_i.
\end{eqnarray*}
By Lemma \ref{LH3}, $\mu^{\rm alg}(\ell^{\theta^0}_3(\alpha_i))=0$ and clearly $\mu^{\rm alg}(A_iA_iA_i)=0$.
Therefore,
\begin{eqnarray*}
\mu^{\rm alg}(\theta^0_3(\alpha_i)) &=& \frac{1}{4}\mu^{\rm alg}(A_i[A_i,B_i]+[A_i,B_i]A_i) \\
&=& \frac{1}{4}(1\otimes A_i+1\otimes A_i)=\frac{1}{2}(1\otimes A_i).
\end{eqnarray*}
Similarly, we obtain $\mu^{\rm alg}(\theta^0_3(\beta_i))=(-1/2)(1\otimes B_i)$. Substituting these equations into
(\ref{mu-A-B}), we have $\mu^{\theta^0}_{(0)}(A_i)=(-1/2)(1\otimes A_i)$ and $\mu^{\theta^0}_{(0)}(B_i)=(-1/2)(1\otimes B_i)$.
This completes the proof of (\ref{mu-zero-1}) for $\theta^0$.

{\it Step 2.} Let $\theta$ be an arbitrarily symplectic expansion. Then there exist
$u_1\in \Lambda^3 H$ and $u_2\in {\rm Hom}(H,\widehat{\mathcal{L}}\cap H^{\otimes 3})$ such that
\begin{equation}
\label{th-th}
\theta_3(x)=\theta^0_3(x)+(u_1\otimes 1+1\otimes u_1)\theta^0_2(x)+u_2([x])
\end{equation}
for any $x\in \pi$. This is proved by a similar way to the proof of Lemma 6.4.2 in \cite{KK1}.
Here $\Lambda^3 H$ is the third exterior power of $H$ and
we regard $u_1$ as an element of ${\rm Hom}(H,H^{\otimes 2})$ by the inclusion
$\Lambda^3 H \subset H^{\otimes 3}\cong H^* \otimes H^{\otimes 2}={\rm Hom}(H,H^{\otimes 2})$,
$X\wedge Y\wedge Z \mapsto XYZ-XZY-YXZ+YZX+ZXY-ZYX$. Notice that $u_1(H)\subset \Lambda^2 H$.

\begin{lem}
\label{muXX}
For any $u_1\in \Lambda^3H$ and $X\in H$ we have
$$\mu^{\rm alg}((u_1\otimes 1+1\otimes u_1)XX)=0.$$
\end{lem}

\begin{proof}
It suffices to consider the case $X\neq 0$. There exists a $\mathbb{Q}$-symplectic
basis $\{ A_j,B_j\}_j \subset H$ such that $X=A_1$. By linearity, we may
assume that $u_1=Y\wedge Z\wedge W$ where $Y,Z,W\in \{ A_j,B_j\}_j$. Now the assertion is proved
by a direct computation.
\end{proof}

Let $x\in \pi$. Since $u_1(H)\subset \Lambda^2 H$, we have $(u_1\otimes 1+1\otimes u_1)\ell^{\theta}_2(x)
\in \widehat{\mathcal{L}}\cap H^{\otimes 3}$. By Lemma \ref{LH3}, $\mu^{\rm alg}((u_1\otimes 1+1\otimes u_1)\ell_2^{\theta}(x))=0$.
By $\theta^0_2(x)=\ell^{\theta^0}_2(x)+(1/2)[x][x]$ and Lemma \ref{muXX}, we conclude
$$\mu^{\rm alg}((u_1\otimes 1+1\otimes u_1)\theta_2^0(x))=0$$
for any $x\in \pi$. Also, we have $\mu^{\rm alg}(u_2([x]))=0$ by Lemma \ref{LH3}. Therefore by (\ref{th-th})
$\mu^{\rm alg}(\theta_3(x))=\mu^{\rm alg}(\theta^0_3(x))$ for any $x\in \pi$.
Now from Step 1 and (\ref{mu-A-B}) we have
$\mu^{\theta}_{(0)}(A_i)=(-1/2)(1\otimes A_i)$ and $\mu^{\theta}_{(0)}(B_i)=(-1/2)(1\otimes B_i)$.
This completes the proof of (\ref{mu-zero-1}). \qed

Next we compute $\mu^{\theta}_{(0)}(X_1\cdots X_m)$ for $X_1,\ldots,X_m\in H$, $m\ge 2$.
Since the constant term of $s(z)$ is $-1/2$, by (\ref{53mu}) we have
\begin{eqnarray}
&& \mu^{\theta}_{(0)}(X_1\cdots X_m) \nonumber \\
&=& -\frac{1}{2}(1\otimes N)\sum_{1\le i<j\le m}(X_1\cdots X_{i-1}\otimes 1)((1\otimes \iota)
\Delta(X_iX_j))(X_{j+1}\cdots X_m\otimes X_{i+1}\cdots X_{j-1}) \nonumber \\
&& +\sum_{i=1}^m (X_1\cdots X_{i-1} \otimes 1)\mu^{\theta}_{(0)}(X_i)(X_{i+1}\cdots X_m\otimes 1).
\label{mu-zero-m}
\end{eqnarray}
Note that we have
\begin{eqnarray}
(1\otimes \iota)(\Delta(X_iX_j))&=& (1\otimes \iota)(X_iX_j\otimes 1+X_i\otimes X_j+X_j\otimes X_i+1\otimes X_iX_j) \nonumber \\
&=& X_iX_j\otimes 1-X_i\otimes X_j-X_j\otimes X_i+1\otimes X_jX_i. \label{DXX}
\end{eqnarray}
Substituting (\ref{mu-zero-1}), (\ref{DXX}) into (\ref{mu-zero-m}), and computing directly, we obtain
$$\mu^{\theta}_{(0)}(X_1\cdots X_m)=-\frac{1}{2}(1\otimes N(X_1\cdots X_m)).$$
This completes the proof of Theorem \ref{mu-zero}. \qed

\appendix

\section{Lie bialgebras and their bimodules}
For the sake of the reader we collect the definitions
of a Lie bialgebra and its bimodules.

We work over the rationals $\mathbb{Q}$.
Let $V$ be a $\mathbb{Q}$-vector space and let
$T=T_V \colon V^{\otimes 2}\to V^{\otimes 2}$ and
$N=N_V \colon V^{\otimes 3} \to V^{\otimes 3}$ be the linear maps
defined by $T(X\otimes Y)=Y\otimes X$ and $N(X\otimes Y\otimes Z)
=X\otimes Y\otimes Z+Y\otimes Z\otimes X+Z\otimes X\otimes Y$
for $X,Y,Z\in V$.

\subsection{Lie bialgebras}\label{LA}
Let $\gfg$ be a $\mathbb{Q}$-vector space equipped with
$\mathbb{Q}$-linear maps $\nabla\colon \gfg\otimes\gfg \to \gfg$ and
$\delta\colon \gfg\to \gfg\otimes\gfg$. Recall that $\gfg$ is called a {\it Lie bialgebra}
with respect to $\nabla$ and $\delta$, if

\begin{enumerate}
\item the pair $(\gfg,\nabla)$ is a Lie algebra,
i.e., $\nabla$ satisfies {\it the skew condition} and {\it the Jacobi identity}
$$\nabla T = -\nabla\colon \gfg^{\otimes 2} \to \gfg, \quad
\nabla(\nabla\otimes 1)N = 0\colon \gfg^{\otimes 3} \to \gfg,$$

\item the pair $(\gfg,\delta)$ is a Lie coalgebra, i.e., $\delta$ satisfies
{\it the coskew condition} and {\it the coJacobi identity}
$$T\delta = -\delta\colon \gfg\to \gfg^{\otimes 2}, \quad
N(\delta\otimes 1)\delta = 0\colon \gfg \to \gfg^{\otimes 3},$$

\item  the maps $\nabla$ and $\delta$
satisfy {\it the compatibility}
$$\forall X, \forall Y \in \gfg, \quad
\delta [X, Y] = \sigma(X)(\delta Y) - \sigma(Y)(\delta X).$$
\end{enumerate}

Here we denote $[X, Y] := \nabla(X\otimes Y)$
and $\sigma(X)(Y\otimes Z)=[X,Y]\otimes Z+Y\otimes [X,Z]$ for $X,Y,Z\in \gfg$.
The map $\nabla$ is called the {\it bracket}, and the map $\delta$ is
called the {\it cobracket}.

Moreover if {\it the involutivity}
$$\nabla\delta = 0\colon \gfg \to \gfg$$
holds, we say $\gfg$ is {\it involutive}.

\subsection{Lie comodules and bimodules}

Let $\gfg$ be a Lie algebra.
Recall that a {\it left $\gfg$-module} is a pair $(M,\sigma)$
where $M$ is a $\mathbb{Q}$-vector space and $\sigma$ is
a $\mathbb{Q}$-linear map $\sigma\colon \gfg\otimes M \to M$,
$X\otimes m\mapsto Xm$, satisfying
$$\forall X, \forall Y \in \gfg, \forall m \in M, \quad 
[X, Y]m = X(Ym) - Y(Xm).$$
This condition is equivalent to the commutativity of the diagram
$$
\begin{CD}
\gfg \otimes \gfg \otimes M @>{((1-T)\otimes 1_M)(1_{\gfg}\otimes \sigma})
>> \gfg\otimes M\\
@V{\nabla \otimes 1_M}VV @V{\sigma}VV\\
\gfg \otimes M @>{\sigma}>> M.
\end{CD}
$$
If we define $\overline{\sigma}\colon M\otimes \gfg \to M$ by 
$\overline{\sigma}(m\otimes X):=-\sigma(X\otimes m)=-Xm$ for $m\in M$ and $X\in \gfg$,
then the pair $(M,\overline{\sigma})$
is a {\it right $\gfg$-module}, i.e., the following diagram commutes:
$$
\begin{CD}
M \otimes \gfg \otimes \gfg @>{(1_M \otimes (1-T))(\overline{\sigma}}\otimes 1_{\gfg})
>> M\otimes \gfg \\
@V{1_M \otimes \nabla}VV @V{\overline{\sigma}}VV\\
M\otimes \gfg @>{\overline{\sigma}}>> M.
\end{CD}
$$

Next let $(\gfg, \delta)$ be a Lie coalgebra and $M$ a
$\mathbb{Q}$-vector space equipped with a $\mathbb{Q}$-linear map 
$\mu\colon M \to M\otimes\gfg$. We say the pair $(M,\mu)$ is a
{\it right $\gfg$-comudule} if the following diagram commutes:
\begin{equation}
\begin{CD}
M @>{\mu}>> M\otimes\gfg\\
@V{\mu}VV @V{1_M\otimes\delta}VV\\
M\otimes\gfg @>{(1_M\otimes(1-T))(\mu\otimes
1_\gfg)}>> M\otimes\gfg\otimes\gfg.
\end{CD}
\label{A6coaction}
\end{equation}
Similarly, we say a pair $(M,\overline{\mu})$ is a {\it left $\gfg$-comodule}
if $\overline{\mu}$ is a $\mathbb{Q}$-linear map
$\overline{\mu}\colon M \to \gfg\otimes M$ and the following diagram commutes:
\begin{equation*}
\begin{CD}
M @>{\overline{\mu}}>> \gfg\otimes M\\
@V{\overline{\mu}}VV @V{\delta\otimes 1_M}VV\\
\gfg\otimes M @>{((1-T)\otimes 1_M)(1_\gfg\otimes\overline{\mu}
)}>> \gfg\otimes\gfg\otimes M.
\end{CD}
\end{equation*}
If we denote the switch map by $T_{\gfg,M}\colon \gfg \otimes M \to M\otimes
\gfg$, $X\otimes m \mapsto m\otimes X$, then it is easy to see that
$(M, \overline{\mu})$ is a left $\gfg$-comodule if and only if 
$(M, -T_{\gfg,M}\overline{\mu})$ is a right $\gfg$-comodule. 

Finally let $\gfg$ be a Lie bialgebra with $\nabla$ the bracket and $\delta$
the cobracket, $(M, \sigma)$ a left $\gfg$-module,
and $(M, \mu)$ a right $\gfg$-comodule
with the same underlying vector space $M$. 
We define $\overline{\sigma}\colon M \otimes\gfg \to M$ by
$\overline{\sigma}(m\otimes X) :=-Xm$, as before. Then $(M, \overline{\sigma})$ is a
right $\gfg$-module. We say the triple $(M,\overline{\sigma},\mu)$ is
a {\it right $\gfg$-bimodule} if $\sigma$ and $\mu$ satisfy
{\it the compatibility}
\begin{equation}
\forall m \in M, \forall Y \in \gfg, \quad
\sigma(Y)\mu (m) - \mu(Ym) -(\overline{\sigma}\otimes
1_\gfg)(1_M\otimes\delta)(m\otimes Y) = 0.
\label{A6compatible+}
\end{equation}
Here $\sigma(Y)\mu(m)=(\sigma \otimes 1_M)(Y\otimes \mu(m))
+(1_M \otimes {\rm ad}(Y))\mu(m)$ and ${\rm ad}(Y)(Z)=[Y,Z]$
for $Z\in \gfg$. Then we also call the triple $(M, \sigma, \overline{\mu})$ 
defined by $\overline{\mu} := -T_{M,\gfg}\mu\colon M \to \gfg\otimes M$ 
a {\it left $\gfg$-bimodule}. Moreover, if $\gfg$ is involutive and
the condition
\begin{equation}
\overline{\sigma}\mu = 0\colon M \to M
\label{A6involutive}
\end{equation}
holds, we say $M$ is {\it involutive}.

\noindent \textsc{Nariya Kawazumi\\
Department of Mathematical Sciences,\\
University of Tokyo,\\
3-8-1 Komaba Meguro-ku Tokyo 153-8914 JAPAN}\\
\noindent \texttt{E-mail address: kawazumi@ms.u-tokyo.ac.jp}

\vspace{0.5cm}

\noindent \textsc{Yusuke Kuno\\
Department of Mathematics,\\
Tsuda College,\\
2-1-1 Tsuda-Machi, Kodaira-shi, Tokyo 187-8577 JAPAN}\\
\noindent \texttt{E-mail address: kunotti@tsuda.ac.jp}

\end{document}